\documentclass[a4paper,11pt]{article}
\usepackage{amsmath,amscd,amssymb}
\usepackage{pstricks}
\usepackage{pst-plot}

\def \vs {\vskip}
\def \hs {\hskip}
\def \noi {\noindent}

\def \oo {{\cal O}}

\def \L {{\cal L}}
\def \M {{\cal M}}

\def \a {{\alpha}}
\def \b {{\beta}}
\def \ga {{\gamma}}

\def \pd {\varpi}

\def \Del {\Delta}

\def \p {{\mathbb P}}
\def \k {{\mathbb K}}

\def \G {{\mathbb G}}
\def \N {{\mathbb{N}}}

\def \Z {{\mathbb{Z}}}

\def \supp {{\rm Supp}}
\def \pic {{\rm Pic}}
\def \weil {{\rm Weil}}
\def \stab {{\rm Stab}}

\def \ne {N\!E}
\def \card {{\rm Card}}
\def \codim {{\rm Codim}}

\def \pu {\p^1}

\def \Mor #1#2{{\bf{Hom}}_{#1}(\pu,#2)}

\def \wt {{\widetilde{w}}}

\def \Xt {{\widetilde{X}}}

\def \pit {{\widetilde{\pi}}}

\def \at {{\widetilde{\a}}}

\def \Ch {{\widehat{C}}}
\def \Qh {{\widehat{Q}}}
\def \Xh {{\widehat{X}}}
\def \Yh {{\widehat{Y}}}

\def \wh {{\widehat{w}}}
\def \vh {{\widehat{v}}}
\def \pih {{\widehat{\pi}}}

\def \comp {{\mathfrak{ne}}}

\def \tha #1#2{\noi{\bf#1{\uppercase{\footnotesize{#2}}}}}

\newtheorem{theor}{\tha{T}{heorem}}[section]
\newenvironment{theo}{
  \begin{theor}\hs -0.2 cm {\bf .} ---  }
{  \end{theor}}

\newtheorem{propo}[theor]{\tha{P}{roposition}}
\newenvironment{prop}{
  \begin{propo}\hs -0.2 cm {\bf .} ---  }
{  \end{propo}}

\newtheorem{lemma}[theor]{\tha{L}{emma}}
\newenvironment{lemm}{
  \begin{lemma}\hs -0.2 cm {\bf .} ---  }
{  \end{lemma}}

\newtheorem{fait}[theor]{\tha{F}{act}}
\newenvironment{fact}{
  \begin{fait}\hs -0.2 cm {\bf .} ---  }
{  \end{fait}}

\newtheorem{defini}[theor]{\tha{D}{efinition}}
\newenvironment{defi}{
  \begin{defini}\hs -0.2 cm {\bf .} ---  }
{  \end{defini}}

\newtheorem{corollaire}[theor]{\tha{C}{orollary}}
\newenvironment{coro}{
  \begin{corollaire}\hs -0.2 cm {\bf .} ---  }
{  \end{corollaire}}

\newtheorem{exemple}[theor]{\sc{Example}}
\newenvironment{exem}{
  \begin{exemple}\hs -0.2 cm {\bf .} ---  }
{  \end{exemple}}

\newtheorem{remarq}[theor]{\sc{Remark}}
\newenvironment{rema}{
  \begin{remarq}\hs -0.2 cm {\bf .} ---  }
{  \end{remarq}}

\newtheorem{construc}{\sc{Construction}}
\newenvironment{cons}{
  \begin{construc}\hs -0.2 cm {\bf .} ---  }
{  \end{construc}}

\newenvironment{preu}{\noi\sc{Proof}\hs 0.1 cm 
--- \rm }{\hfill$\Box$\vs 0.2 cm}

\def \scal #1#2{\langle #1,#2 \rangle}
\def \sca #1#2{\langle #1,#2 \rangle}

\begin{document}
\title{{\bf Small resolutions of minuscule Schubert varieties}}
\author{\Large{\sc{Nicolas Perrin}}}
\date{}

\maketitle

\vs 1 cm

\begin{abstract}
Let $X$ be a minuscule Schubert variety. In this article, we use the
combinatorics of quivers to define new quasi-resolutions of $X$. We
describe in particular all relative minimal models $\pih:\Xh\to X$ of
$X$ and prove that all the morphisms $\pih$ are small (in the sense of
intersection cohomology). In particular, all small resolutions of $X$
are given by the smooth relative minimal models $\Xh$ and we describe
all of them. As another application of this decription of relative
minimal models, we give a more intrinsic statement of the main result
of \cite{Perrin3}.
\end{abstract}

\vs 2 cm

\centerline{\large{\textsc{Introduction}}}

\vs 1 cm

Schubert varieties have been intensively studied and are of great
importance in representation theory. There are several ways to
understand the geometry and the singularities of Schubert varieties. One
way is to describe the singular locus, its irreducible components and
the generic singularity in each of these components. This has been
completely achieved for $GL_n$ only recently (see \cite{Manivel1},
\cite{Manivel2}, \cite{BilleyW}, \cite{KasselLR} and
\cite{Cortez2}). For the general case, there are only partial results
(for an account, see 
\cite{LakshBill}). This description of singularities enables in
particular to calculate the Kazhdan-Lusztig polynomials. For this
study, combinatoric tools are usefull but it is also usefull to
construct special resolutions of Schubert varities (see for example
\cite{Cortez1} and \cite{Cortez2}). Another way to study Schubert
varieties is to calculate the cohomology of line bundles and
especially to prove some vanishing theorems. This has been done in many
ways, one of them thanks to the Bott-Samelson resolutions
(\cite{Kempf3}, \cite{MethaRama}, \cite{RamaRama} or \cite{Rama1}). In this
article, we want to study the geometry of Schubert varieties thanks to
the study of some particular resolutions of this variety. 

A nice resolution of Schubert varieties is the Bott-Samelson
resolution (see \cite{Demazure1} or \cite{Hansen}). This resolution is
usefull to study the singularities and standard monomial theory of
Schubert varities (see \cite{Litt-et-ali}) and also to study the
geometry of Schubert varieties (for curves on minuscule Schubert
varieties, see for example \cite{Perrin3} and
\cite{Perrin-ell}). However these resolutions are big in the sense that
the fibers have big dimensions and there are many
contracted subvarieties. Another class of resolutions is of particular
importance, the $IH$-small resolution (that is to say the small
resolution in the sense of Intersection Cohomology, see definition
\ref{small-def}). These are well suited for the
calculation of Kazhdan-Lusztig polynomials. In particular
A. Zelevinsky in \cite{Zelevinsky} constructed some $IH$-small
resolutions for grassmannian Schubert varieties and gave a geometric
interpretation of the combinatoric computation of Kazhdan-Lusztig
polynomials by A. Lascoux and M.-P. Sch{\"u}tzenberger in
\cite{Lasc}. Later P. Sankaran and
P. Vanchinathan \cite{SanVan} and \cite{SanVan2} constructed small
resolution of some minuscule and cominuscule Schubert varities for
$SO_{2n}$ and $Sp_{2n}$ and 
calculated the corresponding Kazhdan-Lusztig polynomials. In their
article \cite{SanVan}, they construct some minuscule Schubert
varieties not admiting small
resolutions. These examples are locally factorial Schubert variety
(more precisely with singularities in codimension 2) for which the
theorem of purity (see \cite{EGA} theorem 21.12.12) says that there is
no $IH$-small resolution. 

In this article we study the $IH$-small resolutions of minuscule
Schubert varities. We will generalise the constructions of
A. Zelevinsky and P. Sankaran and P. Vanchinathan to any minuscule
Schubert variety and describe all $IH$-small resolutions. We do not
adress the problem of calculating Kazhdan-Lusztig polynomials. This
has been done in a combinatoric way for all minuscule Schubert
varities by B. D. Boe in \cite{Boe} and we hope that our construction
will lead to a geometric interpretation of these results.

In order to define our resolutions, we introduce a combinatorial
objet: a quiver associated to a minuscule Schubert variety
$X(w)$. This quiver is defined thanks to a reduced writing $\wt$ of
$w$ which is unique for minuscule Schubert variety. These quivers seem
to be the same as the quivers defined by S. Zelikson in \cite{Zelikson}
but we did not check this. To this quiver, we associate a configuration
variety $\Xt(\wt)$ which is simply the Bott-Samelson resolution. The
fact that the Bott-Samelson resolution can be seen as a confuguration
variety was already known by P. Magyar \cite{Magyar-conf1} but the use of
the quiver in the situation is very usefull. In the minuscule case,
the quiver will have a
very special and rigid geometry and in particular we will define the
pics of the quiver and the height of a pic. As in the case of
A. Zelevinsky's
construction \cite{Zelevinsky}, the choice of an order on the pics
will lead to a partial resolution $\Xh(\wh)\to X(w)$ of the Schubert
variety. However, the variety $\Xh(\wh)$ will in general be locally
factorial but not smooth.

These varieties are however interesting for the relative minimal model
program. We study the relative minimal models of $X(w)$ and prove the
following theorem:

\begin{theo}
  The relative minimal models of $X(w)$ are the varieties $\Xh(\wh)$
  obtained thanks to an order on the pics preserving the order on the
  heights of the pics (see construction \ref{lespetites} for more
  details).
\end{theo}

For grassmannian Schubert varieties, there are always $IH$-small
resolution \cite{Zelevinsky}. As P. Sankaran and P. Vanchinathan
proved this is not true for a general minuscule Schubert variety. The
following theorem proves that we need to replace $IH$-small
resolutions by relative minimal models to generalise the result of
A. Zelevinsky to any minuscule Schubert variety (that is so say allow
locally factorial singularities):

\begin{theo}
  The morphism $\pih:\Xh(\wh)\to X(w)$ from a relative minimal model
  to a minuscule Schubert variety $X(w)$ is small.
\end{theo}

In other words, the relative minimal models play the role for general
minuscule Schubert varieties of $IH$-small resolution for grassmannian
Schubert varieties. However they do not share the same
nice properties and in particular are not as well fitted as $IH$-small
resolutions for the computation of Kazhdan-Lusztig polynomials. We
also describe all relative canonical models of $X(w)$.
Furthermore, the following result of B. Totaro \cite{Totaro} using a
key result of J. Wisniewski \cite{wisniewski} tells us
to look for $IH$-small resolutions in the class of relative minimal
models:

\begin{theo}[Totaro-Wisniewski]
  Any $IH$-small resolution of a normal variety $X$ is a 
  relative minimal model for $X$.
\end{theo}

In particular, in our situation, all $IH$-small resolutions of $X(w)$
are given by the smooth relative minimal models. We then give a
combinatorial criterion on the quiver for the relative minimal model
$\Xh(\wh)$ to be smooth. We get the following:

\begin{theo}
  The variety $\Xh(\wh)$ is an $IH$-small resolution if and only if
  the order on the
  pics preserves the order on the heights of the pics and if at step
  $i$ the pic $p_i$ is minuscule for the quiver (see definition
  \ref{pic-min}).
\end{theo}

In particular, we are able to say which minuscule Schubert variety
admits an $IH$-small resolution. At the end of the article we sketch
another way the prove this: an $IH$-small resolution has to factor
through the relative canonical minimal model $X_{can(w)}$ of $X(w)$ and is a
  crepant resolution of $X_{can}(w)$. In particular, the stringy Euler
    number $e_{\rm st}$ defined by V. Batyrev \cite{batyrev} of
    $X_{can}(w)$ has to be an integer. We give a formula for $e_{\rm
      st}(X_{can}(w))$ and it is an easy verification that we recover
    in this way all minuscule Schubert varieties not admitting an
    $IH$-small resolution.

Another motivation for the study of resolutions and partial
resolutions of Schubert varieties (and in fact our motivation at the
begining of the study) is the following reinterpretation of our result
in \cite{Perrin3}. Let $X$ be a minuscule Schubert variety and
$\pi:\Xh\to X$ any relative minimal model. For a 1-cycle class $\a\in
A_1(X)$ define the set 
$$\comp(\a)=\{\b\in \ne(\Xh)\ /\ \pi_*\b=\a\}$$
where $\ne(\Xh)$ is the cone of effective 1-cycles in $\Xh$.

\begin{theo}
  The irreducible components of $\Mor{\a}{X}$ of the scheme of
  morphisms from $\pu$ to $X$ of class $\a$ are indexed by
  $\comp(\a)$. 
\end{theo}

The same kind of results are true for other special Schubert
varieties (for example for cominuscule Schubert varieties, this will
be studied in \cite{moivirtuel}). It is also the case for cones over
homogeneous varieties (see \cite{Perrin-cone}).

Let us give an overview of the article. In paragraphs 1 and 2 we
recall some basic notations, definitions and results on elements of
the Weyl group and minuscule Schubert varieties. In the third
paragraph, we define the quiver associated to a reduced writing $\wt$
and the corresponding configuration variety. We study basic properties
of these varieties (Weil, Cartier and canonical divisors, 1-cycles and
intersection formulae), and link them with the geometry of the
quiver. In the fourth paragraph, we describe the particular geometry
of a quiver asociated to a minuscule Schubert variety and study the
link with the geometry of the Schubert variety (Weil, Cartier and
canonical divisors, 1-cycles and intersection formulae). In the fifth
paragraph, we construct and study a generalisation of Bott-Samelson
resolution which is between the Schubert variety and the Bott-Samelson
resolution. In the sixth paragraph, we describe the geometry of this
generalisation (one more time Weil, Cartier and canonical divisors,
1-cycles and intersection formulae) and describe all the relative
Mori theory for a minuscule Schubert variety. In the last paragraph,
we prove that all relative minimal models are $IH$-small and describe
all $IH$-small resolutions of minuscule Schubert varieties. In this
paragraph, in constrast with the rest of the paper, we use a case by
case analysis. A general proof could be possible but we think it would
be more complicated and would lead to much more combinatorics. In an
appendix, we describe the quivers of minuscule Schubert variety. We
use this description intensively in the last paragraph.

\newpage

\begin{small}
\tableofcontents
\end{small}

\newpage

\section{Notations}

Let $G$ be a semi-simple algebraic group, fix $T$ a maximal torus and
$B$ a Borel subgroup containing $T$. Let us denote by $\Del$ the set
of all roots, by $\Del^+$ (resp. $\Del^-$) the set of positive
(resp. negative) roots, by $S$ the set of simple roots associated to
the data $(G,T,B)$ and by $W$ the associated Weyl group. If $P$ is a
parabolic subgroup containing $B$ we note $W_P$ the subgroup of $W$
corresponding to $P$. We will also denote by $\Sigma(P)$ the set of
simple roots $\b$ such that $U_{-\b}\not\subset P$.

\begin{defi}
Let $w\in W$, let us denote by $P^w$ the largest parabolic
subgroup of $G$ such that the morphism
$\overline{BwB}/B\to\overline{BwP^w}/P^w$ is a $P^w/B$
fibration. 

Let us also denote by $P_w$ the stabiliser of
$X(w)=\overline{BwP^w}/P^w$ in $G/P^w$.
\end{defi}

\begin{defi}
(\i) Let $w\in W$, we define the support of $W$ denoted by $\supp(w)$ to
  be the set of simple roots $\b$ such that $s_\b$ appears in a
  reduced writing of $w$. This set is indepedent of the reduced
  writing and only depends on $w$.

(\i\i) We will denote by $G_w$ the smallest reductive subgroup of $G$
containing all the groups $U_\b$ for $\b\in\supp(w)$. It is easy to
see that there is an isomorphism 
$$X(w)=\overline{P_wwP^w}/P^w\simeq \overline{(P_w\cap G_w)w(P^w\cap
  G_w)}/(P^w\cap G_w).$$

(\i\i\i) Let us also define the boundary of $G_w$ denoted by
$\partial(G_w)$ to be the set of simple roots $\b$ not contained in
$\supp(w)$ and non commuting with $w$.
\end{defi}

\section{Minuscule Schubert varieties}
\label{minuscule}

In this paragraph we recall the notion
of minuscule weight and study the related homogeneous
and Schubert varieties. Our basic reference will be \cite{GofG/P3}.

\begin{defi}
Let $\pd$ be a fundamental weight, 

(\i) we say that $\pd$ is minuscule if we have
$\sca{\a^\vee}{\pd}\leq1$ for all positive root $\a\in\Del^+$;

(\i\i) we say that $\pd$ is cominuscule if $\sca{\a_0^\vee}{\pd}=1$
where $\a_0$ is the longest root.
\end{defi}

With the notation of N. Bourbaki \cite{bourb}, the minuscule and
cominuscule weights are:

\begin{center}
\begin{tabular}{|c|c|c|}
\hline
Type&minuscule&cominuscule\\
\hline
$A_n$&$\pd_1\cdots\pd_n$&same weights\\
\hline
$B_n$&$\pd_n$&$\pd_1$\\
\hline
$C_n$&$\pd_1$&$\pd_n$\\
\hline
$D_n$&$\pd_1$, $\pd_{n-1}$ and $\pd_n$&same weights\\
\hline
$E_6$&$\pd_1$ and $\pd_6$&same weights\\
\hline
$E_7$&$\pd_7$&same weight\\
\hline
$E_8$&none&none\\
\hline
$F_4$&none&none\\
\hline
$G_2$&none&none\\
\hline
\end{tabular}
\end{center}

\begin{defi}
  Let $\pd$ be a  minuscule weight and let $P_\pd$
be the associated parabolic subgroup. The homogeneous variety
$G/P_\pd$ is then said to be minuscule. The Schubert varieties of a
minuscule homogeneous variety are called minuscule Schubert
varieties.

An element $w\in W$ is said to be minuscule if $G/P^w$ is a
minuscule homogeneous variety.
\end{defi}

\begin{rema}
  To study minuscule homogeneous varieties and their
  Schubert varieties, it is sufficent to restrict ourselves to
  simply-laced groups.

In fact the variety $G/P_{\pd_n}$ with $G={\rm Spin}_{2n+1}$ is
isomorphic to the variety $G'/P'_{\pd_{n+1}}$ with
$G'={\rm Spin}_{2n+2}$ and there is a one to one correspondence
between Schubert varieties thanks to this isomorphism. The same
situation occurs with $G/P_{\pd_1}$, $G={\rm Sp}_{2n}$ and
$G'/P'_{\pd_1}$, $G'={\rm SL}_{2n}$.
\end{rema}

\section{Quivers and configuration varieties}
\label{BS}

\subsection{Quiver associated to a reduced writing}

Let $P_\pd$ be a parabolic subgroup of $G$ associated to a minuscule
weight $\pd$, let us consider an element $\bar w\in W/W_{P_\pd}$ and
let $w$ be the shortest element in the class $\bar w$. To any reduced
writing 
\begin{equation}
\label{ecri-red}
  w=s_{\b_1}\cdots s_{\b_r}
\end{equation}
of $w$ in terms of simple reflections (for all $i\in[1,r]$, we have
$\b_i\in S$) we associate a quiver with colored vertices. Let us first
give the following

\begin{defi}
\label{quiver}
  For a fixed reduced writing (\ref{ecri-red}) of $w$, we define the
  successor $s(i)$ (resp. the predecessor $p(i)$) of an element
  $i\in[1,r]$ by $\displaystyle{s(i)=\min\{j\in[1,r]\ /\ j>i\
  \textrm{\textit{and}}\ \b_j=\b_i\}}$ resp. by
  $\displaystyle{p(i)=\max\{j\in[1,r]\ /\ j<i\ \textrm{\textit{and}}\
  \b_j=\b_i\}}$.
\end{defi}

Now we can define the quiver $Q_{\wt}$ associated to the reduced
writing (\ref{ecri-red}) of $w$ (we will denote by $\wt$ the data of
$w$ with such a reduced writing):

\begin{defi}
  Let us denote by $Q_\wt$ the quiver whose set of vertices is in
  bijection with the set $[1,r]$ and whose arrows are given in the
  following way: there is an arrow from $i$ to $j$ if
  $\scal{\b_j^\vee}{\b_i}\neq0$ and $i<j<s(i)$.

This quiver comes with a coloration of its edges by simple roots
thanks to the application $\b:[1,r]\to S$ such that $\b(i)=\b_i$.
\end{defi}

\begin{rema}
  It is equivalent to give the reduced writing $\wt$ or the quiver
  $Q_\wt$. This quiver seems to be the same as the one defined by
  S. Zelikson \cite{Zelikson} for ADE types.
\end{rema}

\subsection{Configuration varieties and Bott-Samelson resolution}

In this paragraph, we associate to any quiver $Q_\wt$, a configuration
variety $\Xt(\wt)$ (see also \cite{Magyar-conf1}). 
Let $x$ be an element in $G/B$ and $\b_i$ a
simple root. Let us denote by  $B_{\b_i}$ the parabolic subgroup
generated by $B$ and $U_{-\b_i}$. We have a projection morphism
$\pi_{\b_i}:G/B\to G/B_{\b_i}$ whose fibers are isomorphic to $\pu$ and
we denote by $\p(x,\b_i)$ the projective line
$\pi_{\b_i}^{-1}(\pi_{\b_i}(x))$.

\begin{defi}
  Let $Q_\wt$ a quiver associated to a writing $\wt$ of $w$. We define
  the configuration variety $\Xt(\wt)$ by
$$\Xt(\wt)=\left\{(x_1,\cdots,x_r)\in\prod_{i=1}^rG/B\ /\ x_0=1\ {\rm
    and}\ x_i\in\p(x_{i-1},\b_i)\ \textrm{for all}\
    i\in[1,r]\right\}.$$
\end{defi}

\begin{rema}
(\i) If we denote by $P_{\b_i}$ the maximal parabolic not containing
  $U_{-\b_i}$, the restriction of the morphism $G/B\to G/P_{\b_i}$ to
  $\p(x,\b_i)$ is an isomorphism so that $\Xt(\wt)$ is isomorphic to
$$\left\{(x_1,\cdots,x_r)\in\prod_{i=1}^rG/P_{\b_i}\ /\ x_0=1\ {\rm
    and}\ x_i\in\p(x_{i-1},\b_i)\ \textrm{for all}\
    i\in[1,r]\right\}.$$

(\i\i) In his article \cite{Magyar-conf1}, P. Magyar shows that the variety
$\Xt(\wt)$ is isomorphic to the classical Bott-Samelson variety
described for example in \cite{Demazure1}.
\end{rema}

Because the element $w$ is the shortest in the class $\bar w$
and because the writing $\wt$ is reduced, the last root $\b_r$ has to
be the only simple root $\b$ such that $\scal{\pd^\vee}{\b}=1$ and
$P_r=P_\pd$. The image of the projection morphism $\pi:\Xt(\wt)\to G/P_r$ is
the minuscule Schubert variety $X(\bar w)$. The morphism
$\pi:\Xt(\wt)\to X(\bar w)$ is birational. 

\subsection{Cycles on the configuration variety}

\subsubsection{A basis of the Chow ring}
\label{base-chow}
\label{canonique}
\label{canonique-coord}

In this paragraph, we describe some particular elements in the Chow
ring $A_*(\Xt(\wt))$. We describe a basis of this ring a the cones of
amples divisors and effective curves. We calculate the canonical
divisor in terms of these basis.

For all $k\in[1,r]$, let us denote by $X_k$ the image of $\Xt(\wt)$ in
the product $\displaystyle{\prod_{i=1}^kG/P_{\b_i}}$ and
$X_0=\{x_0\}$. We have natural projection morphism $f_k:X_k\to
X_{k-1}$ for all $k\in[1,r]$ which are $\pu$-fibration. The morphism
$\sigma_k:X_{k-1}\to X_k$ defined by
$\sigma_{k-1}(x_1,\cdots,x_{k-1})=(x_1,\cdots,x_{k-1},x_{p(k)})$ with
$x_{p(k)}=1$ if $p(k)$ does not exist is a section of $f_k$. We
recover in this way the structure of $\Xt(\wt)$ as a tower of
$\pu$-fibrations with sections described in \cite{Demazure1}.

Let us define the divisors $Z_i=f_r^{-1}\cdots
f_{i+1}^{-1}\sigma_i(X_{i-1})$. The divisors $(Z_i)_{i\in[1,r]}$
have normal crossing (cf. for example \cite{Demazure1}). We have
$$Z_i=\left\{(x_1,\cdots,x_r)\in \Xt(\wt)\ /\ x_i=x_{p(i)}\right\}$$
in the configuration variety with $x_{p(i)}=1$ if $p(i)$ does not
exist. Then for any subset $K$ of $[1,r]$, one defines
$\displaystyle{Z_K=\bigcap_{i\in K}Z_i}$. Let us
recall the following (see for example \cite{Demazure1}):

\begin{fact}
\label{image}
The image by $\pi$ of $Z_K$ is the Schubert subvariety $X(y)$ where
  $y$ is the longuest element that can be writen as a subword of $\wt$
  without the terms $s_i$ for $i\in K$.
\end{fact}

Denote by $\xi_i$ the class of $Z_i$ in $A^*(\Xt(\wt))$. These
classes form a basis of the Chow ring.

\begin{theo}
{\rm (Demazure \cite{Demazure1} Par. 4. prop. 1)}
  The Chow ring $A^*(\Xt(\wt))$ of $\Xt(\wt)$ is isomorphic over $\Z$ to
$$\frac{\Z[\xi_1,\cdots,\xi_r]}{\left(\xi_i\cdot\sum_{j=1}^i
  \scal{\a_j^\vee}{\a_i}\xi_j\ \ {\rm for}\ \ {\rm all}\ \
  i\in[1,r]\right)}.$$ 
\end{theo}

Let us recall the following result (which is no longer true if the
writing is not reduced):

\begin{prop}\cite{LauritzThomsen}
  The divisors $(\xi_i)_{1\leq i\leq r}$ form a basis of the cone of
  effective divisors. 
\end{prop}

Let us denote by $T_i$ the pull-back on $\Xt(\wt)$ of the relative tangent
sheaf of the fibration $f_i$. Define the classical sequence of roots
$(\a_i)_{i\in[1,r]}$ associated to $\wt$ by $\a_1=\b_1$,
$\a_2=s_{\b_1}(\b_2)$, $\dots$, $\a_r=s_{\b_1}\cdots
s_{\b_{r-1}}(\b_r)$. We denote by $C_i$ the curve
$Z_K$ with $K=[1,r]\setminus\{i\}$ and recall some formulae in the ring
$A^*(\Xt(\wt))$ given in \cite{Perrin3} corollary 3.8 and propositions
3.3 and 3.11:

\begin{prop}
\label{intersection}
\label{intersectiontgt}
  We have the following formulae:
$$[C_i]\cdot\xi_j=\left\{\begin{array}{cc}
0&{\rm for}\ i>j\\
1 &{\rm for}\ i=j\\
\sca{\b_i^\vee}{\b_j} &{\rm for}\ i<j
\end{array}\right.,\ \ \ [C_i]\cdot T_j=\left\{\begin{array}{cc}
0&{\rm for}\ i>j\\
\sca{\b_i^\vee}{\b_j} &{\rm for}\ i\leq j
\end{array}\right.\ {\rm and}\ \ T_i=\sum_{k=1}^{i}\sca{\a_k^\vee}{\a_i}\cdot
\xi_k.$$
\end{prop}

This proposition and the formula
$\displaystyle{c_1(\Xt(\wt))=\sum_{i=1}^rT_i}$ gives:
$$c_1(\Xt(w))=-K_{\Xt(w)}=\sum_{k=1}^r\left(\sum_{i=k}^r
  \sca{\a_k^\vee}{\a_i}\right)\xi_k.$$

\subsubsection{Ample divisors}
\label{effectif}

In this paragraph, we describe the ample divisors on $\Xt(\wt)$. This
has already been done in \cite{LauritzThomsen} but we rephrase it in
terms of configuration varieties. We also get a description of the
Mori cone (see \cite{Matsuki} for references on this cone).

We have natural morphisms $p_i:\Xt(\wt)\to G/P_{\b_i}$ and as
$P_{\b_i}$ is maximal, the Picard group of $G/P_{\b_i}$ is generated
by a very ample divisor $\oo_{G/P_{\b_i}}(1)$ and we define on
$\Xt(\wt)$ the invertible sheaf
$\L_i=p_i^*(\oo_{G/P_{\b_i}}(1))$. These sheaves will form a basis of
the ample cone.

We also define particular curves $Y_i$ for $i\in[1,r]$ on $\Xt(\wt)$ by:
$$Y_i=\left\{(x_1,\cdots,x_r)\in\Xt(\wt)\ /\ x_j=1\ {\rm for}\ j\neq i
\right\}.$$
The  following lemma shows that $Y_i$ is a curve isomorphic to
$\p(1,\b_i)$.

\begin{lemm}
  For any $x_i\in\p(\bar1,\b_i)$, the element $(x_j)_{j\in[1,r]}$ of
    $\prod_{j=1}^rG/P_{\b_j}$ such that $x_j=\bar1$ for all $j\neq i$
    is in the configuration variety $\Xt(\wt)$.
\end{lemm}

\begin{preu}
We only have to prove that for any $x_i$ in
$\p(\bar1,\b_i)=B_{\b_i}/B$, we have $\bar1\in\p(x_i,\b_{i+1})$.

 The element $x_i$ can be lifted in $b_i\in B_{\b_i}$. The elements of
 $\p(x_i,\b_{i+1})$ are the classes of elements of the form
 $b_ib_{i+1}\in B_{\b_i}B_{\b_{i+1}}$. If $\b_{i+1}\neq\b_i$ (this is
 always the case if the writing is reduced), then $B_{\b_i}\subset
 P_{\b_{i+1}}$. In this case we set $b_{i+1}=1$ so that the class of
 $b_ib_{i+1}\in\p(x_i,\b_{i+1})$ is $\bar 1$ in $G/P_{\b_{i+1}}$. If
 $\b_{i+1}=\b_i$, then we set $b_{i+1}=b_i^{-1}$ to get the result.
\end{preu}

The definitions of the curves $Y_i$ and the line bundles $\L_i$ yield
to following:

\begin{prop}
\label{duale}
  We have the formula $\L_i\cdot [Y_j]=\delta_{i,j}$. In other words
  the families $(\L_i)_{i\in[1,r]}$ and $([Y_i])_{i\in[1,r]}$ are dual
  to each other.
\end{prop}

Let us prove that the family $([Y_i])_{i\in[1,r]}$ forms a basis of
$A_1(\Xt(\wt))$.

\begin{prop}
\label{courb-eff}
  For all $i\in[1,r]$, we have $[Y_i]=[C_i]-[C_{s(i)}]$ (where
  $[C_{s(i)}]=0$ if $s(i)$ doesn't exist).

In consequence, the classes $([Y_i])_{i\in[1,r]}$ form a basis of
  $A_1(X_Q)$ over $\Z$ and the classes $(\L_i)_{i\in[1,r]}$ form a
  basis of $A^1(X_Q)$ over $\Z$.
\end{prop}

\begin{preu}
On the fist hand, the curve $C_i$ is given by the equations
$x_j=x_{p(j)}$ for $j\neq i$. This means that for $j<i$, we have
$x_j=1$ and for all $j$ with $\b_j\neq\b_i$, we also have $x_j=1$. The
only indices $k$ for which $x_k$ may be different from $1$ are such
that $k=s^n(i)$ for some $n\in\N$. For such a $k$, we have the
equality $x_k=x_i$. Denote by $n(i)$ the biggest integer $n$ tsuch
that $s^n(i)$ exists. The curve $C_i$ (resp. $C_{s(i)}$) is the
diagonal in the product 
$$\prod_{k=0}^{n(i)}\p(1,\b_{s^k(i)})\ \ \ \ {\rm resp.}\ \ \ \
\prod_{k=1}^{n(i)}\p(1,\b_{s^k(i)}).$$

On the other hand, the curve $Y_i$ corresponds to the first factor of
the first product. In this product we thus have the required equality.
\end{preu}

We can now describe the ample cone (see also \cite{LauritzThomsen})
and the cone of effective curves.

\begin{coro}
  The cone of ample divisors is generated by the classes $\L_i$ and
  the cone of effective curves is generated by the classes
  $[Y_i]$. All ample divisors are very ample.
\end{coro}

\begin{preu}
  Let $D$ be ample on $\Xt(\wt)$, then $a_i=D\cdot [Y_i]$ is a
  positive integer. Because of proposition \ref{duale}, we have
  ${D=\sum_{i=1}^ra_i\L_i}$ 
and $D$ lies in the cone generated by the $\L_i$. 

Conversely, any
divisor ${\sum_{i=1}^ra_i\L_i}$ with $a_i>0$ gives the
embedding of $\Xt(\wt)$ obtained by composing the inclusion in the
product ${\prod_{i=1}^rG/P_{\b_i}}$ with the Veronese
morphism given by the very ample sheaf
${\bigotimes_{i=1}^r\oo_{G/P_{\b_i}}(a_i)}$.  

In the same way we get the result on effective curve.
\end{preu}

Finaly we calculate the divisors $\L_i$ in terms of the basis
$(\xi_k)_{k\in[1,r]}$:

\begin{prop}
\label{coord}
  The $k^{th}$ coordinate of $\L_i$ in the basis $(\xi_i)_{i\in[1,r]}$
  is 0 if $k>i$, 1 if $k=i$ and is given by the following formulae if
  $k<i$ and $\b_k=\b_i$ (resp. $\b_k\neq\b_i$):
$$1+\sum_{j=k+1,\ \b_j=\b_i}^i\sca{\a_k^\vee}{\a_j}\ \ \left({\rm resp.} \ 
\sum_{j=k+1,\ \b_j=\b_i}^i\sca{\a_k^\vee}{\a_j}\right).$$

In particular we have the following simple formula
$$\L_r=\sum_{k=1}^r\xi_k.$$
\end{prop}

\begin{preu}
Let us recall from \cite{Perrin3} lemma 4.5 that the following classes
  of curves $[\Ch_i]=[C_i]+\sum_{k=i+1}^n\sca{\a_i^\vee}{\a_k}[C_k]$
  form a dual basis to $(\xi_i)_{i\in[1,r]}$. The $k^{th}$ coordinate
  is thus given by the intersection $\L_i\cdot[\Ch_k]$.

For this we will need the formula coming directely from propositions
\ref{duale} and \ref{courb-eff}
$$\L_i\cdot[C_j]=\left\{\begin{array}{cc}
1&{\rm for}\ i>j\ {\rm and}\ \b_i=\b_j\\
0&{\rm otherwise.}
\end{array}\right.$$
Applying this gives the first formula. For the case of $\L_r$, the
formula is a consequence of the following formula from \cite{Perrin3}
corollary 2.18:
$$\sum_{j=k+1,\ \b_j=\b_r}^r\sca{\a_k^\vee}{\a_j}=\left\{
  \begin{array}{cc}
1&\textit{if}\ \b_k\neq\b_r\\
0&\textit{if}\ \b_k=\b_r.
  \end{array}\right.$$
\end{preu}

\section{Geometry of the quiver}

In this paragraph, we give an explicit description of the quiver
$Q_\wt$ given by the reduced writing $w=s_{\b_1}\cdots s_{\b_r}$ of
the shortest element in the class $\bar w\in W/W_{P_\pd}$. We also define
some invariants of the quiver and deduce some consequences on the
geometry of the Schubert variety.

\subsection{Minuscule conditions on the quivers}

Set $w_i=s_{\b_i}\cdots s_{\b_r}$ for $i\in[1,r]$ and $w_{r+1}=1$ and
let us first recall the following result from \cite{GofG/P3} (proof of
theorem 3.1):

\begin{fact}
We have $\sca{\b_{i-1}^\vee}{w_i(-\pd)}=-1$ for all $i\in[2,r+1]$. In
consequence we have for all $i\in[2,r]$:
$$w_i(-\pd)=-\pd+\b_r+\cdots+\b_{i}.$$
\end{fact}

The following proposition describes all possible quivers for minuscule
Schubert varieties.

\begin{prop}
\label{geom-carq}
Geometry of the quiver.
  \begin{itemize}
  \item[(\i)] There is no arrow from the vertex $r$ and $\b_r$ is the
  unique simple root with $\sca{\b_r^\vee}{\pd}=1$.
  \item[(\i\i)] If a vertex $i<r$ of the quiver is such that $s(i)$ does not
  exist, then there is a unique arrow from $i$. If $k$ is the end of
  the arrow we have $\sca{\b_i^\vee}{\b_k}=-1$.
  \item[(\i\i\i)] If a vertex $i$ of the quiver is such that $s(i)$ exists,
  then there are exactely two arrows from $i$. If $k_1$ and $k_2$ are
  the end of these arrows we have
  $\sca{\b_i^\vee}{\b_{k_1}}=\sca{\b_i^\vee}{\b_{k_2}}=-1$.
  \end{itemize}
\end{prop}

\begin{preu}
(\i) The previous fact shows that we have $\sca{\b_r^\vee}{\pd}=1$.

(\i\i) Let $i$ be such a vertex. In particular we have $\b_i\neq\b_r$ and
   $\sca{\b_i^\vee}{\pd}=0$. The previous fact gives
   $\scal{\b_i^\vee}{-\pd+\b_r+\cdots+\b_{i}}
   =\sca{\b_{i}^\vee}{w_{i+1}(-\pd)}=-1$ and thus
$$\sum_{k=i+1}^r\sca{\b_i^\vee}{\b_k}=-1.$$
We conclude because every term of this sum have to be either 0 or $-1$.

(\i\i\i) Let $i$ be such a vertex. The same calculation as above shows that
   $\sum_{k=i+1}^r\sca{\b_i^\vee}{\b_k}=
   \sum_{k=s(i)+1}^r\sca{\b_i^\vee}{\b_k}=-1$ 
   if $\b_i\neq\b_r$ and $\sum_{k=i+1}^r\sca{\b_i^\vee}{\b_k}=
   \sum_{k=s(i)+1}^r\sca{\b_i^\vee}{\b_k}=0$ if $\b_i=\b_r$. In
   particular we always have
$$\sum_{k=i+1}^{s(i)}\sca{\b_i^\vee}{\b_k}=0.$$
We conclude because the only positive term is
$\sca{\b_i^\vee}{\b_{s(i)}}=2$ and every other term of this sum have
to be either 0 or $-1$.
\end{preu}

\begin{rema}
\label{unicite}
(\i) One can deduce from this result (see \cite{moivirtuel}) the fact
  already proved by J. R. Stembridge \cite{stem1} that the reduced
  writing $\wt$ is unique modulo commutation relations (there is no
  braid relation). We can thus write $Q_w$ instead of $Q_\wt$ and call
  it the quiver associated to the minuscule Schubert variety $X(w)$.

Furthermore, because of this unicity and the fact that between
$i$ and $s(i)$ there are always two vertices, the writing deduced from
a quiver satisfying the conditions of the preceding proposition is
always reduced and the quivers satisfying the conditions are always
quivers associated to a minuscule Schubert variety.

(\i\i) A minuscule quiver is always connected: there is a path from
any vertex $i$ to the last vertex $r$.
\end{rema}

\subsection{Combinatoric decription of the minuscule quivers}

It is easy from proposition \ref{geom-carq} to describe the quivers
$Q_\pd$ of a minuscule homogeneous variety $G/P_\pd$ (see appendix for
a list of these quivers). We now describe the quivers of minuscule
Schubert varieties in $G/P_\pd$ as subquivers of $Q_\pd$.
Define a natural partial order on the quiver:

\begin{defi}
\label{order}
(\i) We denote by $\preccurlyeq$ the partial order on the vertices of
the quiver generated by the relations $i\preccurlyeq j$ if there
exists an arrow from $i$ to $j$.

(\i\i) Let $A$ be a totaly unordered set of vertices of the
quiver $Q_\pd$ for the partial order
$\preccurlyeq$. We denote by
$Q_{\bar A}$ the full subquiver of $Q_\pd$ with vertices $i\in Q_\pd$
such that there exists $a\in A$ with
$i\preccurlyeq a$ and by $Q_{A}$ the full subquiver of $Q_\pd$ whose
vertices are not vertices of $Q_{\bar A}$.
\end{defi}

\begin{prop}
\label{description}
The quiver of Schubert varieties in $G/P_\pd$ are in one to one
correspondence with the subquivers $Q_A$ of $Q_\pd$ for $A$ any totaly
unordered set of vertices of $Q_\pd$.

With this correspondence, the bruhat order is given by the inclusion
of quivers.
\end{prop}

\begin{preu}
Let $X(w)\subset G/P_\pd$ a Schubert variety and denote by $w_0\in W$
the only element such that $X(w_0)=G/P_\pd$. There exists a sequence
$(\b_1,\cdots\b_i)$ of simple roots such that $w_0=s_{\b_1}\cdots
s_{\b_i}w$. Taking a reduced writing $w=s_{\b_{i+1}}\cdots s_{\b_r}$
of $w$ we get a reduced writing of $w_0$ which is unique modulo
commutation relations.

The vertices of the quiver $Q_\pd$ are indexed by $[1,r]$. Denote by
$A=\{i_1,\cdots,i_k\}$ the set of maximal elements for the partial
order $\preccurlyeq$ of the set $[1,i]$. The set $A$ is totaly
unordered and the quiver associated to $X(w)$ is $Q_A$. 

The fact that this is a one to one correspondence comes from the
unicity of the reduced writing.
\end{preu}

\begin{rema}
  As a corollary one can prove (see \cite{moivirtuel}) the following
  classical result on the Bruhat order 
  (see for example \cite{GofG/P3}): the Bruhat order in $W/W_P$ is
  generated by  simple reflexions. In other words, all Schubert
  divisors are mobile.
\end{rema}

Finaly let us define some particular vertices of these quivers. In the
following definition $Q_w$ is the quiver of a minuscule Schubert
variety $X(w)$.

\begin{defi}
(\i) We call pic any vertex of $Q_w$ minimal for the partial order
  $\preccurlyeq$. We denote by $p(Q_w)$ the set of pics of $Q_w$.

(\i\i) We call hole of the quiver $Q_w$ any vertex $i$ of $Q_w$ such 
that $p(i)$ does not exist and there are exactely two
vertices $j_1\preccurlyeq i$ and $j_2\preccurlyeq i$ in $Q_w$ with
$\sca{\b_i^\vee}{\b_{j_k}}\neq0$ for $k=1,2$. 

We will also call a hole of $Q_w$ any $i\in Q_\pd\setminus Q_w$
  such that $s(i)$ does not exist in $Q_\pd$ and
  $\b_i\in\partial(G_w)$. Such a hole will be called a virtual
  hole. We denote by $t(Q_w)$ the set of holes of $Q_w$.

(\i\i\i) The height $h(i)$ of a vertex $i$ is the largest $n$ such that
there exists a sequence $(i_k)_{k\in[1,n]}$ of vertices with $i_1=1$,
$i_n=r$ and such that there is an arrow from $i_k$ to $i_{k+1}$ for
all $k\in[1,n-1]$.
\end{defi}

\begin{rema}
(\i) If $Q_w=Q_A$ as in definition \ref{order} then $t(Q_w)=A$. 

(\i\i) The height is well defined because there is at least one path
from any vertex $i$ to the last vertex $r$.
\end{rema}

The following proposition gives a recursive way to calculate the
height of a vertex.

\begin{prop}
\label{hauteur}
  Let $Q$ be a quiver associated to a minuscule Schubert variety and
  $i$ a vertex of this quiver. 
We have the following cases:
\begin{itemize}
\item if $s(i)$ does not exist, then there exists a unique
  $k\succcurlyeq i$ with $\sca{\b_k^\vee}{\b_{i}}=-1$ and we have
$$h(i)=h(k)+1.$$
\item If $s(i)$ exists, then there exists a non negative integer $n$
  and a sequence $(j_k,j'_k)_{k\in[0,n+1]}$ of vertices with
  $j_0=i$, $j'_k=s(j_k)$ for $k\in[0,n]$,
  $\b_{j_{n+1}}\neq\b_{j'_{n+1}}$, $\sca{\b_{j_k}}{\b_{j_{k+1}}}=-1$
  and $\sca{\b'_{j_k}}{\b'_{j_{k+1}}}=-1$ for $k\in[0,n]$ and
  $j_0\preccurlyeq\cdots\preccurlyeq j_n\preccurlyeq j_{n+1},
  j'_{n+1}\preccurlyeq j'_n\preccurlyeq\cdots\preccurlyeq j'_0$. In
  this case we have 
$$h(j_k)=2n+2-k+h(s(i))\ \ {\rm and}\ \ h(j'_k)=k+h(s(i))\ \
\textrm{foll all}\ k\in[0,n+1].$$ 
\end{itemize}
\end{prop}

\begin{preu}
We prove these formulae by descending induction on $i$. If $i=r$ then
$h(i)=1$. Assume the proposition is true for all $j>i$. In the first
case, any sequence of arrows from $i$ to $r$ has to pass through the
vertex $k$ and we have $h(i)=h(k)+1$.

The first to prove in the second case is the existence of the sequence
$(j_k,j'_k)_{k\in[0,n+1]}$. Let us denote by $j_1$ and $j'_1$ the two
vertices $j$ such that there is an arrow from $i$ to $j$. If
$\b_{j_1}\neq\b_{j'_1}$ then set $n=0$ and we are done. Otherwise,
assume (for example) that $j_1<j'_1$ then $j'_1=s(j_1)$. Indeed, if it was not
the cas there would exist $k\in[j_1,j'_1]$ and in particular $k<s(i)$
with $\b_k=\b_{j_1}$ thus $\sca{\b_i^\vee}{\b_k}=-1$. By construction
of the quiver, there must be an arrow from $i$ to $k$ and thus at
least three arrows from $i$. This is impossible by proposition
\ref{geom-carq}. We can construct with $(j_1,j'_1)$ a couple
$(j_2,j'_2)$ in the same way and by induction a sequence
$(j_k,j'_k)_{k}$. As long as $\b_{j_k}=\b_{j'_k}$ we can go on. This
has to stop because the quiver is finite.

The formula is now clear because any sequence from $i$ to $r$ as to go
through $j_1$ or $j'_1$. As the height of $j_1$ is bigger than the one
of $j'_1=s(j_1)$ by induction we must have $h(i)=h(j_1)+1$ and we
conclude one more time by induction. 
\end{preu}

\begin{rema}
\label{order-comm}
  By changing the order of commuting factors, we may assume in the
  preceding proposition that $k=i+1$ in the first case and that
  $j_k=i+k$ and $j'_k=i+2n+3-k$ for all $k\in[0,n+1]$ in the second
  one.
\end{rema}

We can now describe the stabiliser of a Schubert varietie $X(w)$
thanks to its quiver $Q_w$.

\begin{prop}
\label{stab}
  Let $J$ be the set of simple roots not in $\b(t(Q_w))$. The stabiliser
  of $X(w)$ is the parabolique subgroup $P_J$ generated by $B$ and the
  groups $U_{-\b}$ with $\b\in J$.
\end{prop}

\begin{preu}
A simple root $\b$ is such that $U_{-\b}\subset P$ if and only if
$s_\b w<w$ (for the Bruhat order). But from unicity of reduced writing
and our characterisation (proposition \ref{geom-carq}) of quivers
associated to a reduced writing, we see that this implies that $\b$
has to be in $\b(t(Q_w))$.
\end{preu}

\begin{coro}
\label{svds}
  Let $Q_{w'}$ be the quiver of a Schubert subvariety $X(w')$ of
  $X(w)$ stable under $P_J={\rm Stab}(X(w))$. Then
  $\b(t(Q_{w'}))\subset\b(t(Q_w))$. 
\end{coro}

\subsection{Weil and Cartier divisors}

In this paragraph, we describe thanks to the quiver the Weil and
Cartier divisors of a minuscule Schubert variety $X(w)$. We also
compute the canonical sheaf of $X(w)$.

\begin{prop}
\label{ample-div}
  The group $\weil(X(w))$ of Weil-divisors is the free $\Z$-module
  generated by the classes $D_i:=\pi_*\xi_i$ for $i\in p(Q_w)$

The Picard group $\pic(X(w))\subset\weil(X(w))$ is isomorphic to $\Z$
and is generated by the element
$\L(w):=\pi_*\L_r=\oo_{G/P_{\b_r}}(1)\vert_{X(w)}$. We have the formula
$$\L(w)=\sum_{i\in p(Q_w)}D_i.$$
\end{prop}

\begin{preu}
It is well known (see for example \cite{Brion}) that the Picard group
is isomorphic to $\Z$ and generated by $\L(w)$ and that the group of
Weil divisors is free generated by the divisorial Schubert
varieties. These varieties are the images by $\pi:\Xt(\wt)\to X(w)$ of
the non contracted divisors $Z_i$. Now the divisor $Z_i$ is the
configuration variety obtained from the quiver $Q$ with the vertex $i$
removed (the arrows are reoganised as in definition
\ref{quiver}). According to the fact \ref{image}, the image of $Z_i$
is not contracted if and only if the quiver is reduced (ie correspond
du a reduced writing). It is clear that this can only be the case for
$i\in p(Q_w)$.

The last formula is an application of corollary \ref{coord}.
We recover the well known fact that all Schubert divisors are of
multiplicity one in the minuscule case.
\end{preu}

\begin{coro}
\label{goren}
  A minuscule Schubert variety is locally factorial if and only if its
  quiver has a unique pic.
\end{coro}

\subsection{Canonical divisor}

The Schubert varieties are singular and in general not
Gorenstein (see \cite{Alexalex} for a characterisation of Gorenstein
Schubert varieties for $GL_n$). We can therefore not define the
canonical divisor as a Cartier divisor.

The canonical divisor $K_{X(w)}$ of a Schubert variety $X(w)$ is well
defined as a Weil divisor thanks to the divisor of a 1-form on
$X(w)$. The properties of Schubert varieties (they are normal,
Cohen-Macaulay with rational singularities) and the Bott-Samelson
resolution $\pi:\Xt(\wt)\to X(w)$ enables however to calculate
$K_{X(w)}$ by $K_{X(w)}=\pi_*(K_{\Xt(\wt)})$ (see for example
\cite{BrionKumar} paragraph 3.4).

Let us denote by $h(w)$ the lowest height of a pic in $Q_w$ (the
quiver associated to $X(w)$), we have the following:

\begin{prop}
\label{cano-ample}
  We have the formula
$$-K_{X(w)}=\sum_{i\in p(Q_w)}(h(i)+1)D_i=(h(w)+1)\L(w)+\sum_{i\in
  p(Q_w)}(h(i)-h(w))D_i.$$
\end{prop}

\begin{preu}
The second part of the formula comes from the first one and
proposition \ref{ample-div}.

To prove the first part, we use the fact that
$K_{X(w)}=\pi_*(K_{\Xt(\wt)})$ and the formula of paragraph
\ref{canonique}. We are left to prove the following: 

\begin{lemm}
\label{som-haut}
We have the formula:
$$\sum_{k=i}^r\sca{\a_i^\vee}{\a_k}=\sum_{k\succcurlyeq i}
\sca{\a_i^\vee}{\a_k}=h(i)+1$$
\end{lemm}

\begin{preu}
We proceed by descending induction and use the proposition
\ref{hauteur} and the remark \ref{order-comm}. We have the following two
cases:
\begin{itemize}
\item if $s(i)$ does not exist, then  with $\sca{\b_{i+1}^\vee}{\b_{i}}=-1$.
\item If $s(i)$ exists, then there exists a non negative integer $n$
  and a sequence $(j_k,j'_k)_{k\in[0,n+1]}$ of vertices with $j_0=i$,
  $j'_k=s(j_k)$ for $k\in[0,n]$, $\b_{j_{n+1}}\neq\b_{j'_{n+1}}$,
  $\sca{\b_{j_k}}{\b_{j_{k+1}}}=-1$ and
  $\sca{\b'_{j_k}}{\b'_{j_{k+1}}}=-1$ for $k\in[0,n]$ and
  $j_0\preccurlyeq\cdots\preccurlyeq j_n\preccurlyeq j_{n+1},
  j'_{n+1}\preccurlyeq j'_n\preccurlyeq\cdots\preccurlyeq
  j'_0$. Furthermore we may assume that $j_k=i+k$ and $j'_k=i+2n+3-k$
  for all $k\in[0,n+1]$.
\end{itemize}

We proceed by descending induction on $i$. If $i=r$, there is only
one term and the sum is $\sca{\a_r^\vee}{\a_r}=2=h(r)+1$. We assume
that the formula is true for all $j\geq i+1$. Let us use the following
sequence $\at_j=s_{\b_r}\cdots s_{\b_{r-j+2}}(\b_{r-j+1})$ already
introduced in \cite{Perrin3} and satisfying the equality
$\sca{\at_{r+1-k}^\vee}{\at_{r+1-i}}=\sca{\a_k^\vee}{\a_i}$.

Calculating $\at_{r+1-i}=s_{\b_r}\cdots s_{\b_{i+1}}(\b_i)$ we find
$$\at_{r+1-i}=\left\{
  \begin{array}{ll}
\at_{r-i}+s_{\b_r}\cdots s_{\b_{i+2}}(\b_i)& \textrm{in the first case;}\\
\at_{r-i-n}+\at_{r-i-n-1}-\at_{r-i-2n-2}& \textrm{in the second case;}\\
  \end{array}\right.$$

Let us now calculate the sum
$$\sum_{k\succcurlyeq i} \sca{\a_k^\vee}{\a_i}=\sum_{k\succcurlyeq i}
\sca{\at_{r+1-k}^\vee}{\at_{r+1-i}}.$$

In the first case, denote $\a=s_{\b_r}\cdots s_{\b_{i+2}}(\b_i)$. If
$j\geq i+2$, then we have $\sca{\b_i^\vee}{{\b_j}}=0$ so that
$\a=\b_i$. Furthermore, the root $\at_{r+1-j}$ is a sum of simple
roots contained in the set $\{\b_j,\cdots,\b_r\}$ so for $j\geq i+2$,
we have $\sca{\at_{r+1-j}^\vee}{\a}=0$. In this case the sum equals:
$$\sum_{k\succcurlyeq i}
\sca{\at_{r+1-k}^\vee}{\at_{r+1-i}}=\sum_{k\succcurlyeq i}
\sca{\at_{r+1-k}^\vee}{\at_{r-i}}+\sum_{k\succcurlyeq
  i}\sca{\at_{r+1-k}^\vee}{\a}$$
$$=\sum_{k\succcurlyeq i+1} \sca{\at_{r+1-k}^\vee}{\at_{r-i}}+
\sca{\at_{r+1-i}^\vee}{\at_{r-i}}+ \sca{\at_{r+1-i}^\vee}{\a}+
\sca{\at_{r-i}}{\a}.$$ 
$$=h(i+1)+1-\sca{\b_i^\vee}{\b_{i+1}}+
\sca{(\b_i+\b_{i+1})^\vee}{\b_i}+\sca{\b_{i+1}^\vee}{\b_{i}}$$
$$=h(i+1)+1+1+1-1=h(i+1)+2.$$ 

In the second case, it is an easy exercice to see that the roots
$\{\b_k\}_{k\in[i,i+n+2]}$ form a diagram of type $D_{n+2}$ (with the
notations of \cite{bourb} the root $\b_{i+k}$ is the $(k+1)$-th root of
the diagram). We can then calculate 
$$\at_{r+1-i-k}=\left\{
\begin{array}{cc}
\displaystyle{s_{\b_r}\cdots s_{\b_{i+2n+4}}
  \left(\sum_{j=0}^{2n+2-k}\b_{i+j}\right)}
& \textrm{for all } k\in[0,n+1]\\
\displaystyle{s_{\b_r}\cdots s_{\b_{i+2n+4}}
  \left(\sum_{j=0}^{2n+3-k}\b_{i+j}\right)}
& \textrm{for all } k\in[n+2,2n+3]\\
\end{array}\right.$$

The sum is in this case:
$$\sum_{k\succcurlyeq i}
\sca{\at_{r+1-k}^\vee}{\at_{r+1-i}}=\sum_{k\succcurlyeq i}
\sca{\at_{r+1-k}^\vee}{\at_{r-i-n}}+\sum_{k\succcurlyeq i}
\sca{\at_{r+1-k}^\vee}{\at_{r-i-n-1}}-\sum_{k\succcurlyeq i}
\sca{\at_{r+1-k}^\vee}{\at_{r-i-2n-2}}$$
$$=\sum_{k\succcurlyeq i+n+1}
\sca{\at_{r+1-k}^\vee}{\at_{r-i-n}}+\sum_{k\succcurlyeq i+n+2}
\sca{\at_{r+1-k}^\vee}{\at_{r-i-n-1}}-\sum_{k\succcurlyeq i+2n+3}
\sca{\at_{r+1-k}^\vee}{\at_{r-i-2n-2}}$$
$$+\sum_{k=0}^{n}\sca{\at_{r+1-k}^\vee}{\at_{r-i-n}}+
\sum_{k=0}^{n+1}\sca{\at_{r+1-k}^\vee}{\at_{r-i-n-1}}-
\sum_{k=0}^{2n+2}\sca{\at_{r+1-k}^\vee}{\at_{r-i-2n+2}}.$$ 
The description of $\at_{r+1-i-k}$ shows that all roots
$\at_{r+1-i-k}$ for $k\in[0,n]$ (resp. $k\in[1,2n+2]$) have
intersection 1 with $\at_{r-i-n}$ and $\at_{r-i-n-1}$
(resp. $\at_{r-i-2n-2}$) and we have the intersections
$\sca{\at_{r-i-n}^\vee}{\at_{r-i-n-1}}=0$ and
$\sca{\at_{r+1-i}^\vee}{\at_{r-i-2n-2}}=0$. This gives us the
following formulae:
$$\sum_{k=0}^{n}\sca{\at_{r+1-k}^\vee}{\at_{r-i-n}}=n+1\ ;\ \ 
\sum_{k=0}^{n+1}\sca{\at_{r+1-k}^\vee}{\at_{r-i-n-1}}=n+1\ ;\ \ 
\sum_{k=0}^{2n+2}\sca{\at_{r+1-k}^\vee}{\at_{r-i-2n+2}}=2n+2.$$
Using the induction hypothesis we get:
$$\sum_{k\succcurlyeq i}
\sca{\at_{r+1-k}^\vee}{\at_{r+1-i}}=h(i+n+1)+1+h(i+n+2)+1-h(i+2n+3)-
1+n+1+n+1-2n+2$$
$$=h(i+n+1)+h(i+n+2)+1-h(i+2n+3).\ \ \ \ \ \ $$ 
We conclude in both cases thanks to proposition \ref{hauteur}.
\end{preu}
This lemma completes the proof of the proposition.
\end{preu}

\begin{coro}
\label{gorenstein}
  The Schubert variety $X(w)$ is Gorenstein if and only if all the
  pics of its quiver have the same height. In this case we have
  $K_{X(w)}=(h(w)+1)\L(w)$.
\end{coro}

\begin{preu}
  The variety $X(w)$ is Gorenstein if and only if its canonical
  divisor is Cartier. The preceding formula shows that this is
  equivalent to the fact that all the pics of its quiver have the same
  height.
\end{preu}

\begin{rema}
For $Gl(n)$, we recover a particular case of the result of A. Woo
  and A. Yong \cite{Alexalex} on Gorenstein Schubert varieties.
\end{rema}

\section{Generalisation of Bott-Samelson's construction}

In this paragraph, we are going to construct some varieties $\Xh(\wh)$
together with a birational morphism $\pih:\Xh(\wh)\to X(w)$. These
constructions generalise the $IH$-small resolutions constructed by 
A. Zelevinsky \cite{Zelevinsky} and P. Sankaran and P. Vanchinathan
\cite{SanVan}. 

Recall that the Bott-Samelson varieties can be seem as a tower of
$\pu$-fibrations coming from a reduced writing $\wt$ of an element
$w$. Many generalisations of this construction (for example in
\cite{Zelevinsky} and \cite{SanVan} but also in
\cite{Perrin-resol} and even in the general construction of
C. Contou-Carr{\`e}re \cite{Contoucarrere2}) are
constructed as towers of locally trivial fibrations with
fibers isomorphic to a fixed homogeneous variety thanks to a more
general decomposition $\wh$ of $w$ as a product of elements in the
Weyl group. For the varieties $\Xh(\wh)$ we make the same construction
but we allow locally trivial fibrations with fiber a locally factorial
or gorenstein Schubert variety.

\subsection{Elementary construction}

Let us explain the following elementary construction. As in
\cite{Demazure1}, the construction of $\Xh(\wh)$ will simply be a
successive application of this elementary construction.
Let $u\in W$ and $Y$ a variety with action of $P_Y$ a parabolic
subgroup of $G$ containing $B$ and assume that $P^u\cap G_u\subset
P_Y$. We define
$$\Yh(u)=\overline{(P_u\cap G_u)u(P^u\cap G_u)}\times^{(P^u\cap
  G_u)}Y.$$

\begin{lemm}
\label{stabilisor}
(\i) The variety $\Yh(u)$ is a locally
trivial fibration over $X(u)$ with fibers isomorphic to $Y$.

(\i\i) Define the parabolic $P_{\Yh(u)}$ by
  $\Sigma(P_{\Yh(u)})=(\Sigma(P_Y)\cap\supp(u)^c)\cup
  \partial(G_u)\cup(\Sigma(P_u) \cap\supp(u))$.
There is an action of $P_{\Yh(u)}$ on $\Yh(u)$.
\end{lemm}

\begin{preu}
(\i) The first part of the proposition comes from the isomorphism
beetwen $X(u)$ and $\overline{(P_u\cap G_u)u(P^u\cap G_u)}/{(P^u\cap
  G_u)}$. 

(\i\i) For the second part, let us remark that one can replace the
groups $P_u\cap G_u$ and $P^u\cap G_u$ by bigger groups $A$ and $B$
such that the natural map $\overline{(P_u\cap G_u)u(P^u\cap
  G_u)}/{(P^u\cap G_u)}\to \overline{AuB}/B$ is an isomorphism and
$B\subset P_Y$. For
example, we take $A$ and $B$ such that
$\Sigma(A)=(\Sigma(P_Y)\cap\supp(u)^c)\cup
\partial(G_u)\cup(\Sigma(P_u) \cap\supp(u))$ and
$\Sigma(B)=(\Sigma(P_Y)\cap\supp(u)^c)
\cup(\Sigma(P^u) \cap\supp(u))$. We have the required isomorphism and
$B\subset P_Y$ (simply because $\Sigma(P_Y)\subset\Sigma(B)$). Then
$\Yh(u)$ is isomorphic to $\overline{AuB}\times^BY$ and $A$ acts on
$\Yh(u)$.
\end{preu}

\subsection{Construction of the resolution}

Let us give the following:

\begin{defi}
(\i) Let $w\in W$, a writing $w=w_1\cdots w_n$ where for all
$i\in[1,n]$, we have $w_i\in W$ is called a generalised decomposition
and denoted by $\wh$. If moreover we have the equality
$\displaystyle{l(w)=\sum_{i=1}^nl(w_i)}$ then we will say that the
generalised decomposition is reduced.

(\i\i) Let us associate to any generalised decomposition a sequence of
parabolic subgroups $(P_i)_{i\in[1,n]}$ defined by $P_n=P_{w_n}$ and 
$$\Sigma(P_{i})=(\Sigma(P_{i+1})\cap\supp(w_{i})^c)
\cup\partial(G_{w_i})\cup(\Sigma(P_{w_i})\cap\supp(w_i)).$$

(\i\i\i) We will say that a generalised reduced decomposition is
admissible if for all $i\in[1,n-1]$ we have $P^{w_i}\cap G_{w_i}\subset
P_{i+1}$. 

(\i v) We will say that such a generalised reduced decomposition is
good if for all $i\in[1,n-1]$ we have $P^{w_i}\cap G_{w_i}\subset
P_{w_{i+1}\cdots w_n}$ and
$\partial(G_{w_i})\subset\Sigma(P_{w_i\cdots w_n})$. 
\end{defi}

\begin{prop}
Let $\wh$ an admissible reduced generalised decomposition of a
minuscule element $w$ of $W$. One can
define by descending induction on $n$ the varieties $\Xh_i(\wh)$ by
$\Xh_n(\wh)=X(w_n)$ and for $i<n$ by $\Xh_i(\wh)=\Yh(w_i)$ where
$Y=\Xh_{i+1}(\wh)$. The group $P_{\Xh_i(\wh)}$ is the group $P_{i}$.

Furthermore, if the decomposition is good, then the group
$P_{\Xh_i(\wh)}$ is the group $P_{w_i\cdots w_n}$ and in particular
any good generalised reduced decomposition is admissible.
\end{prop}

\begin{preu}
  We proceed by induction. The variety $\Xh_n(\wh)$ is well defined
  and we have $P_{\Xh_n(\wh)}=P_{w_n}$. Assume that $\Xh_{i+1}(\wh)$
  is well defined and that $P_{\Xh_{i+1}(\wh)}=P_{i+1}$. To prove that
  $\Xh_i(\wh)$ exists, we have to prove that $P^{w_i}\cap
  G_{w_i}\subset P_{\Xh_{i+1}(\wh)}$ but it is the case by
  hypothesis. The fact that $P_{\Xh_i(\wh)}=P_i$ comes from lemma
  \ref{stabilisor}.

Now in the case of a good generalised reduced decomposition, we have
  to prove that $P_i=P_{\Xh_{i}(\wh)}=P_{w_{i}\cdots w_n}$. We know
  from lemma \ref{stabilisor} that
  $\Sigma(P_{\Xh_i(\wh)})=(\Sigma(\Xh_{i+1}(\wh))\cap\supp(w_i)^c)\cup
  \partial(G_{w_i})\cup(\Sigma(P_{w_i}) \cap\supp(w_i))$. Let
  $\b\in\Sigma(P_{w_i\cdots w_n})$. If $\b\in\supp(w_i)$ then $\b$ has
  to be a hole of the quiver of $w_i$ so that $\b\in\Sigma(P_{w_i})$
  and $\b\in\Sigma(P_{\Xh_i(\wh)})$. If $\b\not\in\supp(w_i)$ and
  $\b\not\in\partial(G_{w_i})$ then $s_\b$ commutes with $w_i$ and we
  have $\b\in\Sigma(P_{w_{i+1}\cdots w_n})$ and
  $\b\in\Sigma(P_{\Xh_i(\wh)})$. Finaly if $\b\in\partial(G_{w_i})$ we
  also have $\b\in\Sigma(P_{\Xh_i(\wh)})$.

Conversely, let $\b\in\Sigma(P_{\Xh_i(\wh)})$. If $\b\in\supp(w_i)$
then $\b\in\Sigma(P_{w_i})$ so $\b$ corresponds to a hole of the
quiver of $w_i$ and thus have to be a hole of the quiver of $w_i\cdots
w_n$. If $\b\in\partial(G_{w_i})$ we are done by hypothesis and finaly
if $\b$ is neither in $\supp(w_i)$ nor in $\partial(G_{w_i})$ then
$\b\in\Sigma(P_{w_{i+1}\cdots w_n})$ by induction hypothesis. Thus
$\b$ corresponds to a hole of the quiver of $w_{i+1}\cdots w_n$ and
does not appear in the quiver of $w_i$. It is thus still a hole of the
quiver of $w_{i}\cdots w_n$.
\end{preu}

\begin{defi}
  We denote by $\Xh(\wh)$ the variety $\Xh_1(\wh)$.
\end{defi}

\begin{coro}
The variety $\Xh(\wh)$ is a tower of localy trivial
 fibrations $f_i$ with fibers isomorphic to $X(w_i)$.  
\end{coro}

\begin{lemm}
\label{commute}
  Let $\wh$ an admissible reduced generalised decomposition of a minuscule
element $w$ of $W$. Let $i\in[1,n-1]$ and assume that for any couple
$(\b,\b')\in\supp(w_i)\times\supp(w_{i+1})$ we have
$\sca{\b^\vee}{\b'}=0$.

Then
  $w_iw_{i+1}=w_{i+1}w_i$, the generalised writing $\wh'$ given by
  $w=w'_1\cdots w'_n$ where $w'_k=w_k$ for $k\not\in\{i;i+1\}$,
  $w'_i=w_{i+1}$ and $w'_{i+1}=w_i$ is admissible and reduced and the
  morphisms $\pih:\Xh(\wh)\to X(w)$ and $\pih':\Xh(\wh')\to X(w)$ are
  the same.
\end{lemm}

\begin{preu}
  We simply have to look at the following situation. Let $A$ and $B$
  be parabolic subgroups of a group $G$ and $C$ and $D$
  be parabolic subgroups of a group $G'$. Assume that $B$ and $D$ act
  on a variety $X$ and consider the variety
  $\overline{AuB}\times^B\overline{CvD}\times^DX$ ($B$ acts on
  $\overline{CvD}\times^DX$ thanks to its action on $X$). It is
  isomorphic to
  $(\overline{AuB}\times\overline{CvD})\times^{B\times D}X$ and the
  construction is completely symetric. 

Let us remark that the variety $\Xh(\wh)$ is also isomorphic to the
variety $\Xh(\wh'')$ where $\wh''$ is such that $w''_k=w_k$ for $k<i$,
$w''_i=w_iw_{i+1}$ and $w''_k=w_{k+1}$ for $k>i+1$.
\end{preu}

Thanks to this lemma, we may assume that the support of any element
$w_i$ is connected (otherwise replace it by a product of elements
having a connected support).

\subsection{Link with the Bott-Samelson resolution}

In this paragraph, we show that the Bott-Samelson resolution
$\Xt(\wt)$ of a minuscule element factorises through any
pseudo-resolution $\Xh(\wh)$ constructed above. And we view this
variety as a projection from the Bott-Samelson resolution
$\Xt(\wt)$.

Let $\wh$ be an admissible generalised reduced writing of $w$ and let
us fix for any $i\in[1,n]$ a (unique) reduced writing
$w_i=s_{1,i}\cdots s_{r_i,i}$ denoted $\wt_i$.

\begin{lemm}
For any $i\in[1,n]$, the writing 
$$w=\left(\prod_{k=1}^i\prod_{j=1}^{r_k}s_{j,k}\right)
\cdot\prod_{k=i+1}^nw_k$$
denoted $\wh'_i$ is an admissible generalised reduced writing of $w$.
\end{lemm}

\begin{preu}
  Let us denote by $w'_k$ for $k\in[1,N]$ the terms of the generalised
  writing. It is clear that it is reduced. Let us prove that it is
  admissible. Because the writing $\wh$ is admissible, it is clear
  that the inclusion $P^{w'_{N-k}}\cap G_{w'_{N-k}}\subset P_{N-k}$
  holds for $k\leq i+2$. But for $k\geq i+1$ then $w'_{N-k}$ is a simple
  reflexion $s_\b$  and we have $G_{w'_{N-k}}=SL_2(\b)$ and
  $P^{w'_{N-k}}\cap G_{w'_{N-k}}$ is contained in the Borel $B$ so
  that the inclusion in $P_{N-k+1}$ is trivial.
\end{preu}

\begin{rema}
  Let us remark that the classical Bott-Samelson $\Xt(\wt)$ resolution
  is given by $\Xh(\wh'_1)$.
\end{rema}

\begin{prop}
  There is a morphism $\Xh(\wh'_i)\to\Xh(\wh'_{i+1})$ for all
  $i\in[1,n]$ (for $i=n$, let us fix $\Xh(\wh'_{n+1})=\Xh(\wh)$) and
  the morphism $\pi:\Xt(\wt)\to X(w)$ from the Bott-Samelson resolution
  to the Schubert variety factors through these morphisms. In
  particular we will denote by $\pit$ the morphism from $\Xt(\wt)$ to
  $\Xh(\wh)$.
\end{prop}

\begin{preu}
  The variety $\Xh(\wh'_i)$ is the quotient of the product 
$$\left(\prod_{k=1}^i\prod_{j=1}^{r_k}\overline{(P_{s_{j,k}}\cap
      G_{s_{j,k}})s_{j,k}(P^{s_{j,k}}\cap G_{s_{j,k}})}\right)
\times\prod_{k=i+1}^n\overline{(P_{w_k}\cap G_{w_k})w_k(P^{w_k}\cap
  G_{w_k})}$$
by the product
$$\left(\prod_{k=1}^i\prod_{j=1}^{r_k}P^{s_{j,k}}\cap G_{s_{j,k}}\right)
\times\prod_{k=i+1}^nP^{w_k}\cap
  G_{w_k}.$$
The action respects multiplication and in particular the
multiplication map on the $i^{\textrm{th}}$ factor
$$\prod_{j=1}^{r_i}\overline{(P_{s_{j,i}}\cap
      G_{s_{j,i}})s_{j,i}(P^{s_{j,i}}\cap G_{s_{j,i}})}\to
      \overline{(P_{w_i}\cap G_{w_i})w_i(P^{w_i}\cap 
  G_{w_i})}$$
and the identidy map on all the other factors is still defined modulo the
      action giving a map $\Xh(\wh'_i)\to\Xh(\wh'_{i+1})$. This map is
      simply the identity on all but one fibration (the one of fiber
      $X(w_i)$) and on this fibration is it given by the map
      $\Xt(\wt_i)\to X(w_i)$ from the
      Bott-Samelson resolution to the Schubert variety.

The last affirmation is a simple consequence of the associativity of
the product: the morphism from $\Xt(\wt)=\Xh(\wh'_1)$ is given by the
product of all the terms and the factorisations are given by making
the product in a certain order.
\end{preu}

\begin{rema}
\label{rem-proj}
(\i) The morphism $\pih:\Xh(\wh)\to X(w)$ is
$P_{\Xh(\wh)}$-equivariant and in particular if the generalised
reduced decomposition $\wh$ is good, it is $P_w$-equivariant.

(\i\i)  Let us explain this construction with the quiver $Q_w$ of $w$ and the
  configuration variety. 
Denote by $m(Q_w)$ the set of maximal elements of $Q_w$ for the
  partial order $\preccurlyeq$.

The classical morphism $\pi:\Xt(\wt)\to X(w)$ is given by the
projection from the configuration variety
$\Xt(\wt)\in\prod_{i\in Q_w}G/P_{\b_i}$ on the product $\prod_{i\in
  m(Q_w)}G/P_{\b_i}$.

To give a generalised decomposition $\wt$ of $w$ is equivalent to give
a partition of the vertices of the quiver by subquivers
$(Q_{w_i})_{i\in[1,n]}$. Then the morphism from $\Xt(\wt)$ to
$\Xh(\wh'_i)$ is given by the projection on the product $\prod_{k=1}^i
\prod_{j\in Q_{w_k}}G/P_{\b_j}\times\prod_{k=i+1}^n \prod_{j\in
  m(Q_{w_k})}G/P_{\b_j}$.
\end{rema}

We can give a generalisation of fact \ref{image} for the morphism
$\pit$. Recall that we have an admissible generalised reduced
decomposition $\wh$ given by $w=w_1\cdots w_n$ and a reduced
decomposition of each $w_i$ giving a reduced writing $\wt$ of $w$
given by $s_{\b_1}\cdots s_{\b_r}$. In terms of quiver, the quiver
$Q_w$ of the reduced writing $\wt$ has a partition by subquivers
$Q_{w_i}$ isomorphic to quivers of the elements $w_i$. We denote by
$p_\wh(Q_w)$ the set of vertices of $Q_w$ which are pics for the
quiver $Q_{w_i}$ to which they belong.

\begin{coro}
  The variety $Z_K$ is not contracted by $\pit:\Xt(\wt)\to\Xh(\wh)$ if
  and only if for each $j\in[1,n]$ the part of the subword
  $\displaystyle{\prod_{i\in[1,r]\setminus K}s_i}$ corresponding to a
  subword of $w_j$ is reduced.

In particular, assume that all the $w_i$ are minuscule elements, the
group of divisors of $\Xh(\wh)$ has a basis given by $\pit_*[Z_i]$ for
$i\in p_\wh(Q_w)$ and the group of 1-cycles of $\Xh(\wh)$ has a basis
indexed by $[1,n]$ given by $\pit_*[C_i]$ for $i$ the maximal vertex
of the quiver $Q_{w_i}$ for $i\in[1,n]$.
\end{coro}

\begin{preu}
  This comes from the description fiberwise of the morphism from
  $\Xt(\wt)$ to $\Xh(\wh)$ and the lemma \ref{image}. To verify that
  the described elements form a basis, let us recall that $\Xh(\wh)$ is
  a tower of locally trivial fibration with fibers $X(w_i)$ for
  $i\in[1,n]$ so we know from the case of minuscule Schubert varieties
  that these elements form a basis.
\end{preu}

\subsection{Constructing generalised reduced decomposition}
\label{construction}

In this paragraph, we are going to give a way of constructing good
generalised reduced decomposition of an element $w$ in a product of
minuscule elements $(w_i)_{i\in[1,n]}$. 

\begin{defi}
 Let $A\subset p(Q_w)$ be a subset of the set of pic of $Q_w$, we
denote by $Q(A_w)$ the full subquiver of $Q_w$ containing the vertices
$i$ of $Q_w$ such that $i\not\succcurlyeq j$ for all $j\in p(Q_w)\setminus
A$. 

It is
different from $Q_w$ as soon as $A$ is different from $p(Q_w)$.
\end{defi}

\begin{prop}
(\i) Each connected component $C$ of the quiver $Q_w(A)$ is isomorphic to
  the quiver of a minuscule Schubert variety and in particular has a
  unique maximal element $m(C)$ for the partial order $\preccurlyeq$.

(\i\i) When $A$ has a unique element then $Q_w(A)$ is connected.

(\i\i\i) The quiver $\Qh_w(A)$ obtained from $Q_w$ by removing the
vertices of $Q_w(A)$ is also the quiver of a minuscule Schubert
variety. 

(\i v) The set $p(Q_w(A))$ is $A$ and the set $p(\Qh_w(A))$ is
$p(Q_w)\setminus A$.
\end{prop}

\begin{preu}
(\i) (a) Let us prove that in any connected component $C$ there is only
one maximal element for the partial order $\preccurlyeq$. Let $j_1$
and $j_2$ be two such maximal elements. By connectivity, there exists
a sequence of vertices $i_0=j_1,i_1,\cdots,i_n=j_2$ such that there is
an arrow linking $i_k$ and $i_{k+1}$ (in one sense or another). Let us
take a minimal such sequence (that is to say $n$ is minimal) an let
$x$ be the smallest integer in $[0,n]$ such that $i_{x}\preccurlyeq
i_{x-1}$ and $i_x\preccurlyeq  i_{x+1}$. Such an element exists
because $j_1$ and $j_2$ are maximal. By minimality of $n$ we have
$i_{x-1}\neq i_{x+1}$ are different and thanks to proposition
\ref{geom-carq} the vertex $s(i_x)$ exists. The arrows arriving to
$s(i_x)$ come from $i_{x-1}\in C$, $i_{x+1}\in C$ and maybe from a
third vertex 
$k\preccurlyeq i_x$ (and thus $k\in C$). The vertex $s(i_x)$ has to be
in $C$. If we replace $i_x$ by $s(i_x)$, we get a new sequence of
length $n$ but with $i_{x-1}$ the first term such that $i_{x-1}\preccurlyeq
i_{x-2}$ and $i_{x-1}\preccurlyeq  i_{x}$. By induction we get a
sequence of length $n$ such that the smallest $x\in[0,n]$ with
$i_{x}\preccurlyeq
i_{x-1}$ and $i_x\preccurlyeq  i_{x+1}$ is $x=1$. This tels us that
$s(i_1)\in C$ and $j_1=i_0\preccurlyeq s(i_1)$ in $C$ which is a
contradiction to the maximality of $j_1$.

(b) Let us now prove that any connected component $C$ of $Q(A)$ satisfies
the conditions of \ref{geom-carq}. 
Let $k$ be a vertex of $C$ such that $s(k)$ does not exist or is not
in $C$. 

In the first case, this means that there is at most one arrow from $k$
and denote by $j$ the end vertex of this arrow. If $j$ is not in $C$
(or does not exist) then $k$ is the maximal element of $C$. Otherwise
$j$ is in $C$ and there is exactely one arrow from $k$ in $C$.

In the second case, we have two vertices $k_1$ and $k_2$ such that the
arrows arriving to $s(k)$ come from $k_1$,
$k_2$ and eventually a third one $k_3\preccurlyeq k$ which has to be
in $C$. As $s(k)\not\in C$, at least one of the two vertices $k_1$ and
$k_2$ has to be out of $C$. If both are out of $C$ then $k$ is the
unique maximal element of $C$. Otherwise exactely one vertex from
$\{k_1,k_2\}$ is in $C$.

If $k$ is a vertex of $C$ such that $s(k)\in C$, then we have two
vertices $k_1$ and $k_2$ such that the arrows arriving to $s(k)$ come
from $k_1$, $k_2$ and eventually a third one $k_3\preccurlyeq k$ which
has to be in $C$. These two elements have to be in $C$ otherwise
$s(k)$ would not be in $C$.

(c) We are left to prove that if $m(C)$ is the maximal element of $C$
then $\b(m(C))$ is a simple minuscule root for some group. 

If the Dynkin diagram is of type $A_n$ this is always true because any
simple root is minuscule. Even if the set $\b(C)$ of simple roots is
of type $A_n$ we are done. Let us assume that $\b(C)$ contains a
trivalent root $\ga$ and a root on each branch of the Dynkin diagram
(remark that because $C$ is connected, so is $\b(C)$).

Denote by $m(C)$ the maximal element of $C$ and by $\b$ the simple
root $\b(m(C))$. If this root was not a minuscule root of the Dynkin
diagram $\b(C)$ (a sub-Dynkin diagram of the one of $G$) then we would
have the following situation: 

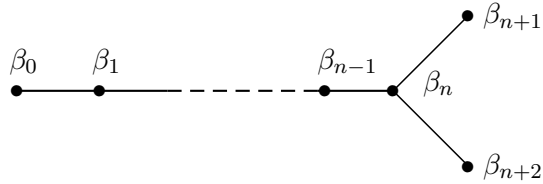
\psset{xunit=1cm}
\psset{yunit=1cm}
\centerline{\hskip 3 cm\begin{pspicture*}(-5.3,-1.5)(5,1.5)
\psline(-5,0)(-3,0)
\psline[linestyle=dashed](-3,0)(-1,0)
\psline(-1,0)(0,0)
\psline[linewidth=0.02](0,0)(1,1)
\psline[linewidth=0.02](0,0)(1,-1)
\put(-5.1,-0.1){$\bullet$}
\put(-5.1,0.3){$\b_0$}
\put(-4,-0.1){$\bullet$}
\put(-4,0.3){$\b_1$}
\put(-1,-0.1){$\bullet$}
\put(-1,0.3){$\b_{n-1}$}
\put(-0.1,-0.1){$\bullet$}
\put(0.4,0){$\b_n$}
\put(0.9,0.9){$\bullet$}
\put(1.2,0.9){$\b_{n+1}$}
\put(0.9,-1.1){$\bullet$}
\put(1.2,-1.1){$\b_{n+2}$}
\end{pspicture*}}
such that $\b_1=\b$, $\b_n=\gamma$ and all the simple roots $\b_i$ are
in $\b(C)$. There are two distinct simple roots $\b_{n+1}$ and
$\b_{n+2}$ in $\b(C)$ and not in $[\b,\ga]$ such that
$\sca{\ga^\vee}{\b_{n+i}}\neq0$ for $i=1,2$. There is a simple root
$\b_0$ in $\b(C)$ different from all the $\b_i$ for $i\in[1,n+2]$ such
that $\sca{\b^\vee}{\b_0}\neq0$.

Denote by $i_k$ the biggest (for $\preccurlyeq$) vertex in $C$ such
that $\b(i_k)=\b_k$ for all $k\in[0,n+2]$. Then because $C$ satisfies
the properties of proposition \ref{geom-carq}, we see that in $C$, for
all $k\in[2,n+1]$ there exists a unique arrow from $i_k$ and it goes
to $i_{k-1}$. In the same way there exists a unique arrow from $i_{n+2}$
and it goes to $i_{n}$ and from $i_0$ to $i_1=m(C)$. This means that
in $C$ we have the following subquiver:

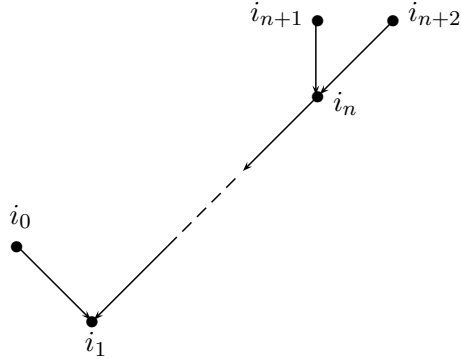
\psset{xunit=1cm}
\psset{yunit=1cm}
\centerline{\begin{pspicture*}(-5,-2.5)(5,2.5)
\psline[linewidth=0.02]{->}(-3,-1)(-2.02,-1.98)
\psline[linestyle=dashed,linewidth=0.02](-0.1,-0.1)(-1,-1)
\psline[linewidth=0.02]{->}(0.97,2)(0.97,1.01)
\psline[linewidth=0.02]{->}(2,2)(1.02,1.02)
\psline[linewidth=0.02]{->}(-1,-1)(-1.98,-1.98)
\psline[linewidth=0.02]{<-}(0,0)(1,1)
\put(-2.1,-2.1){$\bullet$}
\put(-2.1,-2.4){$i_1$}
\put(-3.1,-1.1){$\bullet$}
\put(-3.1,-0.7){$i_0$}
\put(1.9,1.9){$\bullet$}
\put(2.2,2){$i_{n+2}$}
\put(0.9,1.9){$\bullet$}
\put(.1,2){$i_{n+1}$}
\put(0.9,0.9){$\bullet$}
\put(1.2,0.8){$i_n$}
\end{pspicture*}}

Now let us consider the subquiver $Q'$ of $Q$ corresponding to the
vertices $i$ such that $i\succcurlyeq i_0$ or $i\succcurlyeq i_0$ or
$i\succcurlyeq i_0$. This is a quiver corresponding to a minuscule
Schubert variety. Each time there is a hole $i$ in the quiver, we can
add a new vertex $j$ such that $\b(j)=\b(i)$ to obtain a quiver which
still corresponds to a minuscule Schubert variety. We can thus add a
vertex $i_{n+3}$ with $\b(i_{n+3})=\b$ and by induction vertices
$i_{n+2+k}$ with $\b(i_{n+2+k})=\b(i_k)$ for all $k\in[1,n-1]$. In
this new quiver we have the following subquiver:

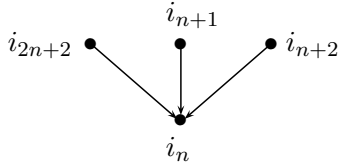
\psset{xunit=1cm}
\psset{yunit=1cm}
\centerline{\begin{pspicture*}(-5,-1)(5,2)
\psline[linewidth=0.02]{<-}(0.05,0)(1.2,1)
\psline[linewidth=0.02]{<-}(0,0.05)(0,1)
\psline[linewidth=0.02]{<-}(-0.05,0)(-1.2,1)
\put(-1.3,0.9){$\bullet$}
\put(-.1,0.9){$\bullet$}
\put(1.1,0.9){$\bullet$}
\put(-0.1,-0.1){$\bullet$}
\put(-0.2,-0.5){$i_n$}
\put(-2.3,0.9){$i_{2n+2}$}
\put(-0.2,1.3){$i_{n+1}$}
\put(1.4,0.9){$i_{n+2}$}
\end{pspicture*}}

But we also could have choosen $\b(i_{2n+2})=\b(i_n)$ proving that in
the quiver $Q_\pd$ of the minuscule homogeneous variety there exists a
vertex $j$ with $s(j)=i$. But then between $j$ and $i_n$ there would
be three vertices (namely $i_{n+1}$, $i_{n+2}$ and $i_{2n+2}$) having
an arrow to $i_n$. This contradicts proposition \ref{geom-carq} for
$Q_\pd$.

(\i\i) If $A$ has a unique element then all vertices of $Q_w(A)$ are
bigger than the unique element of $A$ and $Q_w(A)$ is connected.

(\i\i\i) The quiver $\Qh_w(A)$ is obtained from $Q_w$ by removing all the
vertices smaller than $m(C)$ for any connected component $C$ of
$Q_w(A)$. It is thus (see proposition \ref{description}) the quiver of a
minuscule Schubert variety.

(\i v) Clear from the definitions.
\end{preu}

To construct a partition of the quiver $Q_w$ of a minuscule element
$w$ into quivers $(Q_{w_i})_{i\in[1,n]}$ with $w_i$ a minuscule
elements, it suffices to give a partition $(A_i)_{i\in[1,n]}$ of the
set $p(Q_w)$ of the pics of the quiver. Indeed, given such a partition
$(A_i)_{i\in[1,n]}$, we define by induction a sequence
$(Q_i)_{i\in[0,n]}$ of quivers with $Q_0=Q$ and
$Q_{i+1}=\Qh_i(A_{i+1})$ $Q_{w_1}=Q_w(A_1)$. We then denote by
$Q_{w_i}$ the quiver $Q_{i-1}(A_i)$. The quivers
$(Q_{w_i})_{i\in[1,n]}$ form a partition of $Q_w$ by quivers
associated to minuscule elements $w_i$.

\begin{rema}
\label{sommets}
(\i) For such partitions (giving a reduced generalised decomposition
$\wh$), we have $p_\wh(Q_w)=p(Q_w)$. 

(\i\i)  For such constructions, the vertices of $Q_{w_i}$ are the
  vertices $x$ such that there exists a pic $p\in A_i$ with
  $p\preccurlyeq x$ and $p'\not\preccurlyeq x$ for any pic $p'$ in
  $A_j$ with $j>i$.
\end{rema}

For such partitions $(Q_{w_i})_{i\in[1,n]}$ of the quiver $Q_w$ coming
from partitions $(A_{i})_{i\in[1,n]}$ of $p(Q_w)$ we have a reduced
generalised decomposition $w=w_1\cdots w_n$ denoted $\wh$.

\begin{prop}
  The reduced generalised writing $\wh$ is good.
\end{prop}

\begin{preu}
  We have to prove that the inclusions $\Sigma(P_{w_{i+1}\cdots
  w_n})\subset \Sigma(P^{w_i}\cap G_{w_i})$ and
  $\partial(G_{w_i})\subset\Sigma(P_{w_{i}\cdots w_n})$. But the set
  $\Sigma(P_{w_i\cdots w_n})$ is the set $\b(\{i\in Q_w\ /\ i\
  \textrm{is a hole of } Q_{w_i\cdots w_n} \})$ where $Q_{w_i\cdots
  w_n}$ is the subquiver of $Q_w$ whose vertices are in $\cup_{k\geq
  i}Q_{w_k}$ (see proposition \ref{stab}).

The set $\Sigma(P^{w_i}\cap G_{w_i})$ is the set
$\b(m_i)\cup\supp(w_i)^c$ where $m_i$ is the maximal vertex of
$Q_{w_i}$. So for the first inclusion we only have to prove that for
any simple root $\b\in\supp(w_i)\cap \Sigma(P_{w_{i+1}\cdots w_n})$ we
have $\b=\b(m_i)$. But as $\b\in\supp(w_i)$, there exists $j\in
Q_{w_i}$ such that $\b(j)=\b$. Let $j$ be the biggest such vertex. If
$j$ was not the biggest element $m_i$ in $Q_{w_i}$ then there would
exist in $Q_{w_i}$ an element $k$ with an arrow from $j$ to $k$. But
then there are two cases, if $s(j)$ exists, it is a hole of
$Q_{w_{i+1}\cdots w_n}$ and there are two vertices $k_1$ and $k_2$
having an arrow to $s(j)$ so between $j$ and $s(j)$ there are three
vertices $k$, $k_1$ and $k_2$ and this is impossible thanks to
proposition \ref{geom-carq}. If $s(j)$ does not exist, then $j$ is a
virtual hole of $Q_{w_{i+1}\cdots w_n}$ and there is a vertex $k'$
such that $\sca{\b(k')}{\b(j)}\neq0$ so that $s(j)$ does not exist
but there are two vertices
$k$ and $k'$ having an arrow coming from $j$ and this is impossible
thanks to proposition \ref{geom-carq}.

For the second inclusion, let $\b$ be a simple root in
$\partial(G_{w_i})$, then there exists a vertex $j\in Q_{w_i}$ with
$\sca{\b(j)}{\b}\neq0$. If $\b$ is not in the support of $w_i\cdots
w_n$ then $\b$ is the simple root of a virtual hole and
$\b\in\Sigma(P_{w_i\cdots w_n})$. If $\b$ is in this support then
there exists a vertex $k$ such that $\b(k)=\b$. Let $k$ be the
smallest such vertex, we have an arrow from $j$ to $k$ thus $k$ is not
a pic of $Q_{w_i\cdots w_n}$ and thus not a pic of $Q_{w_{i+1}\cdots
  w_n}$ (see the previous proposition). In particular there exists a
vertex $x\in Q_{w_i\cdots w_n}$ with an arrow from $x$ to $k$. But
then $k$ is the smallest vertex with $\b(k)=\b$ in $Q_{w_i\cdots w_n}$
and there are two arrows arriving to $k$ thus it is a hole of
$Q_{w_i\cdots w_n}$ and we are done.
\end{preu}

We now give here three types of partitions of $Q_w$ constructed in
this way.

\begin{cons}
\label{tous}
  Choose any order $\{i_1,\cdots,i_n\}$ on the set $p(Q_w)$ of the
  pics of $Q_w$ and set $A_k=\{i_k\}$.
\end{cons}

\begin{cons}
\label{lacanonique}
  Define a partition $(A_i)_{i\in[1,n]}$ by induction: $A_1$ is the
  set of pics with minimal height and $A_{i+1}$ is the set of pics in
  $p(Q_w)\setminus\bigcup_{k=1}^iA_k$ with minimal height.
\end{cons}

Before giving the last construction let us fix some notations and
prove the following proposition.
Recall that $p_\wh(Q_w)=p(Q_w)$ (with these constructions) is the set
of all vertices $j$ of $Q_w$ such that there exists an integer
$i\in[1,n]$ with $j\in p(Q_{w_i})$. Let us denote by $m_\wh(Q_w)$ the
set of vertices $j$ of $Q_w$ such that $j$ is a maximal element of
$Q_{w_i}$ for some $i\in[1,n]$.

The partial order $\preccurlyeq$ induces a partial order on
$m_\wt(Q_w)$. Let us finally prove the following:

\begin{prop}
  Let $i\in m_\wh(Q_w)$, there exists a unique minimal element $f(i)$ in
  $m_\wh(Q_w)$ for $\preccurlyeq$ such that $i\preccurlyeq f(i)$. 
\end{prop}

\begin{preu}
  Let us prove the following 

  \begin{lemm}
\label{ordre-max}
    Let $j$ and $k$ in $m_\wh(Q_w)$ such that there exists $x\in Q_w$
    with $x\preccurlyeq j$ and $x\preccurlyeq k$. Then we have either
    $j\preccurlyeq k$ or $k\preccurlyeq j$
  \end{lemm}

  \begin{preu}
    We proceed by induction on $a+b$ where $a$ and $b$ are the indexes
    in $[1,n]$ such that $j\in Q_{w_a}$ and $k\in Q_{w_b}$. Let $x$ be
    a maximal element (for $\preccurlyeq$) such that $x\preccurlyeq j$
    and $x\preccurlyeq k$ and suppose that $x$ is different from $j$
    and $k$.

If there is only one
    arrow from $x$ say going to a vertex $y$, then we must have
    $y\preccurlyeq
    j$ and $y\preccurlyeq k$ contradicting the maximality. Let $y_1$
    and $y_2$ be the two arriving elements of the two arrows from
    $x$. If $\b(y_1)=\b(y_2)$ then $y_1\preccurlyeq y_2$ (or the
    converse) and we have $y_1\preccurlyeq j$ and $y_1\preccurlyeq k$
    contradicting the maximality. We thus have $\b(y_1)\neq\b(y_2)$
    and $y_1\preccurlyeq j$ but $y_1\not\preccurlyeq k$ and
    $y_2\not\preccurlyeq j$ but $y_2\preccurlyeq k$. This also implies
    that $s(x)$ exists because $y_1$ and $y_2$ are connected to the
    biggest element $r$ of the quiver and so the segments
    $[\b(y_1),\b(r)]$ and $[\b(y_1),\b(r)]$ are contained in
    $\b(\{z\in Q_w\ /\ z\succcurlyeq x,\ z\neq x\}$ and thus $\b(x)$
    is in this set. So $s(x)$ exists and we have $s(x)\not\preccurlyeq
    j$ and $s(x)\not\preccurlyeq k$.

Now let $c$, $d$ and $e$ in $[1,n]$ such that $s(x)\in Q_{w_c}$,
$j_1\in Q_{w_d}$ and $j_2\in Q_{w_e}$. We must have $c\geq d$ and $c\geq
e$ because $s(x)\succcurlyeq j_1,j_2$. We must also have $a\geq d$ and
$b\geq e$. But if $p$ is a pic in $A_c$ such that $s(x)\succcurlyeq
p$, we must have $p\preccurlyeq j_1$ or $p\preccurlyeq j_2$ which
implies (see remark \ref{sommets}) that $d\geq c$ or $e\geq c$. We
thus have $c=d$ or $c=e$. Assume for example that $c=d$ and denote by
$m$ the maximal element of $Q_{w_c}$. If $c=d=a$ then
$j_2\preccurlyeq s(x)\preccurlyeq m=j$ and $j_2\preccurlyeq k$ a
contradiction to the maximality of $x$. So $c=d<a$, but we have
$x\preccurlyeq m$ and $x\preccurlyeq j$ and by induction, we must have
$m\preccurlyeq j$. Then we have
$j_2\preccurlyeq s(x)\preccurlyeq m\preccurlyeq j$ and
$j_2\preccurlyeq k$ one more time a contradiction to the maximality of
$x$.
  \end{preu}

The preceding lemma proves that for $i\in m_\wh(Q_w)$ (and even for
any $i\in Q_w$) the set $\{j\in m_\wh(Q_w)\ /\ j\succcurlyeq i\}$ is
totally ordered and thus there exists a minimal element $f(i)$.
\end{preu}

We can now give the last construction which is a particular case of
construction \ref{tous}. Because in construction \ref{tous} there is a
bijection between $p(Q_w)$ and $m_\wh(Q_w)$ we define thanks to the
preceding proposition the function $f$ on the set $p(Q_w)$ simply by the
following: if $p\in p(Q_w)$ is in $Q_{w_i}$ and if $m\in m_\wh(Q_w)$
is the maximal element of $Q_{w_i}$ then $f(m)$ is the maximal element
of some $Q_{w_j}$.
There is only one pic $q$ in $Q_{w_{j}}$ and we define $f(p)=q$.

\begin{cons}
\label{lespetites} 
 Choose an order $\{i_1,\cdots,i_n\}$ on the set $p(Q_w)$ of the
  pics of $Q_w$ such that $h(i_k)\leq h(f(i_{k}))$ for all
  $k\in[1,n-1]$ and set $A_k=\{i_k\}$.

This choice on the order is equivalent to choose an order
$\{i_1,\cdots,i_n\}$ on the set $p(Q_w)$ of the pics of $Q_w$ such
that if $i_k$ and $i_{k+1}$ are adjacent in the quiver then
$h(i_k)\leq h(i_{k+1})$.
\end{cons}

These constructions may produce non connected subquivers $Q_{w_i}$
but thanks to lemma \ref{commute} we may assume (replacing these
quivers by their connected components) that all the quivers $Q_{w_i}$
are quivers of minuscule Schubert varieties.

Construction \ref{lespetites} will give all relative minimal models of
$X(w)$ and construction \ref{lacanonique} will give the relative
canonical model of $X(w)$.

\section{Relative Mori theory of minuscule Schubert varieties}

In this paragraph we prove the above assertion on the relative
canonical and minimal models of a minuscule Schubert variety
$X(w)$. We only consider generalised reduced writing $\wh$ of $w$
obtained thanks one of the three previous constructions.

\subsection{Ample divisors and effective curves}

Recall that we described a basis of divisors and 1-cycles on
$\Xh(\wh)$ in the following way: the group of divisors of $\Xh(\wh)$
has a basis given by $D_i=\pit_*[Z_i]$ for $i\in p_\wh(Q_w)=p(Q_w)$
and the group of 1-cycles of $\Xh(\wh)$ has a basis given by
$\pit_*[C_i]$ for $i\in m_\wh(Q_w)$. Recall also that
(cf. \cite{FS..}) the Chow groups of $\Xh(\wh)$ are generated by
$B$-orbits, free over $\Z$ and the Picard group is dual to the group
of 1-cycles.

We have seen in  remark \ref{rem-proj} that the morphism
$\pit:\Xt(\wt)\to\Xh(\wh)$ is the projection from $\Xt(\wt)$ to the
product $\prod_{i\in m_\wh(Q_w)}G/P_i$. We have on
$\Xh(\wh)$ a projection $p_i:\Xh(\wh)\to G/P_i$ for all $i\in
m_\wh(Q_w)$. Let us define the Cartier divisor
$\M_i=p_i^*(\oo_{G/P_{i}}(1))$ for all $i\in m_\wh(Q_w)$. We have
$\L_i=\pit^*\M_i$ for all $i\in m_\wh(Q_w)$. Because of the
description of $\Xh(\wh)$ as a of tower of locally trivial
fibration with fibers isomorphic to minuscule Schubert varieties
$X(w_i)$, we have the following:

\begin{fact}
  The familly $(\M_i)_{i\in m_\wh(Q_w)}$ is a basis of
  $\pic(\Xh(\wh))$.
\end{fact}

Recall that we gave a basis $(Y_i)_{i\in[1,r]}$ of the cone of
effective 1-cycles $\Xt(w)$ in paragraph \ref{effectif}. Because
$[Y_i]=[C_i]-[C_{s(i)}]$ we see that we have the following:

\begin{fact}
  The familly $(\pit_*[Y_{i}])_{i\in m_\wh(Q_w)}$ is a basis of the
  group of 1-cycles on $\Xh(\wh)$ and is dual to the basis
  $(\M_i)_{i\in m_\wh(Q_w)}$.
\end{fact}

\begin{preu}
  The pull-back of $\M_i$ by $\pit$ is $\L_i$ and we have
  $\L_{i}\cdot[Y_j]=\delta_{i,j}$ so by projection formula we get the
  result.
\end{preu}

\begin{prop}
  The familly $([\pit_*Y_{i}])_{i\in m_\wh(Q_w)}$ is a basis of the cone
  of effective classes and the familly $(\M_i)_{i\in m_\wh(Q_w)}$
  is a basis of the closure of the ample cone.
\end{prop}

\begin{preu}
  The embedding of $\Xh(\wh)$ in $\prod_{i\in m_\wh(Q_w)}G/P_i$ is
  given by $\bigotimes_{i\in m_\wh(Q_w)}\M_i$. The cone generated by
  the $\M_i$ is thus contained in the ample cone.

  Conversely, let $A$ be an ample sheaf and let $a_i=A\cdot[\pit_*Y_{i}]$
  for $i\in m_\wh(Q_w)$. We must have $a_i>0$ and the divisor
  $A-\sum_ia_i\M_i$ is numerically trivial and we get the result.

  By duality we have the result on curves.
\end{preu}

Proposition \ref{ample-div} has the following relative version.

\begin{prop}
\label{ample-div-rel}
We have the formula
$$\M_i=\sum_{k\in
p(Q_{w}),\ k\preccurlyeq i}D_k.$$
\end{prop}

\begin{preu}
Because the variety $\Xh(\wh)$ is normal with rational singularities,
  we have $\M_i=\pit_*\L_i$ and we obtain the formula in the same way
  as in proposition \ref{ample-div} thanks to the fact that all
  quivers $Q_{w_i}$ are associated to minuscule Schubert varieties.
\end{preu}

Let us generalise corollary \ref{goren} and give a criterion for a
variety $\Xh(\wh)$ to be locally factorial.

\begin{coro}
  The variety $\Xh(\wh)$ is locally factorial if and only if for all
  $i\in[1,n]$ the quiver $Q_{w_i}$ has a unique pic.
\end{coro}

\subsection{Canonical divisor of $\Xh(\wh)$}

As in proposition \ref{cano-ample}, we have
$K_{\Xh(\wh)}=\pit_*K_{\Xt(\wt)}$. The same calculus gives the: 

\begin{fact}
  We have
$$-K_{\Xh(\wh)}=\sum_{k\in p(Q_{w})}(h(k)+1)D_k.$$
\end{fact}

For the three constructions, the pics of a fixed quiver $Q_{w_i}$ have
all the same height so we can define $h(w_i)$ to be the height of any
pic of $Q_{w_i}$. Set $h(w_{n+1})=-1$. Proposition \ref{ample-div-rel}
gives us the following (by induction on the order of $m_\wh(Q_w)$):

\begin{coro}
\label{cano-ample-rel}
We have the formula
$$-K_{\Xh(\wh)}=\sum_{i\in m_\wh(Q_w)}(h(w_i)-h(w_{f(i)}))\M_i$$
and in particular $\Xh(\wh)$ is Gorenstein.
\end{coro}

\subsection{Types of singularities}

In this paragraph we are going to prove that the variety $\Xh(\wh)$
has terminal singularities in case of constructions \ref{tous},
\ref{lacanonique} and \ref{lespetites}.

For this we use the resolution $\Xt(\wt)$ of $\Xh(\wh)$ and compare
the canonical divisor $K_{\Xt(\wt)}$ to the pull-back of the canonical
divisor$K_{\Xh(\wh)}$ by $\pit$. We need the following fact coming
directely from the formula of paragraph \ref{canonique} and lemma
\ref{som-haut}:

\begin{fact}
We have the formula 
$$-K_{\Xt(\wt)}=\sum_{i=1}^r(h(i)+1)\xi_i.$$
\end{fact}

We can now prove the following:

\begin{prop}
The variety $\Xh(\wh)$ has terminal (and hence canonical)
singularities.
\end{prop}

\begin{preu}
Let us calculate the pull-back of $K_{\Xh(\wh)}$ by $\pit$: 
$$-\pit^*K_{\Xh(\wh)}=\sum_{i\in m_\wh(Q_w)}(h(w_i)-h(w_{f(i)}))\L_i$$
but thanks to proposition \ref{coord} and the fact that the quivers
$Q_{w_i}$ are the quivers of minuscule Schubert varieties we have
$$\L_{i}=\sum_{k\preccurlyeq i}\xi_k$$
giving
$$-\pit^*K_{\Xh(\wh)}=\sum_{i=1}^n\left((h(w_i)+1)
  \sum_{k\in  Q_{w_i}}\xi_k\right)$$
We get for the difference:
$$K_{\Xt(\wt)}-\pit^*K_{\Xh(\wh)}=\sum_{i=1}^n\sum_{k\in
  Q_{w_i}}(h(w_i)-h(k))\xi_k.$$
But $h(w_i)$ is the highest height of an element in $Q_{w_i}$ so
  $h(w_i)-h(k)\geq0$ with equality for $k$ a pic of the quiver that is
  to say if and only if $\xi_i$ is not contracted by $\pit$.
\end{preu}

\subsection{Description of the relative minimal and canonical models}

We can now prove our results on relative minimal and canonical models
of minuscule Schubert varieties.

\begin{theo}
  (\i) The varieties $\Xh(\wh)$ obtained from construction
  \ref{lespetites} are relative minimal models of $X(w)$.

  (\i\i) The variety $\Xh(\wh)$ obtained from construction
  \ref{lacanonique} is the relative canonical model of $X(w)$.
\end{theo}

\begin{preu}
(\i) We have to prove that any curve $C$ contracted by $\pih:\Xh(\wh)\to
  X(w)$ satisfies $[C]\cdot K_{\Xh(\wh)}\geq0$. The class $[C]$ can be
  writen as $[C]=\sum_ia_i[\pit_*Y_{i}]$ with $a_i\geq0$ and $a_n=0$
  (because the curve is contracted). We just have to prove the non
  negativity of the intersections $K_{\Xh(\wh)}\cdot[\pit_*Y_{j}]$
  for $j\in[1,n-1]$.
We have $K_{\Xh(\wh)}\cdot[\pit_*Y_{j}]=h(w_{f(j)})-h(w_j)$ and by
  construction \ref{lespetites} this intersection is non negative.

(\i\i) It suffices to prove that the contracted curves have a positive
intersection with the canonical divisor and this comes from the
previous calculus and construction \ref{lacanonique}.
  \end{preu}

Let $\Xh(\wh)$ a variety obtained from one of the construction
\ref{tous}, \ref{lacanonique} or \ref{lespetites}. We can completely
describe the extremal rays of the relative cone of effective
1-cycles on $\Xh(\wh)$ (i.e. the effective 1-cycles contracted by
$\pih$).

\begin{fact}
(\i) The extremal rays of $\Xh(\wh)$ are given by the classes
$[\pit_*Y_{j}]$ such that $h(w_{f(j)})< h(w_j)$.

(\i\i) If $\Xh(\wh)$ is obtained from construction \ref{lespetites},
then there is no extremal ray. However, if $D$ is any effective
divisor, then the $(K_{\Xh(\wh)}+D)$-extremal rays are given by the
classes $[\pit_*Y_{j}]$ such that
$(K_{\Xh(\wh)}+D)\cdot[\pit_*Y_{j}]<0$.
\end{fact}

\begin{preu}
(\i) Let $[C]$ the class of an effective curve. Then there exists non
  negative integers $a_i$ such that $[C]=\sum_ia_i[\pit_*Y_i]$. Denote
  by $\mu_j$ (resp. $\nu_k$ and $\omega_l$) the classes $[\pit_*Y_i]$
  such that $K_{\Xh(\wh)}\cdot [\pit_*Y_i]<0$ (resp. $>0$ and
  $=0$). For each $j$ and $k$ there is a linear combination
  with positive coefficient $x_{j,k}\mu_j+y_{j,k}\nu_k$ such that
  $K_{\Xh(\wh)}\cdot(x_{j,k}\mu_j+y_{j,k}\nu_k)=0$. It is easy to
  check that if $K_{\Xh(\wh)}\cdot[C]<0$ then $[C]$ has to be a linear
  combination with non negative coefficient of classes $(\mu_j)$,
  $(x_{j,k}\mu_j+y_{j,k}\nu_k)$ and $(\omega_l)$ proving the result.

(\i\i) The same proof with $K_{\Xh(\wh)}+D$ instead of $K_{\Xh(\wh)}$
works. It works for all varieties $\Xh(\wh)$ obtained from
construction \ref{tous}, \ref{lacanonique} or \ref{lespetites}.
\end{preu}

Let us consider a fixed order $(p_1,\cdots,p_n)$ on the pics of the
quiver $Q_w$ and $\wh$ the good reduced generalised writing it induces
and let us denote by $\wh'$ the good reduced generalised writing
induced by the order $(q_1,\cdots,q_n)$ on the pics where $q_k=p_k$
for $k\not\in\{i,i+1\}$, $q_i=p_{i+1}$ and $q_{i+1}=q_i$. Denote by
$\wh''$ the good reduced generalised writing given by $w''_k=w_k$ for
$k<i$, $w''_i=w_iw_{i+1}$ and $w_k=w_{k+1}$ for $k>i$ that is to say
obtained by the partition $(A_k)_{k\in[1,n-1]}$ of $p(Q_w)$ given by
$A_k=\{p_k\}$ for $k<i$, $A_i=\{p_i,p_{i+1}\}$ and $A_k=\{p_{k+1}\}$ for
$k>i$. Denote by $k_i$ (resp. $k_{i+1}$) the maximal vertex of
$Q_{w_i}$ (resp. $Q_{w_{i+1}}$).

We have two morphisms (for example because of remark \ref{rem-proj})
$$f:\Xh(\wh)\to \Xh(\wh'')\ \ \ {\rm and}\ \ \ f':\Xh(\wh')\to
\Xh(\wh'').$$

\begin{prop} Let $\Xh(\wh)$ obtained from construction \ref{tous},
  \ref{lacanonique} or \ref{lespetites}.

(\i) If $k_i\not\preccurlyeq k_{i+1}$ (i.e. if $i+1\neq f(i)$)
then the morphisms $f$ and $f'$ are isomorphisms.

(\i\i) If $f(i)=i+1$ and $\pit_*[Y_{k_i}]\cdot K_{X}>0$, then $f$
(resp. $f'$) is the small contraction corresponding to the extremal ray
$\mathbb{R}\pit_*[Y_{k_i}]$ (resp. $\mathbb{R}\pit_*[Y_{k'_i}]$) and
$f$ is the flip of $f'$.

\vs 0.2 cm

(\i\i\i)  If $f(i)=i+1$ and $\pit_*[Y_{k_i}]\cdot K_{X}=0$, denote
$D=D_{i+1}$ (resp. $D'=D'_{i}$) then for $\varepsilon>0$, the class
$\pit_*[Y_{k_i}]$ (resp. $\pit_*[Y_{k'_i}]$) is extremal for
$K_{\Xh(\wh)}+\varepsilon D$ (resp. $K_{\Xh(\wh')}+\varepsilon D'$)
and $f$ (resp. $f'$) is the small contraction corresponding to the
extremal ray $\mathbb{R}\pit_*[Y_{k_i}]$
(resp. $\mathbb{R}\pit_*[Y_{k'_i}]$). The morphism $(f',D')$ is the
flop of $(f,D)$.
\end{prop}

\begin{preu}
(\i) Because $k_i$ and $k_{i+1}$ are non comparable for
$\preccurlyeq$, we have thanks to lemma \ref{ordre-max} that $w_i$ and
$w_{i+1}$ satisfy the hypothesis of lemma \ref{commute} and we have
the result.

In the cases where $f(i)=i+1$, we already know
that $\Xh(\wh)$ and $\Xh(\wh')$ are locally factorial with terminal
singularities. We also know that $\Xh(\wh'')$ is normal and the
morphisms are birational and mori-small (the group of Weil divisors
has a basis given by the pics). Because of our description of Picard
groups we also have $\rho(X/X'')=\rho(X'/X'')=1$.

(\i\i) We are left to study the divisors $K_{\Xh(\wh)}$ and
$K_{\Xh(\wh')}$ on the fibers of $f$ and $f'$. We will not describe
these fibers in details in this proposition but more details will be
given in the next paragraph. Because all the sheaves $\M_{k_j}$ for $j\neq
i$ are already defined on $\Xh(\wh'')$ they are trivial on the fibers
of $f$ and the sheaf $\M_{k_i}$ on $\Xh(\wh)$ is relatively ample with
respect to $f$. The restriction of $K_{\Xh(\wh)}$ to the fibers of $f$
is given by $(K_{\Xh(\wh)}\cdot \pit_*[Y_{k_i}])\M_{k_i}$ so that
$K_{\Xh(\wh)}$ is ample, anti-ample or trivial according to the
positivity of the intersection $K_{\Xh(\wh)}\cdot \pit_*[Y_{k_i}]$ and thus
according to the height of the pics. This proves in case (\i\i) that
$K_{\Xh(\wh)}$ is $f$-ample.

Furthermore, the fiber of $f$ is contained in $G/P_{\b(k_i)}$ thus the
classes of contracted curves are proportional to $\pit_*[Y_{k_i}]$.

In the same way we get that $-K_{\Xh(\wh')}$ is $f'$-ample and the
result.

(\i\i) Proposition \ref{ample-div-rel} tel us that $D=D_{i+1}$
  satisfies $\pit_*[Y_{k_i}]\cdot D<0$ so the class $\pit_*[Y_{k_i}]$
  is extremal for $K_{\Xh(\wh)}+\varepsilon D$.

The Bott-Samelson variety $\Xt(\wt)$ is a resolution of the
birational morphism between $\Xh(\wh)$ and $\Xh(\wh')$ and $D$ is the
image of $\xi_{p_{i+1}}$ whose image in $\Xh(\wh')$ is $D'$ so that
$D'$ is the strict transform of $D$.

We are left to study the divisors $K_{\Xh(\wh)}$, $K_{\Xh(\wh')}$, $D$
  and $D'$ on the fibers of the morphisms $f$ and $f'$. The calculus
  in (\i) proves the triviality of $K_{\Xh(\wh)}$ and $K_{\Xh(\wh')}$
  on the fibers. For $D$ and $D'$, the same argument as for the
  canonical sheaves proves that their restriction to the fibers of $f$
  (resp. $f'$) is a positive multiple of $-\M_{k_i}$
  (resp. $\M_{k'_i}$) concluding the proof.
\end{preu}

\begin{rema}
  If $f(i)=i+1$ and $[Y_{k_i}]\cdot K_{\Xh(\wh)}<0$ then by symetry
  $f'$ is the flip of $f$.
\end{rema}

\begin{coro}
(\i) The varieties $\Xh(\wh)$ obtained from construction \ref{tous}
are linked by flips and flops and any variety obtained from $\Xh(\wh)$
by flips and flops comes from this construction.

(\i\i) The varieties $\Xh(\wh)$ obtained from construction
\ref{lespetites} are linked by flops and any variety obtained from
$\Xh(\wh)$ by flops comes from this construction.
\end{coro}

\begin{preu}
(\i) We know that the extremal rays (or more generally the
  $(K_{\Xh(\wh)}+D)$-extremal rays) of $\Xh(\wh)$ are generated by
  the classes $\pit_*[Y_{k_i}]$ with
  $K_{\Xh(\wh)}\cdot\pit_*[Y_{k_i}]<0$
  (resp. $(K_{\Xh(\wh)}+D)\cdot\pit_*[Y_{k_i}]<0$). But the
  associated flip or flop gives a variety obtained by construction
  \ref{tous}.

(\i\i) The same argument works in this case because
$h(w_{f(i)})=h(w_i)$ and we stay in the class of varieties obtained
from construction \ref{lespetites}.
\end{preu}

\begin{rema}
  For any variety obtained from construction \ref{tous} we proved
  existence and terminaison of flips and flops (in the sens of
  K. Matsuki \cite{Matsuki}).
\end{rema}

  \begin{coro}
\label{fin-mori}
The relative minimal models of $X(w)$ are exactely the varieties
    obtained from construction \ref{lespetites}.
  \end{coro}

  \begin{preu}
We use theorem 12-1-8 of \cite{Matsuki}, the fact that
    varieties obtained from construction \ref{lespetites} are relative
    minimal models and existence and termination of flops for these
    varieties.
  \end{preu}

\section{Small $IH$-resolutions of minuscule Schubert varities}

In this section we prove that the morphism $\Xh(\wh)\to X(w)$ obtained
from construction \ref{lespetites} is $IH$-small. We then discuss
the smoothness of $\Xh(\wh)$ and describe all $IH$-small resolutions
of minuscule Schubert varieties. Let us first recall the definition of
an $IH$-small morphism:

\begin{defi}
\label{small-def}
  A morphism $\pi:Y\to X$ is said to be $IH$-small if for al $k>0$, we
  have
$${\rm codim}_X\{x\in X\ /\ \dim(\pi^{-1}(x))=k\}>2k.$$

A small morphism $\pi:Y\to X$ is a small resolution of $X$ if $Y$ is
smooth.
\end{defi}

In this section we will use a case by case analysis instead of a
global proof for all minuscule Schubert varieties in the same time. A
combinatorical direct proof on the quiver is possible but it would lead
to a too complicated combinatorical discussion and we prefer avoiding
it. 

\subsection{Necessary condition}

Let us first prove the following proposition showing that among the
morphisms $\pi:\Xh(\wh)\to X(w)$ obtained from construction
\ref{tous} only the one coming from construction \ref{lespetites} can
be small:

\begin{prop}
Let $\Xh(\wh)$ obtained from construction \ref{tous} but not from
construction \ref{lespetites}. Then the morphism $\pih:\Xh(\wh)\to
X(w)$ is not $IH$-small.
\end{prop}

\begin{preu}
If $\pih:\Xh(\wh)\to X(w)$ was small then, because $\Xh(\wh)$ has
terminal singularities, it would be a relative minimal model
(this is a consequence of the proof by B. Totaro \cite{Totaro} of
theorem \ref{totaro}). This is not the case by corollary
\ref{fin-mori}.

One can give an explicit subvariety in $X(w)$ not satisfying the
  $IH$-small condition for $\pih$. Let $i$ be a pic of $Q_w$ such that
  $h(f(i))<h(i)$ (such a pic exists because the resolution is not
  obtained from contruction \ref{lespetites}). Let us consider the
  smallest (for $\preccurlyeq$) 
  vertex $j\in Q_w$ such that $j\succcurlyeq i$ and $j\succcurlyeq
  f(i)$. Then one can prove that the image $\pit(Z_k)$ of the divisor
  $Z_j\subset\Xt(\wt)$ in $\Xh(\wh)$ is of codimension
  $h(f(i))-h(j)+1$ and that its image $\pi(Z_j)$ in $X(w)$ is of
  codimension $h(i)-h(j)+h(f(i))-h(j)+1$. The fiber above $\pi(Z_j)$
  contains $\pit(Z_j)$ and is of dimension at least $h(i)-h(j)$. But
  we have
$$\codim_{X(w)}(\pi(Z_j))=h(i)-h(j)+h(f(i))-h(j)+1\leq2(h(i)-h(j)).$$
\end{preu}

On the contrary when we choose a good order on the pics then we obtain
the following

\begin{theo}
  \label{petit}
The morphisms $\pih:\Xh(\wh)\to X(w)$ obtained from construction
\ref{lespetites} are $IH$-small.
\end{theo}

In subsections \ref{Fibers}, \ref{A_netD_n} and \ref{casexcep} we are
going to give a proof of this theorem. Let us first give an easy
corollary:

\begin{coro}
  The morphism $\pih:\Xh(\wh)\to X(w)$ obtained from construction
\ref{lacanonique} is $IH$-small.
\end{coro}

\begin{preu}
  Indeed, any morphism $\pih'$ obtained from construction
  \ref{lespetites} factors through the morphism $\pih$ and it is
  easy to verify that this implies that, as $\pih'$ is $IH$-small,
  $\pih$ is $IH$-small. 
\end{preu}

\subsection{Fibers}
\label{Fibers}

We will adapt the technics of \cite{SanVan} in our setting. The idea
that choosing a good order in the pics will produce $IH$-small
morphisms comes from A. Zelevinsky's paper \cite{Zelevinsky}.

Let us recall that the varieties $\Xh(\wh)$ where constructed by
induction. The first step beeing given by the morphism
$p:\overline{PuQ}\times^QX(v)\to X(w)$ where $P$ is the stabiliser of
$X(w)$, $v$ is obtained from $w$ by removing the first pic and $u$
corresponds to the removed vertices. The group $Q$ is the intersection
of $P$ with the stabiliser of $X(v)$. By induction there exists a
resolution $\pi':\Xh(\vh)\to X(v)$ equivariant under the stabiliser
of $X(v)$. The resolution is given by the fiber product
$\pih:\Xh(\wh)=\overline{PuQ}\times^Q\Xh(\vh)\to X(w)$.

Let us prove a lemma giving a description of the fibers of $\pi$
and a formula on their dimension. This lemma is directely inspired by
lemma 2.1 of \cite{SanVan}. If $w'\leq w$ in the Burhat order, let us
denote by $U(w')$ the $P$-orbit ($P$ is the stabiliser of $X(w)$) of
$e_{w'}$ (the fixed point of the torus corresponding to the Schubert cell
of $w'$) in $X(w)$. Because $\pi$ is $P$-equivariant, all the fibers of
points in $U(e_{w'})$ are isomorphic and to calculate $f_{\pi,w'}$ the
dimension of the fiber $\pi^{-1}(e_{w'})$, it is enough to calculate
$\dim(\pi^{-1}(U(e_{w'})))-\dim(U(e_{w'}))$.

\begin{lemm}
\label{fibres}
  Define the set 
$$S(w',w)=\left\{(u',v')\in W\ \bigg/
  \begin{array}{cc}
\ u'\leq u\ \textrm{\textit{and}}\
  v'\leq v,\ \textrm{\textit{in the Bruhat order and}}\\
PX(u'v')=X(w'), PX(u')=X(u')\ \textrm{\textit{and}}\
  QX(v')=X(v')  \end{array}\right\}.$$

(\i) We have $\displaystyle{p^{-1}(U({w'}))=\bigcup_{(u',v')\in
  S(w',w)}Pu'Q\times^QQe_{v'}.}$ 

(\i\i) We have $\displaystyle{\pi^{-1}(U({w'}))=\bigcup_{(u',v')\in
  S(w',w)}Pu'Q\times^Q\pi'^{-1}(Qe_{v'}).}$ 

(\i\i\i) This gives the formula
$$f_{\pi,w'}=\card(Q_{u'})+f_{\pi',v'}+\card(Q_{v'})-\card(Q_{w'})=
\card(Q_{u'})+f_{\pi',v'}-\codim_{X(w')}(X(v'))$$ 
for some $(u',v')\in S(w',w)$.
\end{lemm}

\begin{preu}
(\i) Let $(u',v')\in S(w',w)$, we have the inclusions:
$$p(Pu'Q\times^QQe_{v'})\subset
Pe_{u'}Qe_{v'}\subset Pe_{u'}X({v'})\subset
PX(u'v')=X(w').$$
Furthermore $p(Pu'Q\times^QQe_{v'})$ is a $P$-orbit
an contains $Pe_{u'}e_{v'}=Pe_{w'}$ so $Pu'Q\times^QQe_{v'}$ is
contained in $p^{-1}(U(w'))$. 

Conversely, if $(x,y)\in\overline{PuQ}\times^QX(v)$ is such that
$p(x,y)=xy\in U(w')$, then there are elements $u'$ and $v'$ in the
Weyl group such that we have $PX(u')=X(u')$, $QX(v')=X(v')$ and
$\overline{PxQ}\times^Q\overline{Qy P_\pd}/P_\pd=
\overline{Pu'Q}\times^QX(v')$. But then there exists $(p,q)\in P\times
Q$ such that $pe_{u'}=x$ and $qe_{v'}=y$ so that we have $(x,y)\in
Pu'Q\times^QQe_{v'}$. Furthermore, there exists $p'\in P$ such that
$p'xy=e_{w'}$ thus $p'pe_{u'}qe_{v'}=e_{w'}$. This implies, because
$QX(v')=X(v')$, that $PX(u'v')=X(w')$.

(\i\i) Comes directely from (\i).

(\i\i\i) Because $\dim U(w')=\card(Q_{w'})$, we only need to prove
that 
$$\dim(\pi^{-1}(U(w')))=\card(Q_{u'})+f_{\pi',v'}+\card(Q_{v'})$$
for some $(u',v')\in S(w',w)$, it is true thanks to (\i\i).
\end{preu}

\begin{rema}
\label{u'v'}
(\i) In the case where $\overline{PuQ}/Q$ is an homogeneous variety, we
  recover lemma 2.1 of \cite{SanVan}. In this case we must have
  $u'=u$ because $X(u)$ is the only $P$-stable Schubert subvariety of
  $X(u)$. We then have
  $\card(Q_{u'})=\codim_{X(w)}(X(v))=\card(Q_w)-\card(Q_v)$.

(\i\i) More generaly, the Schubert variety $X(u')$ is a Schubert
subvariety of $X(u)$ with the same stabiliser so that if $i$ is a hole
of its quiver then $\b(i)\in\b(t(Q_u))$. The Schubert variety
$X(v')$ is a Schubert subvariety of $X(v)$ with stabiliser
$\stab(w)\cap\stab(v)$. If $i$ is a hole of $Q_{v'}$ then $\b(i)$ has
to be in the union $\b(t(Q_v))\cup\b(t(Q_w))$.

Let us now describe the condition $PX(u'v')=X(w')$. The quiver of the
Schubert variety $X(u'v')$ is obtained by gluing the quiver of $X(u')$
above the quiver of $X(v')$. Furthermore, the quiver of the Schubert
variety $PX(a)$ if $X(a)$ is a
Schubert subvariety of $X(w)$, is the smallest subquiver $Q$ of $Q_w$
containing the quiver $Q_a$ and such that
$\b(t(Q))\subset\b(t(Q_w))$. In particular, if we denote by $A$ the
set of vertices $i\in Q_a$ such that $i$ is not the succesor of an
element of $Q_a$ and $\b(i)\in\b(t(Q_w))$ then the set non virtual holes
of $PX(a)$ is $A$ and the virtual holes are associated to simple roots
in $\b(t(Q_w))\setminus\b(A)$.

If $i$ is a hole of $Q_{w'}$ such that $\b(i)$ is not in the support
of $u$. Let $j$ be the smallest vertex in $Q_{v'}$ such that
$\b(j)=\b(i)$. Then $j$ will be a vertex of $Q_{v'u'}$ with no
predecessor and has to be a hole of $PX(u'v')=X(w')$. We must thus
have $j=i$. In particular all the holes of $Q_{w'}$ associated to
simple roots not in the support of $u$ have to be holes of $v'$.  
\end{rema}

The proof of the theorem \ref{petit} will go as follows. Because the
morphism $\pi$ is $P$-equivariant where $P$ is the stabilisor of
$X(w)$, we need to prove that for any $w'\in W$ such that $X(w')$ is
stable under $P$ in $X(w)$, we have
$\codim_{X(w)}(X(w'))>2f_{\pi,w'}$. 

For the classical cases ($A_n$ and $D_n$) we will introduce two
functions $\Gamma$ and $q$ such that 
$$\codim_{X(w)}(X(w'))=\Gamma(w',w)+q(w',w)$$
The function $q$ will take only non negative values and will be
positive if $w'\neq w$.

We then proceed by induction on the number of pics of $w$ (or on the
number of fibrations in $\Xh(\wh)$) and prove the more stronger
result:
$$\Gamma(w',w)\geq 2f_{\pi,w'}.$$
Because of the previous lemma, it is enough to prove that for all
$(u',v')\in S(w',w)$ we have
$$\Gamma(w',w)\geq
2(\card(Q_{u'})+f_{\pi',v'}-\codim_{X(w')}(X(v'))).$$
Let $\theta\in W$ such that $X(\theta)$ is the closure of the orbit of
$X(v')$ in $X(v)$ under $\stab(X(v))$. We have
$f_{\pi',v'}=f_{\pi',\theta}$ and by induction hypothesis we have
$2f_{\pi',\theta}\leq\Gamma(\theta,v)$. We are thus reduced to prove:
$$2(\card(Q_{u'})-\codim_{X(w')}(X(v')))\leq\Gamma(w',w)-
\Gamma(\theta,v).$$
We prove this formula in the following paragraph in the $A_n$ and
$D_n$ case.

\subsection{The case of $A_n$ and $D_n$}
\label{A_netD_n}

To prove the result on smallness, we will need to describe the
elements of $S(w',w)$ and calculated the dimension in the formula of
lemma \ref{fibres}. The only two difficult cases of minuscule Schubert
varieties will be the cases of grassmannians (the varieties constructed
are the varieties of Zelevinsky \cite{Zelevinsky}) and of maximal
isotropic subspaces in an even dimensional vector space endowed with a
non degenerate quadratic form (some of these cases have been treated
in \cite{SanVan} and we complete their study). Indeed the other
minuscule Schubert varieties are those of $E_6$ and $E_7$ for which we
will make a case by case analysis (which is tedious by not hard with
lemma \ref{fibres}) and the case of Schubert varieties in a quadric
which are very simple. 

\subsubsection{The $A_n$ case}

Let us consider a minuscule quiver $Q_w$ of a Schubert variety in the
grassmannian $\G(p,q)$ of $p$-dimensional subvector spaces of a
$q$-dimensional vector space. We may assume that all the simple roots
are in the support of $w$ (otherwise we simply restrict the group) so
that there will be no virtual hole in $Q_w$. The set of simple roots
$\b(t(Q_w))$ can 
be writen as $\{\a_{k_1},\cdots,\a_{k_s}\}$ in the notation of
\cite{bourb}. Let us denote by $t_1,\cdots,t_s$ the holes such that
$\b(t_i)=\a_{k_i}$. Because of proposition \ref{geom-carq}, for all
$i\in[2,s]$ there exists exactely one pic $p_i$ between the holes
$t_{i-1}$ and $t_i$. Furthermore there must be a pic $p_1$
(resp. $p_{s+1}$) with $\b(p_1)=\a_{k}$ (resp
$\b(p_{s+1})=\a_{k}$) with $k< k_1$ (resp. $k> k_s$). In
particular we see that the number $n$ of pics equals $s+1$.

Let us now define the following sequences $(a_i(w))_{i\in[1,s+1]}$ and
$(b_i(w))_{i\in[0,s]}$ of integers (we will sometimes simply denote
them by $a_i$ and $b_i$ omiting $w$): 
$$
\left\{\begin{array}{l}
a_i(w)=h(p_i)-h(t_i) \textrm{ and } b_i(w)=h(p_{i+1})-h(t_i) \textrm{ for
} i\in[1,s],\\
a_{s+1}(w)=p-\sum_{i=1}^s a_i(w),\\
b_{0}(w)=q-p-\sum_{i=1}^sb_i(w).
\end{array}\right.
$$

It is an easy game on the quiver and the description of configuration
varieties to verify that if the sequences of integers associated to
the quiver $Q_w$ are $(a_i(w))_{i\in[1,s+1]}$ and $(b_i(w))_{i\in[0,s]}$
then we have:
$$X(w)=\left\{V\in\G(p,q)\ /\ \dim(V\cap\k^{n_i})\geq m_i \textrm{ for
    all }i\in[1,s]\right\}$$
where $n_i=\displaystyle{\sum_{k=1}^i(a_{k}+b_{k-1})}$ and
    $m_i=\displaystyle{\sum_{k=1}^ia_{k}}$.

Let $X(w')$ be a Schubert subvariety of $X(w)$ with the same
stabiliser. Then we must have $\b(t(Q_{w'}))\subset\b(t(Q_w))$. For
any hole $t_i$ of $Q_w$ let us define the depth of $w'$ in $t_i$ to be
the integer
$$c_i=\card\{j\in Q_w\setminus Q_{w'}\ /\ \b(j)=\b(t_i)\}.$$
The same game on the quiver and the description of configuration
varieties shows that the associated sequences are given by
$$\left\{
  \begin{array}{l}
a_i(w')=a_i(w)+c_i-c_{i-1}\ \textrm{for all } i\in[1,s+1]\\
 b_i(w')=b_i(w)+c_{i}-c_{i+1}\ \textrm{for all } i\in[0,s]
  \end{array}\right.$$
with $c_0=c_{s+1}=0$. We have
$$X(w)=\left\{V\in\G(p,q)\ /\ \dim(V\cap\k^{n_i})\geq m_i+l_i \textrm{ for
    all }i\in[1,s]\right\}$$
with $l_i={\sum_{k=1}^ic_k}$. These
    description enables us to give the following fact to calculate
    the codimension of $X(w')$ in $X(w)$ (we set $c_0=c_{s+1}=0$):

    \begin{fact}
      We have the formula:
$$\codim_{X(w)}(X(w'))=\card(Q_w)-\card(Q_{w'})=\Gamma(w',w)+q(w',w)$$
where 
$$\Gamma(w',w)=\sum_{i=1}^sc_i(a_i+b_i)\ \ \textrm{   and   }\ \ 
q(w',w)=\frac{1}{2}\sum_{i=1}^{s+1}(c_i-c_{i-1})^2.$$
    \end{fact}

    \begin{preu}
      This can be seen with a simple calculation. We will describe it
      geometrically on the quiver. We have the following quiver (see
      the appendix for a description of the quivers):

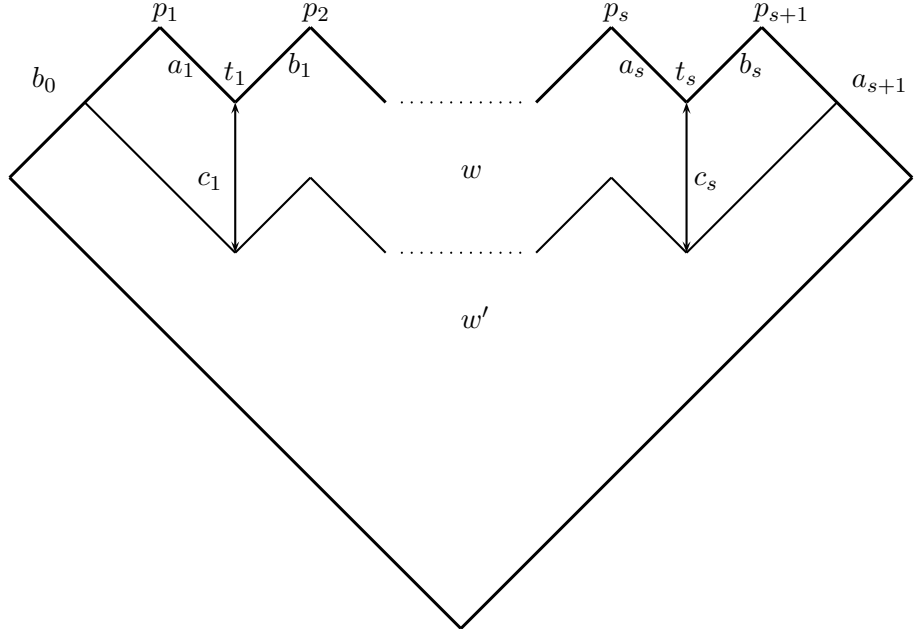
\psset{xunit=1cm}
\psset{yunit=1cm}
\centerline{\begin{pspicture*}(-9,-7)(9,2)
\psline[linewidth=0.04](-6,-1)(-4,1)
\psline[linewidth=0.04](-3,0)(-4,1)
\psline[linewidth=0.04](-3,0)(-2,1)
\psline[linewidth=0.04](-1,0)(-2,1)
\psline[linewidth=0.03,linestyle=dotted](-0.8,0)(0.8,0)
\psline[linewidth=0.04](6,-1)(4,1)
\psline[linewidth=0.04](3,0)(4,1)
\psline[linewidth=0.04](3,0)(2,1)
\psline[linewidth=0.04](1,0)(2,1)
\psline[linewidth=0.04](-6,-1)(0,-7)
\psline[linewidth=0.04](0,-7)(6,-1)
\psline(-5,0)(-3,-2)
\psline(-3,-2)(-2,-1)
\psline(-2,-1)(-1,-2)
\psline[linewidth=0.03,linestyle=dotted](-0.8,-2)(0.8,-2)
\psline(5,0)(3,-2)
\psline(3,-2)(2,-1)
\psline(2,-1)(1,-2)
\psline{<->}(-3,0)(-3,-2)
\put(-3.5,-1.1){$c_1$}
\psline{<->}(3,0)(3,-2)
\put(3.1,-1.1){$c_s$}
\put(-2.1,1.15){$p_2$}
\put(1.9,1.15){$p_s$}
\put(-4.1,1.15){$p_1$}
\put(3.9,1.15){$p_{s+1}$}
\put(-3.15,0.3){$t_1$}
\put(2.85,0.3){$t_s$}
\put(-3.9,0.4){$a_1$}
\put(2.1,0.4){$a_{s}$}
\put(-2.3,0.4){$b_1$}
\put(3.7,0.4){$b_{s}$}
\put(-5.7,0.2){$b_0$}
\put(5.2,0.2){$a_{s+1}$}
\put(0,-1){$w$}
\put(0,-3){$w'$}
\end{pspicture*}}

The codimension beeing given by the difference of the number of
vertices, we find
$$\codim_{X(w)}(X(w'))=\sum_{i=1}^sc_i(a_i+b_i)-
\sum_{i=1}^{s-1}c_ic_{i+1}+\sum_{i=1}^{s}c_i^2$$ 
and simple calculation gives the formula. Remark that we have
$q(w',w)>0$ for $w'\neq w$.
    \end{preu}

Let us assume that $v$ is obtained from $w$ by removing the $k^{th}$
pic of $Q_w$. We obtain on the quivers the following situation:

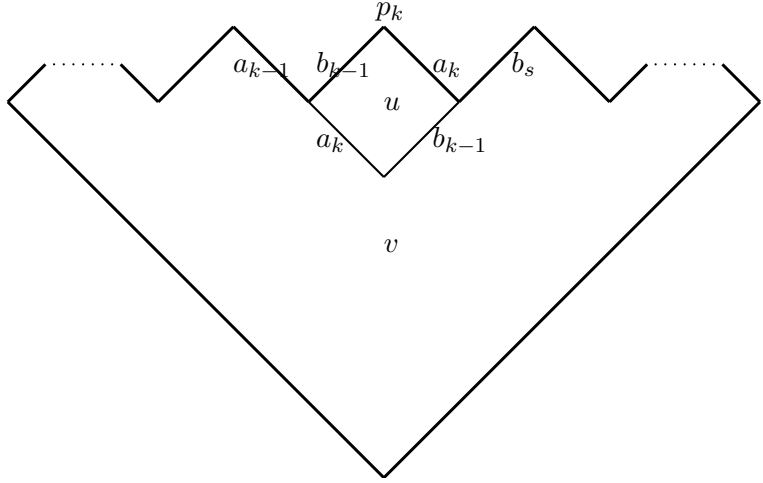
\psset{xunit=1cm}
\psset{yunit=1cm}
\centerline{\begin{pspicture*}(-9,-5)(9,2)
\psline[linewidth=0.04](-3,0)(-3.5,0.5)
\psline[linewidth=0.04](-4.5,0.5)(-5,0)
\psline[linewidth=0.04](-3,0)(-2,1)
\psline[linewidth=0.04](-1,0)(-2,1)
\psline[linewidth=0.03,linestyle=dotted](-4.4,0.5)(-3.6,0.5)
\psline[linewidth=0.04](3,0)(2,1)
\psline[linewidth=0.04](1,0)(2,1)
\psline[linewidth=0.04](3,0)(3.5,0.5)
\psline[linewidth=0.04](4.5,0.5)(5,0)
\psline[linewidth=0.04](-1,0)(0,1)
\psline[linewidth=0.04](1,0)(0,1)
\psline[linewidth=0.04](-5,0)(0,-5)
\psline[linewidth=0.04](5,0)(0,-5)
\psline(-1,0)(0,-1)
\psline(1,0)(0,-1)
\psline[linewidth=0.03,linestyle=dotted](4.4,0.5)(3.6,0.5)
\put(-0.1,1.15){$p_k$}
\put(-2,0.4){$a_{k-1}$}
\put(-0.9,0.4){$b_{k-1}$}
\put(-0.9,-0.6){$a_{k}$}
\put(0.65,0.4){$a_{k}$}
\put(0.65,-0.6){$b_{k-1}$}
\put(1.7,0.4){$b_{s}$}
\put(0,-0.1){$u$}
\put(0,-2){$v$}
\end{pspicture*}}

This simply means that the sequences of integers
$(a_i(v))_{i\in[1,s]}$ and $(b_i(v))_{i\in[0,s-1]}$ are given by 
$$a_i(v)=\left\{\begin{array}{cl}
a_i(w) & \textrm{ for } i<k-1\\
a_k(w)+a_{k-1}(w) & \textrm{for } i=k-1\\
a_{i+1}(w) & \textrm{ for } i\geq k
\end{array}\right.\ \textrm{ and }\ b_i(v)=\left\{\begin{array}{cl}
b_i(w) & \textrm{ for } i<k-1\\
b_k(w)+b_{k-1}(w) & \textrm{for } i=k-1\\
b_{i+1}(w) & \textrm{ for } i\geq k
\end{array}\right.$$
where $a_0(w)=b_0(w)$. Furthermore, the quiver $Q_u$ has no hole
meaning that the variety $\overline{PuQ}/Q$ is smooth and $u'$ has to
be equal to $u$. We only need to determine $v'$.

Let us now consider the quiver $Q$ obtained by intersecting in $Q_w$
the quivers $Q_v$ and $Q_{w'}$. The quivers of $v'$ has to be a
subquiver of this quiver such that (see remark \ref{u'v'}) all the
holes of $Q_{w'}$ are holes of $Q_{v'}$ and $Q_{v'}$ may have one more hole
corresponding to the hole of $v$ which is not a hole of $w$. 

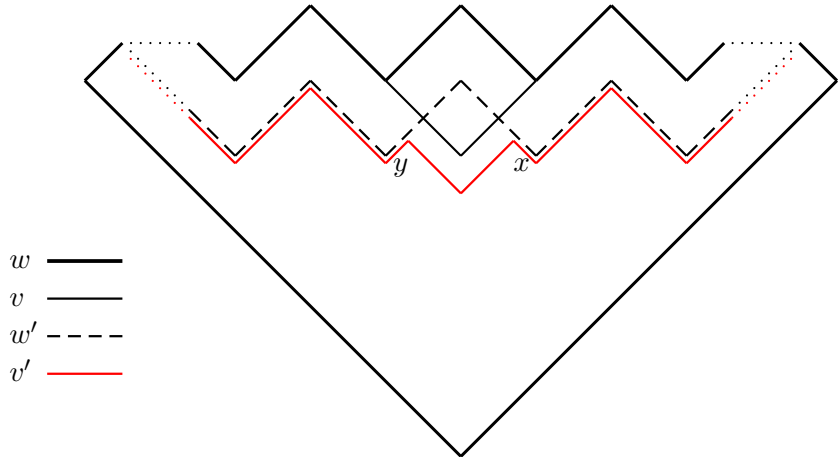
\psset{xunit=1cm}
\psset{yunit=1cm}
\centerline{\begin{pspicture*}(-9,-5)(9,2)
\psline[linewidth=0.04](-3,0)(-3.5,0.5)
\psline[linewidth=0.04](-4.5,0.5)(-5,0)
\psline[linewidth=0.04](-3,0)(-2,1)
\psline[linewidth=0.04](-1,0)(-2,1)
\psline[linewidth=0.03,linestyle=dotted](-4.4,0.5)(-3.6,0.5)
\psline[linewidth=0.04](3,0)(2,1)
\psline[linewidth=0.04](1,0)(2,1)
\psline[linewidth=0.04](3,0)(3.5,0.5)
\psline[linewidth=0.04](4.5,0.5)(5,0)
\psline[linewidth=0.04](-1,0)(0,1)
\psline[linewidth=0.04](1,0)(0,1)
\psline[linewidth=0.04](-5,0)(0,-5)
\psline[linewidth=0.04](5,0)(0,-5)
\psline[linewidth=0.03,linestyle=dotted](4.4,0.5)(3.6,0.5)
\psline(-1,0)(0,-1)
\psline(1,0)(0,-1)
\psline[linewidth=0.03,linestyle=dotted](4.4,0.4)(3.6,-0.4)
\psline[linewidth=0.03,linestyle=dotted](-4.4,0.4)(-3.6,-0.4)
\psline[linewidth=0.03,linestyle=dashed](-3,-1)(-3.6,-0.4)
\psline[linewidth=0.03,linestyle=dashed](3,-1)(3.6,-0.4)
\psline[linewidth=0.03,linestyle=dashed](-3,-1)(-2,0)
\psline[linewidth=0.03,linestyle=dashed](3,-1)(2,0)
\psline[linewidth=0.03,linestyle=dashed](-1,-1)(-2,0)
\psline[linewidth=0.03,linestyle=dashed](1,-1)(2,0)
\psline[linewidth=0.03,linestyle=dashed](-1,-1)(0,0)
\psline[linewidth=0.03,linestyle=dashed](1,-1)(0,0)
\put(-6,-2.5){$w$}
\put(-6,-3){$v$}
\put(-6,-3.5){$w'$}
\put(-6,-4){$v'$}
\psline[linewidth=0.05](-5.5,-2.4)(-4.5,-2.4)
\psline[linewidth=0.03,linestyle=dashed](-5.5,-3.4)(-4.5,-3.4)
\psline[linewidth=0.03](-5.5,-2.9)(-4.5,-2.9)
\psline[linewidth=0.03,linecolor=red](-5.5,-3.9)(-4.5,-3.9)
\psline[linewidth=0.03,linestyle=dotted,linecolor=red](4.4,0.3)(3.6,-0.5)
\psline[linewidth=0.03,linestyle=dotted,linecolor=red](-4.4,0.3)(-3.6,-0.5)
\psline[linewidth=0.03,linecolor=red](-3,-1.1)(-3.6,-0.5)
\psline[linewidth=0.03,linecolor=red](3,-1.1)(3.6,-0.5)
\psline[linewidth=0.03,linecolor=red](-3,-1.1)(-2,-0.1)
\psline[linewidth=0.03,linecolor=red](3,-1.1)(2,-0.1)
\psline[linewidth=0.03,linecolor=red](-1,-1.1)(-2,-0.1)
\psline[linewidth=0.03,linecolor=red](1,-1.1)(2,-0.1)
\psline[linewidth=0.03,linecolor=red](-1,-1.1)(-0.7,-0.8)
\psline[linewidth=0.03,linecolor=red](1,-1.1)(0.7,-0.8)
\psline[linewidth=0.03,linecolor=red](0,-1.5)(-0.7,-0.8)
\psline[linewidth=0.03,linecolor=red](0,-1.5)(0.7,-0.8)
\put(0.7,-1.2){$x$}
\put(-0.9,-1.2){$y$}
\end{pspicture*}}

The
quiver $Q_{v'}$ has $s+1$ holes and the sequences
$(a_i(v'))_{i\in[1,s+2]}$ and $(b_i(v'))_{i\in[0,s+1]}$ are given by
$$a_i(v')=\left\{\begin{array}{cl}
a_i(w') & \textrm{ for } i\leq k-1\\
a_k(w')-x & \textrm{for } i=k\\
x & \textrm{for } i=k+1\\
a_{i-1}(w') & \textrm{ for } i > k+1
\end{array}\right.\ \textrm{ and }\ b_i(v')=\left\{\begin{array}{cl}
b_i(w') & \textrm{ for } i<k-1\\
y & \textrm{for } i=k-1\\
b_{k-1}(w')-y & \textrm{for } i=k\\
b_{i-1}(w) & \textrm{ for } i\geq k+1
\end{array}\right.$$
where $x\in[0,c_{k}]$, $y\in[0,c_{k-1}]$ and
$c_{k-1}-y=c_k-x$. Indeed, the last formula is given by the fact
that the only hole different from those of $Q_{w'}$ has to be
associated to the same root as the hole of $v$ which is not a hole of
$w$. This gives the equality
$$\sum_{i=1}^{k-1}(a_i(v)+b_{i-1}(v))=\sum_{i=1}^{k}(a_i(v')+b_{i-1}(v'))$$
and the equality $c_{k-1}-y=c_k-x$. The fact that
$x\in[0,c_{k}]$ and $y\in[0,c_{k-1}]$ are equivalent to the fact that
$X(v')$ is a Schubert subvariety of $X(v)$.

The Schubert subvariety $X(\theta)$ is contained in $X(v)$, contains
$X(v')$ and is stable by the stabiliser of $X(v)$. It must have the
same hole as $v'$ except for those not corresponding to holes of
$v$. In our case we have to fill the holes $k-1$ and $k+1$ of $v'$ to
obtain $\theta$:

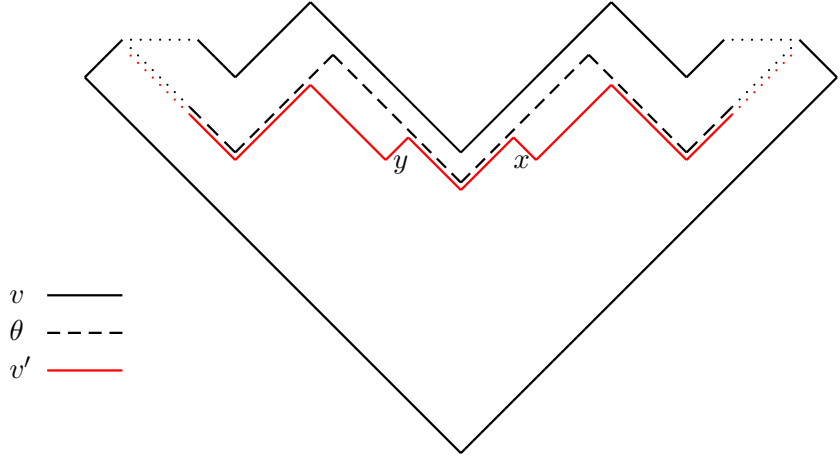
\psset{xunit=1cm}
\psset{yunit=1cm}
\centerline{\begin{pspicture*}(-9,-5)(9,2)
\psline[linewidth=0.03](-3,0)(-3.5,0.5)
\psline[linewidth=0.03](-4.5,0.5)(-5,0)
\psline[linewidth=0.03](-3,0)(-2,1)
\psline[linewidth=0.03](-1,0)(-2,1)
\psline[linewidth=0.03,linestyle=dotted](-4.4,0.5)(-3.6,0.5)
\psline[linewidth=0.03](3,0)(2,1)
\psline[linewidth=0.03](1,0)(2,1)
\psline[linewidth=0.03](3,0)(3.5,0.5)
\psline[linewidth=0.03](4.5,0.5)(5,0)
\psline[linewidth=0.03](-5,0)(0,-5)
\psline[linewidth=0.03](5,0)(0,-5)
\psline[linewidth=0.03,linestyle=dotted](4.4,0.5)(3.6,0.5)
\psline[linewidth=0.03,linestyle=dotted](4.4,0.4)(3.6,-0.4)
\psline[linewidth=0.03,linestyle=dotted](-4.4,0.4)(-3.6,-0.4)
\psline[linewidth=0.03,linestyle=dashed](-3,-1)(-3.6,-0.4)
\psline[linewidth=0.03,linestyle=dashed](3,-1)(3.6,-0.4)
\psline[linewidth=0.03,linestyle=dashed](-3,-1)(-1.7,0.3)
\psline[linewidth=0.03,linestyle=dashed](3,-1)(1.7,0.3)
\psline[linewidth=0.03,linestyle=dashed](0,-1.4)(-1.7,0.3)
\psline[linewidth=0.03,linestyle=dashed](0,-1.4)(1.7,0.3)
\psline(-1,0)(0,-1)
\psline(1,0)(0,-1)
\put(-6,-3){$v$}
\put(-6,-3.5){$\theta$}
\put(-6,-4){$v'$}
\psline[linewidth=0.03,linestyle=dashed](-5.5,-3.4)(-4.5,-3.4)
\psline[linewidth=0.03](-5.5,-2.9)(-4.5,-2.9)
\psline[linewidth=0.03,linecolor=red](-5.5,-3.9)(-4.5,-3.9)
\psline[linewidth=0.03,linestyle=dotted,linecolor=red](4.4,0.3)(3.6,-0.5)
\psline[linewidth=0.03,linestyle=dotted,linecolor=red](-4.4,0.3)(-3.6,-0.5)
\psline[linewidth=0.03,linecolor=red](-3,-1.1)(-3.6,-0.5)
\psline[linewidth=0.03,linecolor=red](3,-1.1)(3.6,-0.5)
\psline[linewidth=0.03,linecolor=red](-3,-1.1)(-2,-0.1)
\psline[linewidth=0.03,linecolor=red](3,-1.1)(2,-0.1)
\psline[linewidth=0.03,linecolor=red](-1,-1.1)(-2,-0.1)
\psline[linewidth=0.03,linecolor=red](1,-1.1)(2,-0.1)
\psline[linewidth=0.03,linecolor=red](-1,-1.1)(-0.7,-0.8)
\psline[linewidth=0.03,linecolor=red](1,-1.1)(0.7,-0.8)
\psline[linewidth=0.03,linecolor=red](0,-1.5)(-0.7,-0.8)
\psline[linewidth=0.03,linecolor=red](0,-1.5)(0.7,-0.8)
\put(0.7,-1.2){$x$}
\put(-0.9,-1.2){$y$}
\end{pspicture*}}

The quiver $Q_\theta$ of $\theta$ has $s-1$ holes and the integers
$(a_i(\theta))_{i\in[1,s]}$ and $(b_i(\theta))_{i\in[0,s-1]}$ are
given by:
$$a_i(\theta)=\left\{\begin{array}{cl}
a_i(v') & \textrm{ for } i\leq k-2\\
a_{k-1}(v')+a_k(v') & \textrm{for } i=k-1\\
a_{k+1}(v')+a_{k+2}(v') & \textrm{for } i=k\\
a_{i+2}(v') & \textrm{ for } i \geq k+1
\end{array}\right.\ \textrm{ and }\ b_i(\theta)=\left\{\begin{array}{cl}
b_i(v') & \textrm{ for } i<k-2\\
b_{k-2}(v')+b_{k-1}(v') & \textrm{for } i=k-2\\
b_{k}(v')+b_{k+1}(v') & \textrm{for } i=k-1\\
b_{i+2}(w) & \textrm{ for } i\geq k
\end{array}\right.$$
It is now an easy calculation (and straightforward on the quiver) that
the depth $c'_i$ of $\theta$ in the holes of $v$ is $c_{k-1}-y=c_k-x$ for the
$(k-1)^{\rm th}$, and $c_i$ for all the holes before the
$(k-1)^{\rm th}$ hole and $c_{i+1}$ for all the holes after the
$(k-1)^{\rm th}$ hole. We can now calculate:
$$\codim_{X(w)}(X(v))-\codim_{X(w')}(X(v'))=
a_k(w)b_{k-1}(w)-(a_{k}(w')-x)(b_{k-1}(w')-y)$$
and finaly:
$$\codim_{X(w)}(X(v))-\codim_{X(w')}(X(v'))=xa_k(w)+yb_{k-1}(w)-xy.$$
On the other hand we calculate
$$\Gamma(w',w)-\Gamma(\theta,v)=
\sum_{i=1}^sc_i(a_i(w)+b_i(w))-\sum_{i=1}^{s-1}c'_i(a_i(v)+b_i(v))$$ 
$$=\sum_{i=1}^sc_i(a_i(w)+b_i(w))- \sum_{i<k-1}c_i(a_i(w)+b_i(w))-
\sum_{i>k-1}c_{i+1}(a_{i+1}(w)+b_{i+1}(w))-c'_{k-1}(a_{k-1}(v)+b_{k-1}(v))$$
$$=c_{k-1}(a_{k-1}(w)+b_{k-1}(w))+c_k(a_k(w)+b_k(w))-
(c_k-x)(a_k(w)+b_{k}(w))- (c_{k-1}-y)(a_{k-1}(w)+b_{k-1}(w))$$
and finaly
$$\Gamma(w',w)-\Gamma(\theta,v)=x(a_k(w)+b_k(w))+y(a_{k-1}(w)+b_{k-1}(w)).$$
Now the fact that the $k^{\rm th}$ pic of $w$ was smaller than the
adjacent pics means that we have $a_{k-1}(w)\geq b_{k-1}(w)$ and
$b_k(w)\geq a_k(w)$ so that we get the inequality (because $x$ and $y$
are non negative):
$$2(\codim_{X(w)}(X(v))-\codim_{X(w')}(X(v')))\leq
\Gamma(w',w)-\Gamma(\theta,v).$$
And the theorem is proved in this case.

\subsubsection{The $D_n$ case}

Let us consider a minuscule quiver $Q_w$ of a Schubert variety in the
grassmannian $\G_{iso}(p,2p)$ of $p$-dimensional isotropic subvector
spaces of a $2p$-dimensional vector space endowed with a non
degenerate quadratic form. We may assume that all the simple roots
are in the support of $w$ (otherwise we simply restrict the group) so
that there will be no virtual hole in $Q_w$. The set of simple roots
$\b(t(Q_w))$ can 
be writen as $\{\a_{k_1},\cdots,\a_{k_s}\}$ in the notation of
\cite{bourb}. Let us denote by $t_1,\cdots,t_s$ the holes such that
$\b(t_i)=\a_{k_i}$. Because of proposition \ref{geom-carq}, for all
$i\in[2,s]$ there exists exactely one pic $p_i$ between the holes
$t_{i-1}$ and $t_i$. Furthermore there must be a pic $p_1$
with $\b(p_1)=\a_{k}$ with $k< k_1$. If
$\a_{k_s}\not\in\{\a_{p-1},\a_p\}$ then there must be pic $p_{s+1}$
with $\b(p_{s+1})=\a_{p-1}$ or $\a_p$ and in this case there are $s+1$
pics, otherwise there are $s$ pics. In the last case, we have
$\a_{k_s}=\a_{p-1+i}$ with $i=0$ or $i=1$. We define $p_{s+1}$ to be
the smallest vertex (for $\preccurlyeq$) of $Q_w$ with
$\b(p_{s+1})=\a_{p-i}$. 

Let us now define the following sequences $(a_i(w))_{i\in[1,s+1]}$ and
$(b_i(w))_{i\in[0,s]}$ of integers (we will sometimes simply denote
them by $a_i$ and $b_i$ omiting $w$): 
$$\textrm{In the first case:}\ 
\left\{\begin{array}{l}
a_i(w)=h(p_i)-h(t_i) \textrm{ for } i\in[1,s]\\
b_i(w)=h(p_{i+1})-h(t_i) \textrm{ for } i\in[1,s-1],\\
b_s(w)=h(p_{s+1})-h(t_s)+\frac{1}{2}\\ 
a_{s+1}(w)=p-\sum_{i=1}^s a_i(w),\\
b_{0}(w)=\frac{1}{2}h(p_{s+1})-\sum_{i=1}^sb_i(w).
\end{array}\right.
$$
$$\textrm{In the second case:}\ 
\left\{\begin{array}{l}
a_i(w)=h(p_i)-h(t_i) \textrm{ and }
b_i(w)=h(p_{i+1})-h(t_i) \textrm{ for } i\in[1,s-1],\\
a_s(w)=h(p_s)-h(p_{s+1}) \textrm{ and }
b_s(w)=h(p_{s+1})-h(p_{s+1})-\frac{1}{2}\\ 
a_{s+1}(w)=p-\sum_{i=1}^s a_i(w),\\
b_{0}(w)=\frac{1}{2}h(p_{s+1})-\sum_{i=1}^sb_i(w).
\end{array}\right.
$$

It is an easy game on the quiver and the description of configuration
varieties to verify that if the associated sequences of integers to
the quiver $Q_w$ are $(a_i(w))_{i\in[1,s+1]}$ and $(b_i(w))_{i\in[0,s]}$
then we have:
$$X(w)=\left\{V\in\G_{iso}(p,2p)\ /\ \dim(V\cap \k^{n_i})\geq m_i \textrm{ for
    all }i\in[1,s]\right\}$$
where the $\k^k$ form a complete flag of fixed isotropic subspaces,
    $n_i=
{\sum_{k=1}^i(a_{k}+b_{k-1})}$ and
    $m_i=
{\sum_{k=1}^ia_{k}}$ for $i\in[1,s]$ in the first case, and $n_i=
{\sum_{k=1}^i(a_{k}+b_{k-1})}$ and
    $m_i=
{\sum_{k=1}^ia_{k}}$ for $i\in[1,s-1]$, $n_s=p$,
$m_s=1+{\sum_{k=1}^sa_{k}}$ and $\k^{n_s}$ is associated to the root
$\a_{k_s}$ in the second case.

Let $X(w')$ be a Schubert subvariety of $X(w)$ with the same
stabiliser. Then we must have $\b(t(Q_{w'}))\subset\b(t(Q_w))$. For
any hole $t_i$ of $Q_w$ let us define the depth of $w'$ in $t_i$ to be
the integer
$$c_i=\left\{
  \begin{array}{l}
\card\{j\in Q_w\setminus Q_{w'}\ /\ \b(j)=\b(t_i)\} \textrm{ for }
\b(t_i)\not\in\{\a_{p-1},\a_p\}\\ 
2\card\{j\in Q_w\setminus Q_{w'}\ /\ \b(j)=\b(t_i)\} \textrm{ for }
\b(t_i)\in\{\a_{p-1},\a_p\} \end{array}\right..$$
The same  game on the quiver and the description of configuration
varieties shows that the associated sequences are given by
$$\left\{
  \begin{array}{l}
a_i(w')=a_i(w)+c_i-c_{i-1}\ \textrm{for all } i\in[1,s+1]\\
 b_i(w')=b_i(w)+c_{i}-c_{i+1}\ \textrm{for all } i\in[0,s]
  \end{array}\right.$$
with 
$c_0=0$ and $c_{s+1}=c_s$. 
We have
$$X(w)=\left\{V\in\G(p,q)\ /\ \dim(V\cap\k^{n_i})\geq m_i+l_i \textrm{ for
    all }i\in[1,s]\right\}$$
with $l_i=
{\sum_{k=1}^ic_k}$. These
    description enables us to give the following fact to calculate
    the codimension of $X(w')$ in $X(w)$ (we set $c_0=c_{s+1}=0$):

    \begin{fact}
      We have the formula:
$$\codim_{X(w)}(X(w'))=\card(Q_w)-\card(Q_{w'})=\Gamma(w',w)+q(w',w)$$
where 
$$\Gamma(w',w)=\sum_{i=1}^sc_i(a_i+b_i)\ \ \textrm{   and   }\ \ 
q(w',w)=\frac{1}{2}\sum_{i=1}^{s}(c_i-c_{i-1})^2.$$
    \end{fact}

    \begin{preu}
      This can be seen with a simple calculation. We will describe it
      geometrically on the quiver. In the first case, we have the
      following quiver:

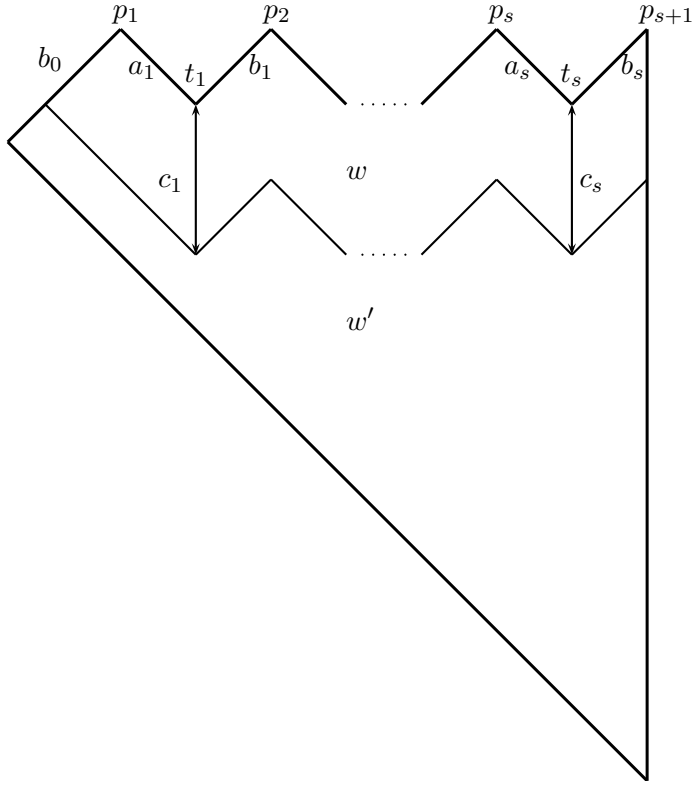
\psset{xunit=1cm}
\psset{yunit=1cm}
\centerline{\begin{pspicture*}(-10,-9)(9,1.5)
\psline[linewidth=0.04](-5.5,-0.5)(-4,1)
\psline[linewidth=0.04](-3,0)(-4,1)
\psline[linewidth=0.04](-3,0)(-2,1)
\psline[linewidth=0.04](-1,0)(-2,1)
\psline[linewidth=0.03,linestyle=dotted](-0.8,0)(-0.2,0)
\psline[linewidth=0.04](3,-9)(3,1)
\psline[linewidth=0.04](2,0)(3,1)
\psline[linewidth=0.04](2,0)(1,1)
\psline[linewidth=0.04](0,0)(1,1)
\psline[linewidth=0.04](-5.5,-0.5)(3,-9)
\psline(-5,0)(-3,-2)
\psline(-3,-2)(-2,-1)
\psline(-2,-1)(-1,-2)
\psline[linewidth=0.03,linestyle=dotted](-0.8,-2)(-0.2,-2)
\psline(3,-1)(2,-2)
\psline(2,-2)(1,-1)
\psline(1,-1)(0,-2)
\psline{<->}(-3,0)(-3,-2)
\put(-3.5,-1.1){$c_1$}
\psline{<->}(2,0)(2,-2)
\put(2.1,-1.1){$c_s$}
\put(-2.1,1.15){$p_2$}
\put(0.9,1.15){$p_s$}
\put(-4.1,1.15){$p_1$}
\put(2.9,1.15){$p_{s+1}$}
\put(-3.15,0.3){$t_1$}
\put(1.85,0.3){$t_s$}
\put(-3.9,0.4){$a_1$}
\put(1.1,0.4){$a_{s}$}
\put(-2.3,0.4){$b_1$}
\put(2.65,0.4){$b_{s}$}
\put(-5.1,0.5){$b_0$}
\put(-1,-1){$w$}
\put(-1,-3){$w'$}
\end{pspicture*}}

The codimension beeing given by the difference of the number of
vertices, we find
$$\codim_{X(w)}(X(w'))=\sum_{i=1}^{s-1}c_i(a_i+b_i)+
c_s(a_s+b_s-\frac{1}{2})-\frac{c_s(c_s-1)}{2}
-\sum_{i=1}^{s-1}c_ic_{i+1}+\sum_{i=1}^{s}c_i^2$$ 
and simple calculation gives the formula. In the second case, we have
the following quiver:

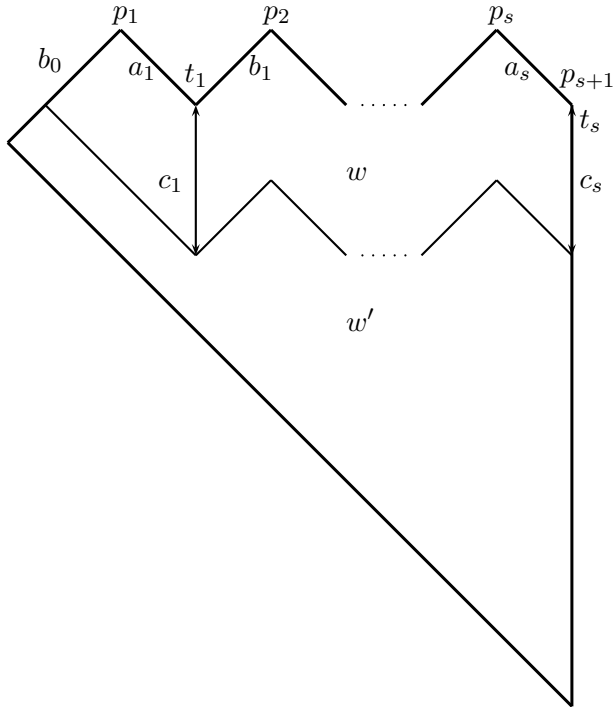
\psset{xunit=1cm}
\psset{yunit=1cm}
\centerline{\begin{pspicture*}(-10,-8)(9,1.5)
\psline[linewidth=0.04](-5.5,-0.5)(-4,1)
\psline[linewidth=0.04](-3,0)(-4,1)
\psline[linewidth=0.04](-3,0)(-2,1)
\psline[linewidth=0.04](-1,0)(-2,1)
\psline[linewidth=0.03,linestyle=dotted](-0.8,0)(-0.2,0)
\psline[linewidth=0.04](2,-8)(2,0)
\psline[linewidth=0.04](2,0)(1,1)
\psline[linewidth=0.04](0,0)(1,1)
\psline[linewidth=0.04](-5.5,-0.5)(2,-8)
\psline(-5,0)(-3,-2)
\psline(-3,-2)(-2,-1)
\psline(-2,-1)(-1,-2)
\psline[linewidth=0.03,linestyle=dotted](-0.8,-2)(-0.2,-2)
\psline(2,-2)(1,-1)
\psline(1,-1)(0,-2)
\psline{<->}(-3,0)(-3,-2)
\put(-3.5,-1.1){$c_1$}
\psline{<->}(2,0)(2,-2)
\put(2.1,-1.1){$c_s$}
\put(-2.1,1.15){$p_2$}
\put(0.9,1.15){$p_s$}
\put(-4.1,1.15){$p_1$}
\put(2.1,-0.3){$t_s$}
\put(-3.15,0.3){$t_1$}
\put(1.85,0.3){$p_{s+1}$}
\put(-3.9,0.4){$a_1$}
\put(1.1,0.4){$a_{s}$}
\put(-2.3,0.4){$b_1$}
\put(-5.1,0.5){$b_0$}
\put(-1,-1){$w$}
\put(-1,-3){$w'$}
\end{pspicture*}}

The codimension beeing given by the difference of the number of
vertices, we find
$$\codim_{X(w)}(X(w'))=\sum_{i=1}^{s-1}c_i(a_i+b_i)+a_sc_s-
\sum_{i=1}^{s-1}c_ic_{i+1}+\sum_{i=1}^{s-1}c_i^2+\frac{c_s(c_s+1)}{2}$$ 
and simple calculation gives the formula. Remark that we have
$q(w',w)>0$ for $w'\neq w$.
    \end{preu}

Let us assume that $v$ is obtained from $w$ by removing the $k^{th}$
pic of $Q_w$. We obtain on the quivers the four following situation:

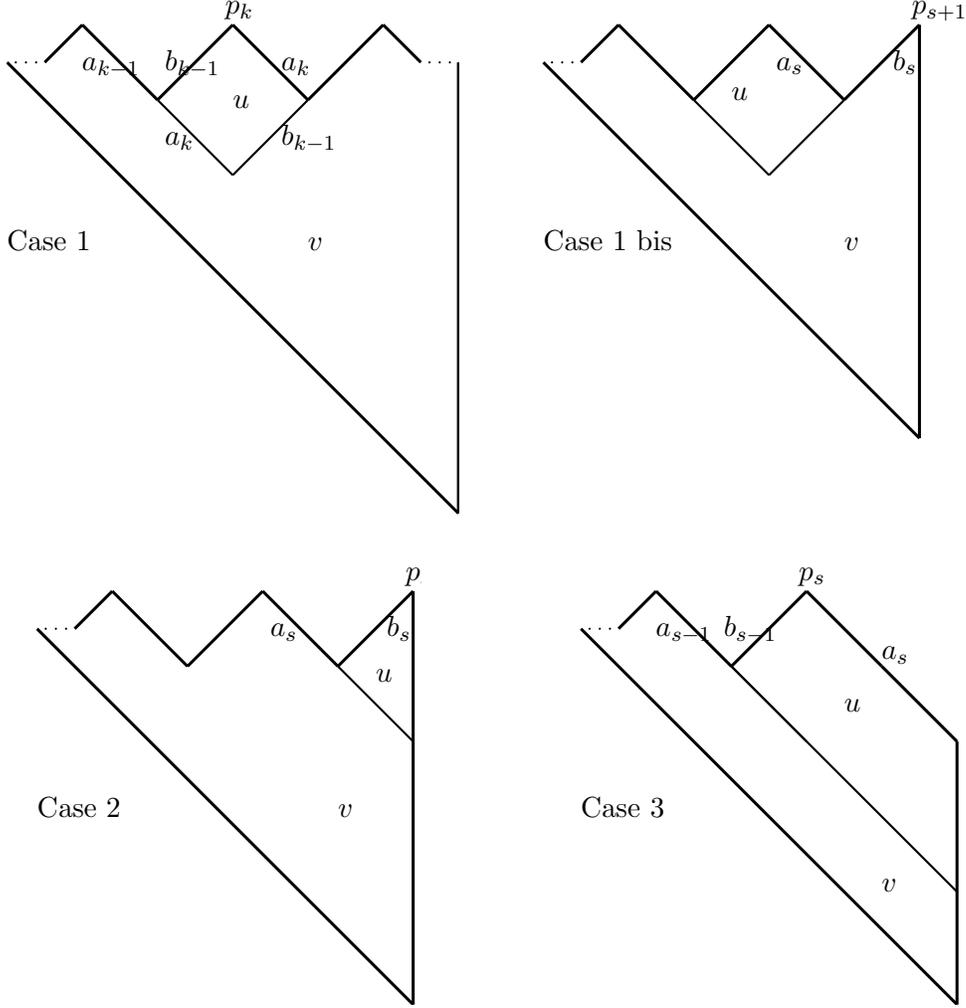
\psset{xunit=1cm}
\psset{yunit=1cm}
\centerline{\begin{pspicture*}(-3,-6)(3,1.5)
\psline[linewidth=0.04](-2.5,0.5)(-2,1)
\psline[linewidth=0.04](-1,0)(-2,1)
\psline[linewidth=0.03,linestyle=dotted](-2.9,0.5)(-2.6,0.5)
\psline[linewidth=0.04](2.5,0.5)(2,1)
\psline[linewidth=0.04](1,0)(2,1)
\psline[linewidth=0.04](-1,0)(0,1)
\psline[linewidth=0.04](1,0)(0,1)
\psline[linewidth=0.04](-3,0.5)(3,-5.5)
\psline[linewidth=0.04](3,0.5)(3,-5.5)
\psline(-1,0)(0,-1)
\psline(1,0)(0,-1)
\psline[linewidth=0.03,linestyle=dotted](2.6,0.5)(2.9,0.5)
\put(-3,-2){Case 1}
\put(-0.1,1.15){$p_k$}
\put(-2,0.4){$a_{k-1}$}
\put(-0.9,0.4){$b_{k-1}$}
\put(-0.9,-0.6){$a_{k}$}
\put(0.65,0.4){$a_{k}$}
\put(0.65,-0.6){$b_{k-1}$}
\put(0,-0.1){$u$}
\put(1,-2){$v$}
\end{pspicture*}
\psset{xunit=1cm}
\psset{yunit=1cm}
\begin{pspicture*}(-4,-6)(3,1.5)
\psline[linewidth=0.04](-2.5,0.5)(-2,1)
\psline[linewidth=0.04](-1,0)(-2,1)
\psline[linewidth=0.03,linestyle=dotted](-2.9,0.5)(-2.6,0.5)
\psline[linewidth=0.04](1,0)(2,1)
\psline[linewidth=0.04](-1,0)(0,1)
\psline[linewidth=0.04](1,0)(0,1)
\psline[linewidth=0.04](-3,0.5)(2,-4.5)
\psline[linewidth=0.04](2,1)(2,-4.5)
\psline(1,0)(0,-1)
\psline(-1,0)(0,-1)
\put(-3,-2){Case 1 bis}
\put(1.9,1.15){$p_{s+1}$}
\put(0.1,0.4){$a_{s}$}
\put(1.65,0.4){$b_{s}$}
\put(-.5,0){$u$}
\put(1,-2){$v$}
\end{pspicture*}}

\psset{xunit=1cm}
\psset{yunit=1cm}
\centerline{\begin{pspicture*}(-3,-4.5)(2.1,1.5)
\psline[linewidth=0.04](-2.5,0.5)(-2,1)
\psline[linewidth=0.04](-1,0)(-2,1)
\psline[linewidth=0.03,linestyle=dotted](-2.9,0.5)(-2.6,0.5)
\psline[linewidth=0.04](1,0)(2,1)
\psline[linewidth=0.04](-1,0)(0,1)
\psline[linewidth=0.04](1,0)(0,1)
\psline[linewidth=0.04](-3,0.5)(2,-4.5)
\psline[linewidth=0.04](2,1)(2,-4.5)
\psline(1,0)(2,-1)
\put(-3,-2){Case 2}
\put(1.9,1.15){$p_{s+1}$}
\put(0.1,0.4){$a_{s}$}
\put(1.65,0.4){$b_{s}$}
\put(1.5,-0.2){$u$}
\put(1,-2){$v$}
\end{pspicture*}
\psset{xunit=1cm}
\psset{yunit=1cm}
\begin{pspicture*}(-5,-4.5)(2.1,1.5)
\psline[linewidth=0.04](-2.5,0.5)(-2,1)
\psline[linewidth=0.04](-1,0)(-2,1)
\psline[linewidth=0.03,linestyle=dotted](-2.9,0.5)(-2.6,0.5)
\psline[linewidth=0.04](1,0)(2,-1)
\psline[linewidth=0.04](-1,0)(0,1)
\psline[linewidth=0.04](1,0)(0,1)
\psline[linewidth=0.04](-3,0.5)(2,-4.5)
\psline[linewidth=0.04](2,-1)(2,-4.5)
\psline(-1,0)(0,-1)
\psline(-1,0)(2,-3)
\put(-3,-2){Case 3}
\put(-0.1,1.15){$p_s$}
\put(-2,0.4){$a_{s-1}$}
\put(-1.1,0.4){$b_{s-1}$}
\put(1,0.1){$a_{s}$}
\put(0.5,-0.6){$u$}
\put(1,-3){$v$}
\end{pspicture*}}

Case 1 is strictly equivalent to the $A_n$ case and we get the result
with the same calculation in this case. Case 1 bis has been done
with these techniques in \cite{SanVan}, we will not make the
calculation one more time. For case 2, the sequences of
integers $(a_i(v))_{i\in[1,s+1]}$ and $(b_i(v))_{i\in[0,s]}$ are given
by:
$$a_i(v)=\left\{\begin{array}{cl}
a_i(w) & \textrm{ for } i<s\\
a_s(w)+b_s(w)-\frac{1}{2} & \textrm{for } i=s\\
a_{s+1}(w)-(b_s(w)-\frac{1}{2}) & \textrm{ for } i=s+1
\end{array}\right.\ \textrm{and}\ b_i(v)=\left\{\begin{array}{cl}
b_{0}(w)+b_s(w)-\frac{1}{2} & \textrm{ for } i=0\\
b_i(w) & \textrm{ for } 0<i<s\\
\frac{1}{2} & \textrm{for } i=s.\\
\end{array}\right.$$
For case 3, the sequences of integers $(a_i(v))_{i\in[1,s]}$ and
$(b_i(v))_{i\in[0,s-1]}$ are given by:
$$a_i(v)=\left\{\begin{array}{cl}
a_i(w) & \textrm{ for } i<s-1\\
a_s(w)+a_{s-1}(w)+b_{s-1}(w) & \textrm{for } i=s-1\\
a_{s+1}(w)-b_{s-1}(w) & \textrm{ for } i=s
\end{array}\right.\ \textrm{ and }\ b_i(v)=\left\{\begin{array}{cl}
b_i(w) & \textrm{ for } i<s-1\\
b_s(w) & \textrm{for } i=s-1.\\
\end{array}\right.$$
Furthermore, in the case 2, the quiver $Q_u$ has no hole
meaning that the variety $\overline{PuQ}/Q$ is smooth and $u'$ has to
be equal to $u$. In this case we only need to determine $v'$.

Let us now consider the quiver $Q$ obtained by intersecting in $Q_w$
the quivers $Q_v$ and $Q_{w'}$. The quiver of $v'$ has to be a
subquiver of this quiver such that (see remark \ref{u'v'}) all the
holes of $Q_{w'}$ are holes of $Q_{v'}$ and $Q_{v'}$ may have one more hole
corresponding to the hole of $v$ which is not a hole of $w$. 

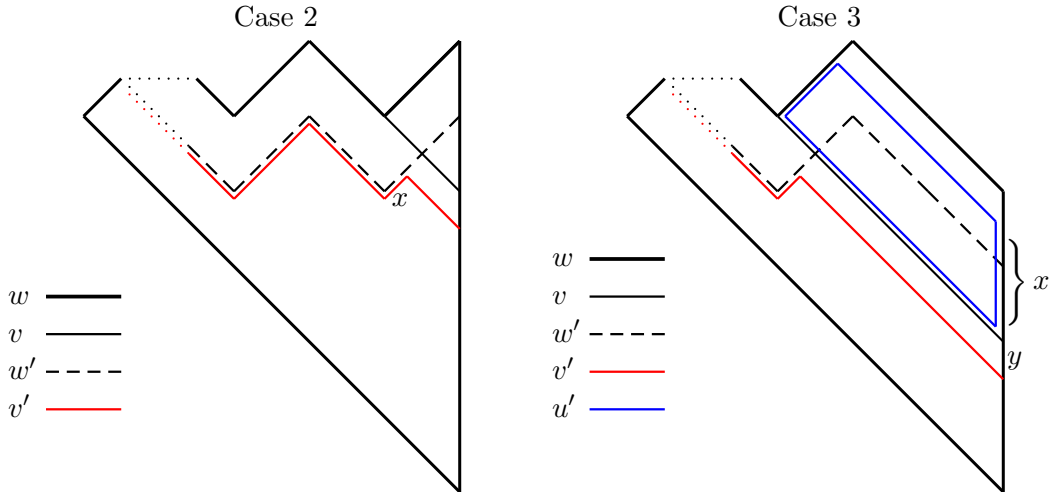
\psset{xunit=1cm}
\psset{yunit=1cm}
\centerline{\begin{pspicture*}(-7,-5)(0.1,2)
\put(-3,1.2){Case 2}
\psline[linewidth=0.04](-3,0)(-3.5,0.5)
\psline[linewidth=0.04](-4.5,0.5)(-5,0)
\psline[linewidth=0.04](-3,0)(-2,1)
\psline[linewidth=0.04](-1,0)(-2,1)
\psline[linewidth=0.03,linestyle=dotted](-4.4,0.5)(-3.6,0.5)
\psline[linewidth=0.04](-1,0)(0,1)
\psline[linewidth=0.04](0,-5)(0,1)
\psline[linewidth=0.04](-5,0)(0,-5)
\psline[linewidth=0.04](-1,0)(0,1)
\psline(-1,0)(0,-1)
\psline[linewidth=0.03,linestyle=dotted](-4.4,0.4)(-3.6,-0.4)
\psline[linewidth=0.03,linestyle=dashed](-3,-1)(-3.6,-0.4)
\psline[linewidth=0.03,linestyle=dashed](-3,-1)(-2,0)
\psline[linewidth=0.03,linestyle=dashed](-1,-1)(-2,0)
\psline[linewidth=0.03,linestyle=dashed](-1,-1)(0,0)
\put(-6,-2.5){$w$}
\put(-6,-3){$v$}
\put(-6,-3.5){$w'$}
\put(-6,-4){$v'$}
\psline[linewidth=0.05](-5.5,-2.4)(-4.5,-2.4)
\psline[linewidth=0.03,linestyle=dashed](-5.5,-3.4)(-4.5,-3.4)
\psline[linewidth=0.03](-5.5,-2.9)(-4.5,-2.9)
\psline[linewidth=0.03,linecolor=red](-5.5,-3.9)(-4.5,-3.9)
\psline[linewidth=0.03,linestyle=dotted,linecolor=red](-4.4,0.3)(-3.6,-0.5)
\psline[linewidth=0.03,linecolor=red](-3,-1.1)(-3.6,-0.5)
\psline[linewidth=0.03,linecolor=red](-3,-1.1)(-2,-0.1)
\psline[linewidth=0.03,linecolor=red](-1,-1.1)(-2,-0.1)
\psline[linewidth=0.03,linecolor=red](-1,-1.1)(-0.7,-0.8)
\psline[linewidth=0.03,linecolor=red](0,-1.5)(-0.7,-0.8)
\put(-0.9,-1.2){$x$}
\end{pspicture*}
\psset{xunit=1cm}
\psset{yunit=1cm}
\begin{pspicture*}(-7,-5)(1,2)
\put(-3,1.2){Case 3}
\psline[linewidth=0.04](-3,0)(-3.5,0.5)
\psline[linewidth=0.04](-4.5,0.5)(-5,0)
\psline[linewidth=0.04](-3,0)(-2,1)
\psline[linewidth=0.04](-1,0)(-2,1)
\psline[linewidth=0.03,linestyle=dotted](-4.4,0.5)(-3.6,0.5)
\psline[linewidth=0.04](-1,0)(0,-1)
\psline[linewidth=0.04](-5,0)(0,-5)
\psline[linewidth=0.04](0,-1)(0,-5)
\psline(-3,0)(0,-3)
\psline[linewidth=0.03,linestyle=dotted](-4.4,0.4)(-3.6,-0.4)
\psline[linewidth=0.03,linestyle=dashed](-3,-1)(-3.6,-0.4)
\psline[linewidth=0.03,linestyle=dashed](-3,-1)(-2,0)
\psline[linewidth=0.03,linestyle=dashed](-1,-1)(-2,0)
\psline[linewidth=0.03,linestyle=dashed](-1,-1)(0,-2)
\put(-6,-2){$w$}
\put(-6,-2.5){$v$}
\put(-6,-3){$w'$}
\put(-6,-3.5){$v'$}
\put(-6,-4){$u'$}
\psline[linewidth=0.05](-5.5,-1.9)(-4.5,-1.9)
\psline[linewidth=0.03,linestyle=dashed](-5.5,-2.9)(-4.5,-2.9)
\psline[linewidth=0.03](-5.5,-2.4)(-4.5,-2.4)
\psline[linewidth=0.03,linecolor=red](-5.5,-3.4)(-4.5,-3.4)
\psline[linewidth=0.03,linecolor=blue](-5.5,-3.9)(-4.5,-3.9)
\psline[linewidth=0.03,linestyle=dotted,linecolor=red](-4.4,0.3)(-3.6,-0.5)
\psline[linewidth=0.03,linecolor=red](-3,-1.1)(-3.6,-0.5)
\psline[linewidth=0.03,linecolor=red](-3,-1.1)(-2.7,-0.8)
\psline[linewidth=0.03,linecolor=red](0,-3.5)(-2.7,-0.8)
\psline[linewidth=0.03,linecolor=blue](-0.1,-1.4)(-2.2,0.7)
\psline[linewidth=0.03,linecolor=blue](-0.10,-1.4)(-0.1,-2.8)
\psline[linewidth=0.03,linecolor=blue](-2.9,-0)(-0.1,-2.8)
\psline[linewidth=0.03,linecolor=blue](-2.9,0)(-2.2,0.7)
\put(0.05,-3.3){$y$}
\put(-0.5,-2.3){$\left.
    \begin{array}{c}
\ \\
\\
    \end{array}
\right\}x$}
\end{pspicture*}}

In case 2, the quiver $Q_{v'}$ has $s+1$ holes and the sequences
$(a_i(v'))_{i\in[1,s+2]}$ and $(b_i(v'))_{i\in[0,s+1]}$ are given by
$$a_i(v')=\left\{\begin{array}{cl}
a_i(w') & \textrm{ for } i\leq s\\
b_s(w')-\frac{1}{2}-x & \textrm{for } i=s+1\\
\end{array}\right.\ \textrm{ and }\ b_i(v')=\left\{\begin{array}{cl}
b_i(w') & \textrm{ for } i<s\\
x & \textrm{for } i=s\\
\frac{1}{2} & \textrm{for } i=s+1\\
\end{array}\right.$$
where $x\in[0,c_{s}]$. In case 3, the quiver $Q_{v'}$ has $s$ holes
and the sequences $(a_i(v'))_{i\in[1,s+1]}$ and
$(b_i(v'))_{i\in[0,s]}$ are given by
$$a_i(v')=\left\{\begin{array}{cl}
a_i(w') & \textrm{ for } i\leq s-1\\
a_s(w)+b_{s-1}(w)-c_{s-1}+y & \textrm{for } i=s\\
\end{array}\right.\ \textrm{ and }\ b_i(v')=\left\{\begin{array}{cl}
b_i(w') & \textrm{ for } i<s-1\\
c_{s-1}-y & \textrm{for } i=s-1\\
b_{s}(w') & \textrm{for } i=s.\\
\end{array}\right.$$
We also get that $u'$ has a unique hole,
$$\left\{\begin{array}{l}
a_0(u')=x \textrm{ and } a_1(u')=a_s(w)+b_{s-1}(w)-x=a_1(u)+b_{0}(u)-x\\
b_0(u')=x \textrm{ and } b_1(u')=b_1(u)=b_s(w)
\end{array}\right.$$
In case 3, we have $y\in[c_{s-1}-c_s,c_{s-1}]$,
$x\in[b_{s-1}(w')-c_{k},b_{s-1}(w')]$ and $b_{s-1}-x=c_s-y$. Indeed,
the last formula is given by the fact that last hole of $Q_{u'v'}$ has
to be the same hole as the last hole of $Q_{w'}$. This implies that
$b_{s-1}(v')+b_0(u')=b_{s-1}(w')$ and the equality. The fact that
$y\in[c_{s-1}-c_s,c_{s-1}]$ comes from the fact that $X(v')$ is a
Schubert subvariety of $X(v)$ and that $X(u'v')=X(w')$.

The Schubert subvariety $X(\theta)$ is contained in $X(v)$, contains
$X(v')$ and is stable by the stabiliser of $X(v)$. It must have the
same hole as $v'$ except for those not corresponding to holes of
$v$. In our case we have to fill the $s^{\rm th}$ hole (in case 2) or
the $(s-1)^{\rm th}$ hole (in case 3) of $v'$ to obtain $\theta$: 

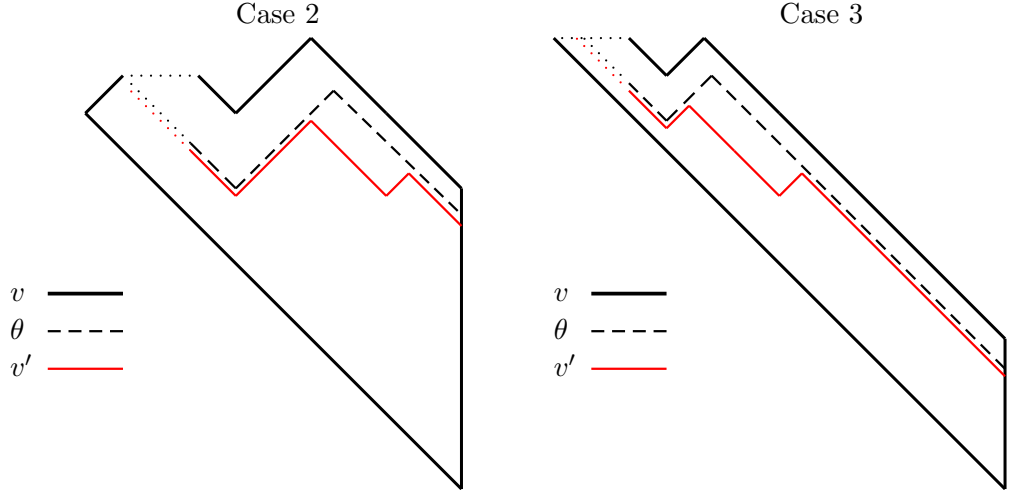
\psset{xunit=1cm}
\psset{yunit=1cm}
\centerline{\begin{pspicture*}(-7,-5)(0.1,2)
\put(-3,1.2){Case 2}
\psline[linewidth=0.04](-3,0)(-3.5,0.5)
\psline[linewidth=0.04](-4.5,0.5)(-5,0)
\psline[linewidth=0.04](-3,0)(-2,1)
\psline[linewidth=0.04](-1,0)(-2,1)
\psline[linewidth=0.03,linestyle=dotted](-4.4,0.5)(-3.6,0.5)
\psline[linewidth=0.04](0,-5)(0,-1)
\psline[linewidth=0.04](-5,0)(0,-5)
\psline[linewidth=0.04](-1,0)(0,-1)
\psline[linewidth=0.03,linestyle=dotted](-4.4,0.4)(-3.6,-0.4)
\psline[linewidth=0.03,linestyle=dashed](-3,-1)(-3.6,-0.4)
\psline[linewidth=0.03,linestyle=dashed](-3,-1)(-1.7,0.3)
\psline[linewidth=0.03,linestyle=dashed](-0.05,-1.3)(-1.7,0.3)
\put(-6,-2.5){$v$}
\put(-6,-3){$\theta$}
\put(-6,-3.5){$v'$}
\psline[linewidth=0.05](-5.5,-2.4)(-4.5,-2.4)
\psline[linewidth=0.03,linestyle=dashed](-5.5,-2.9)(-4.5,-2.9)
\psline[linewidth=0.03,linecolor=red](-5.5,-3.4)(-4.5,-3.4)
\psline[linewidth=0.03,linestyle=dotted,linecolor=red](-4.4,0.3)(-3.6,-0.5)
\psline[linewidth=0.03,linecolor=red](-3,-1.1)(-3.6,-0.5)
\psline[linewidth=0.03,linecolor=red](-3,-1.1)(-2,-0.1)
\psline[linewidth=0.03,linecolor=red](-1,-1.1)(-2,-0.1)
\psline[linewidth=0.03,linecolor=red](-1,-1.1)(-0.7,-0.8)
\psline[linewidth=0.03,linecolor=red](0,-1.5)(-0.7,-0.8)
\end{pspicture*}
\psset{xunit=1cm}
\psset{yunit=1cm}
\begin{pspicture*}(-7,-5)(0.1,2)
\put(-3,1.2){Case 3}
\psline[linewidth=0.04](-3,0)(-4,1)
\psline[linewidth=0.03,linestyle=dotted](-5.1,1)(-5.9,1)
\psline[linewidth=0.04](-4.5,0.5)(-5,1)
\psline[linewidth=0.04](-4,1)(-4.5,0.5)
\psline[linewidth=0.04](-6,1)(0,-5)
\psline[linewidth=0.04](0,-3)(0,-5)
\psline[linewidth=0.04](-3,0)(0,-3)
\psline[linewidth=0.03,linestyle=dotted](-5.1,0.5)(-5.6,1)
\psline[linewidth=0.03,linestyle=dashed](-3.9,0.5)(-0,-3.4)
\psline[linewidth=0.03,linestyle=dashed](-4.5,-0.1)(-3.9,0.5)
\psline[linewidth=0.03,linestyle=dashed](-4.5,-0.1)(-5,0.4)
\put(-6,-2.5){$v$}
\put(-6,-3){$\theta$}
\put(-6,-3.5){$v'$}
\psline[linewidth=0.03,linestyle=dashed](-5.5,-2.9)(-4.5,-2.9)
\psline[linewidth=0.05](-5.5,-2.4)(-4.5,-2.4)
\psline[linewidth=0.03,linecolor=red](-5.5,-3.4)(-4.5,-3.4)
\psline[linewidth=0.03,linestyle=dotted,linecolor=red](-5.1,0.4)(-5.7,1)
\psline[linewidth=0.03,linecolor=red](-4.2,0.1)(-3.6,-0.5)
\psline[linewidth=0.03,linecolor=red](-4.2,0.1)(-4.5,-0.2)
\psline[linewidth=0.03,linecolor=red](-5,0.3)(-4.5,-0.2)
\psline[linewidth=0.03,linecolor=red](-3,-1.1)(-3.6,-0.5)
\psline[linewidth=0.03,linecolor=red](-3,-1.1)(-2.7,-0.8)
\psline[linewidth=0.03,linecolor=red](0,-3.5)(-2.7,-0.8)
\end{pspicture*}}

In case 2, the integers $(a_i(\theta))_{i\in[1,s+1]}$ and
$(b_i(\theta))_{i\in[0,s]}$ associated to $\theta$ are given by 
$$a_i(\theta)=\left\{\begin{array}{cl}
a_i(v') & \textrm{ for } i\neq s\\
a_{s+1}(v')+a_s(v') & \textrm{for } i=s\\
\end{array}\right.\ \textrm{ and }\ b_i(\theta)=\left\{\begin{array}{cl}
b_i(v') & \textrm{ for } i\neq s-1\\
b_{s-1}(v')+b_{s}(v') & \textrm{for } i=s-1\\
\end{array}\right.$$
In case 3, the integers associated to
$\theta$ are given by $(a_i(\theta))_{i\in[1,s]}$ and
$(b_i(\theta))_{i\in[0,s-1]}$ are given by
$$a_i(\theta)=\left\{\begin{array}{cl}
a_i(v') & \textrm{ for } i\neq s-1\\
a_{s-1}(v')+a_s(v') & \textrm{for } i=s-1\\
\end{array}\right.\ \textrm{ and }\ b_i(\theta)=\left\{\begin{array}{cl}
b_i(v') & \textrm{ for } i\neq s-2\\
b_{s-2}(v')+b_{s-1}(v') & \textrm{for } i=s-2\\
\end{array}\right.$$ 
It is now an easy calculation (and straightforward on the quiver) that
in case 2, the depth $c'_i$ of $\theta$ in the holes of $v$ is $c_s-x$ for the
$s^{\rm th}$ hole and $c_i$ for the other holes. In case 3, the depth
$c'_i$ of $\theta$ in the holes of $v$ is $y=c_s+x-b_{s-1}$ for the
$(s-1)^{\rm th}$ hole and $c_i$ for the other holes.

We can now calculate in case 2:
$$\dim{X(u')}-\codim_{X(w')}(X(v'))=
\frac{(b_s(w)+\frac{1}{2})(b_s(w)-\frac{1}{2})}{2}-
\frac{(b_s(w')-x+\frac{1}{2})(b_s(w')-x-\frac{1}{2})}{2}$$
and because $b_s(w')=b_s(w)$ we get:
$$\dim{X(u')}-\codim_{X(w')}(X(v'))=xb_s(w)-\frac{x^2}{2}.$$
On the other hand we have, 
$$\Gamma(w',w)-\Gamma(\theta,v)=
\sum_{i=1}^sc_i(a_i(w)+b_i(w))-\sum_{i=1}^{s}c'_i(a_i(v)+b_i(v)),$$
a simple calculation gives 
$$\Gamma(w',w)-\Gamma(\theta,v)=x(a_s(w)+b_s(w)).$$
Now the fact that the $s^{\rm th}$ pic of $w$ was smaller than the
$(s-1)^{\rm th}$ pic means that $a_{s}(w)\geq b_{s}(w)$ so that we get
the inequality :
$$2(\dim{X(u')}-\codim_{X(w')}(X(v')))\leq
\Gamma(w',w)-\Gamma(\theta,v).$$
The theorem is proved in case 2. In case 3, we have:
$$
\begin{array}{cl}
\dim{X(u')}-\codim_{X(w')}(X(v'))&
=\frac{x(x+1)}{2}+x(a_s(w)+b_{s-1}(w)-x)-\frac{x(x+1)}{2}-xa_s(w')\\
&=x(a_s(w)+ b_{s-1}(w)-x-a_s(w)-c_s+c_{s-1})\\
\end{array}$$
and finaly:
$$\dim{X(u')}-\codim_{X(w')}(X(v'))=x(c_{s-1}-y).$$
On the other hand we have, 
$$\Gamma(w',w)-\Gamma(\theta,v)=
\sum_{i=1}^sc_i(a_i(w)+b_i(w))-\sum_{i=1}^{s-1}c'_i(a_i(v)+b_i(v))$$ 
a simple calculation gives 
$$\Gamma(w',w)-\Gamma(\theta,v)=(c_{s-1}-y)(a_{s-1}(w)+b_{s-1}(w))+
(b_{s-1}-x)(a_{s}(w)+b_{s}(w)).$$
Now the fact that the $(s-1)^{\rm th}$ pic of $w$ is smaller than the
$(s-2)^{\rm th}$ pic means that we have $a_{s-1}(w)\geq
b_{s-1}(w)$. Furthermore, we have $x\leq b_{s-1}(w')=b_{s-1}(w)$ and
$a_s(w)+b_s(w)\geq0$ so that we get the inequality :
$$2(\dim{X(u')}-\codim_{X(w')}(X(v')))\leq
\Gamma(w',w)-\Gamma(\theta,v).$$
The theorem is proved in case 3.

\subsection{Exceptional cases}
\label{casexcep}

We are left to deal with three cases: quadrics and minuscule varieties
for $E_6$ and $E_7$.

\subsubsection{Quadrics}

For quadrics,
let us remark that all Schubert varieties except one are locally
factorial (the quivers have only one pic) so that in all cases except
one we have $\Xh(\wh)=X(w)$ and there is nothing to prove. The only
non locally factorial Schubert variety (we are in $\k^{2p}$ with a
non degenerate quadratic form) is given by:
$$X(w)=\{x\in\p(\k^{2p})\ /\ x \textrm{ is isotropic and } x\in
  W^{\perp}_{p-2}\}$$ 
for a fixed isotropic subspace $W_{p-2}$ of dimension $p-2$. The
  associated quiver has $p$ vertices and is given by

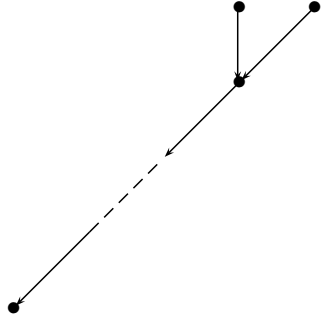
\psset{xunit=1cm}
\psset{yunit=1cm}
\centerline{\begin{pspicture*}(-3,-2.2)(5,2.2)
\psline[linestyle=dashed,linewidth=0.02](-0.1,-0.1)(-1,-1)
\psline[linewidth=0.02]{->}(0.97,2)(0.97,1.01)
\psline[linewidth=0.02]{->}(2,2)(1.02,1.02)
\psline[linewidth=0.02]{->}(-1,-1)(-1.98,-1.98)
\psline[linewidth=0.02]{<-}(0,0)(1,1)
\put(-2.1,-2.1){$\bullet$}
\put(1.9,1.9){$\bullet$}
\put(0.9,1.9){$\bullet$}
\put(0.9,0.9){$\bullet$}
\end{pspicture*}}

In particular, the resolution $\pih:\Xh(\wh)\to X(w)$ is given by
$p:\overline{PuQ}\times^QX(v)\to X(w)$ where $\overline{PuQ}/Q$ is of
dimension 1. The fiber of the morphism $\pih$ is at most 1. On
the other hand, as it is $\stab(X(w))$-equivariant, the fiber is
strictly positive on Schubert subvarieties stable under
$\stab(X(w))$. These subvarieties are of codimension at least 3 and
the result follows.

More generally, if the morphism $\pih$ is of the form
$p:\overline{PuQ}\times^QX(v)\to X(w)$ and its fiber is of dimension
at most one then the morphism $\pih$ has to be $IH$-small.

\subsubsection{The $E_6$ case}

We have seen that if the Schubert variety is locally factorial
(i.e. its quiver has a unique pic) or if the morphism $\pih$ is of the
form $p:\overline{PuQ}\times^QX(v)\to X(w)$ and its fiber is of
dimension at most one then the morphism $\pih$ has to be
$IH$-small. We are now going to list the morphisms $\pih$ obtained
from construction \ref{lespetites} not satisfying these properties and
verify that $\pih$ is small.

For $E_6$, all the morphisms $\pih$ not verifying the preceding
properties are of the form $p:\overline{PuQ}\times^QX(v)\to X(w)$. We
list here the quiver of $X(w)$ indicating on each quiver the
quivers of $u$ and $v$.

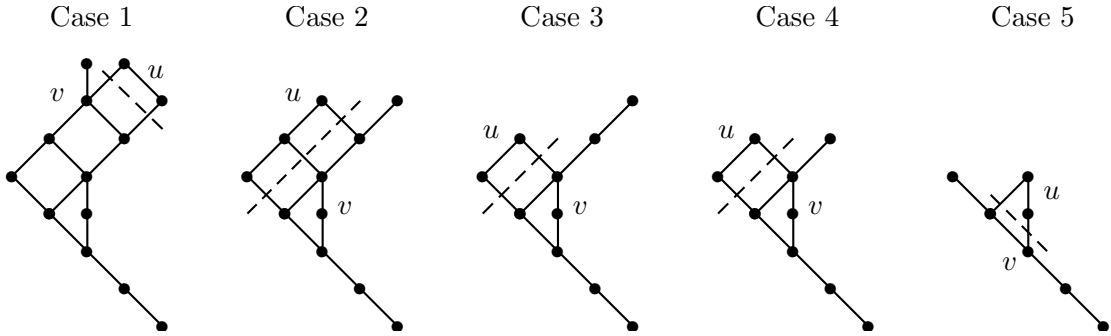
\vs 0.2 cm

\psset{xunit=0.5cm}
\psset{yunit=0.5cm}
\centerline{\begin{pspicture*}(-1,-4.1)(5,4.5)
\put(0.5,2){Case 1}
\psline(0,0)(1,1)
\psline(2,2)(1,1)
\psline(2,2)(3,3)
\psline(1,-1)(4,2)
\psline(2,2)(2,3)
\psline(3,3)(4,2)
\psline(2,2)(3,1)
\psline(1,1)(2,0)
\psline(0,0)(4,-4)
\psline(2,0)(2,-2)
\put(-0.1,-0.1){$\bullet$}
\put(0.4,0.4){$\bullet$}
\put(0.9,0.9){$\bullet$}
\put(1.4,1.4){$\bullet$}
\put(0.9,-0.1){$\bullet$}
\put(1.4,0.4){$\bullet$}
\put(1.9,0.9){$\bullet$}
\put(0.4,-0.6){$\bullet$}
\put(0.9,1.4){$\bullet$}
\put(0.9,-0.6){$\bullet$}
\put(1.9,-2.1){$\bullet$}
\put(1.4,-1.6){$\bullet$}
\put(0.9,-1.1){$\bullet$}
\put(0.4,-0.6){$\bullet$}
\psline[linestyle=dashed](2.4,2.8)(4,1.25)
\put(1.8,1.3){$u$}
\put(0.5,1){$v$}
\end{pspicture*}
\begin{pspicture*}(-1,-4.1)(5,4.5)
\put(0.5,2){Case 2}
\psline(0,0)(1,1)
\psline(2,2)(1,1)
\psline(1,-1)(4,2)
\psline(2,2)(3,1)
\psline(1,1)(2,0)
\psline(0,0)(4,-4)
\psline(2,0)(2,-2)
\put(-0.1,-0.1){$\bullet$}
\put(0.4,0.4){$\bullet$}
\put(0.9,0.9){$\bullet$}
\put(0.9,-0.1){$\bullet$}
\put(1.4,0.4){$\bullet$}
\put(1.9,0.9){$\bullet$}
\put(0.4,-0.6){$\bullet$}
\put(0.9,-0.6){$\bullet$}
\put(1.9,-2.1){$\bullet$}
\put(1.4,-1.6){$\bullet$}
\put(0.9,-1.1){$\bullet$}
\put(0.4,-0.6){$\bullet$}
\psline[linestyle=dashed](0,-1)(3,2)
\put(1.2,-0.5){$v$}
\put(0.5,1){$u$}
\end{pspicture*}
\begin{pspicture*}(-1,-4.1)(5,4.5)
\put(0.5,2){Case 3}
\psline(0,0)(1,1)
\psline(1,-1)(4,2)
\psline(1,1)(2,0)
\psline(0,0)(4,-4)
\psline(2,0)(2,-2)
\put(-0.1,-0.1){$\bullet$}
\put(0.4,0.4){$\bullet$}
\put(0.9,-0.1){$\bullet$}
\put(1.4,0.4){$\bullet$}
\put(1.9,0.9){$\bullet$}
\put(0.4,-0.6){$\bullet$}
\put(0.9,-0.6){$\bullet$}
\put(1.9,-2.1){$\bullet$}
\put(1.4,-1.6){$\bullet$}
\put(0.9,-1.1){$\bullet$}
\put(0.4,-0.6){$\bullet$}
\psline[linestyle=dashed](0,-1)(2,1)
\put(1.2,-0.5){$v$}
\put(0,0.5){$u$}
\end{pspicture*}
\begin{pspicture*}(-1,-4.1)(5,4.5)
\put(0.5,2){Case 4}
\psline(0,0)(1,1)
\psline(1,-1)(3,1)
\psline(1,1)(2,0)
\psline(0,0)(4,-4)
\psline(2,0)(2,-2)
\put(-0.1,-0.1){$\bullet$}
\put(0.4,0.4){$\bullet$}
\put(0.9,-0.1){$\bullet$}
\put(1.4,0.4){$\bullet$}
\put(0.4,-0.6){$\bullet$}
\put(0.9,-0.6){$\bullet$}
\put(1.9,-2.1){$\bullet$}
\put(1.4,-1.6){$\bullet$}
\put(0.9,-1.1){$\bullet$}
\put(0.4,-0.6){$\bullet$}
\psline[linestyle=dashed](0,-1)(2,1)
\put(1.2,-0.5){$v$}
\put(0,0.5){$u$}
\end{pspicture*}
\begin{pspicture*}(-1,-4.1)(5,4.5)
\put(0.5,2){Case 5}
\psline(1,-1)(2,0)
\psline(0,0)(4,-4)
\psline(2,0)(2,-2)
\put(-0.1,-0.1){$\bullet$}
\put(0.9,-0.1){$\bullet$}
\put(0.4,-0.6){$\bullet$}
\put(0.9,-0.6){$\bullet$}
\put(1.9,-2.1){$\bullet$}
\put(1.4,-1.6){$\bullet$}
\put(0.9,-1.1){$\bullet$}
\put(0.4,-0.6){$\bullet$}
\psline[linestyle=dashed](1,-0.5)(2.5,-2)
\put(1.2,-0.3){$u$}
\put(0.65,-1.25){$v$}
\end{pspicture*}}

The dimension of the fiber in these morphisms is at most $f=2$ except
in the second case where it is at most $f=3$. The Schubert
subvarieties $X(w')$ stable under $\stab(X(w))$ of codimension not
bigger than $2f$ are the following:


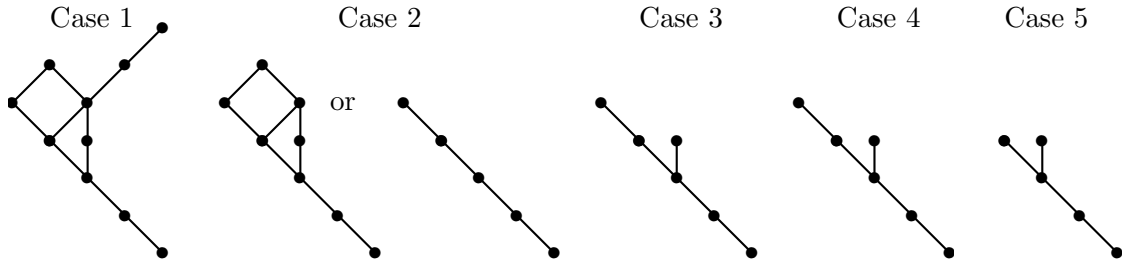
\psset{xunit=0.5cm}
\psset{yunit=0.5cm}
\centerline{\begin{pspicture*}(-0.1,-4.3)(4.7,3.1)
\put(0.5,1){Case 1}
\psline(0,0)(1,1)
\psline(1,-1)(4,2)
\psline(1,1)(2,0)
\psline(0,0)(4,-4)
\psline(2,0)(2,-2)
\put(-0.1,-0.1){$\bullet$}
\put(0.4,0.4){$\bullet$}
\put(0.9,-0.1){$\bullet$}
\put(1.4,0.4){$\bullet$}
\put(1.9,0.9){$\bullet$}
\put(0.4,-0.6){$\bullet$}
\put(0.9,-0.6){$\bullet$}
\put(1.9,-2.1){$\bullet$}
\put(1.4,-1.6){$\bullet$}
\put(0.9,-1.1){$\bullet$}
\put(0.4,-0.6){$\bullet$}
\end{pspicture*}
\begin{pspicture*}(-0.7,-4.3)(5.5,3.1)
\put(1.5,1){Case 2}
\psline(0,0)(1,1)
\psline(1,-1)(2,0)
\psline(1,1)(2,0)
\psline(0,0)(4,-4)
\psline(2,0)(2,-2)
\put(-0.1,-0.1){$\bullet$}
\put(0.4,0.4){$\bullet$}
\put(0.9,-0.1){$\bullet$}
\put(0.4,-0.6){$\bullet$}
\put(0.9,-0.6){$\bullet$}
\put(1.9,-2.1){$\bullet$}
\put(1.4,-1.6){$\bullet$}
\put(0.9,-1.1){$\bullet$}
\put(0.4,-0.6){$\bullet$}
\put(1.4,-0.1){or}
\end{pspicture*}
\hskip -2 cm
\begin{pspicture*}(-3,-4.3)(4.5,3.1)
\psline(0,0)(4,-4)
\put(-0.1,-0.1){$\bullet$}
\put(0.4,-0.6){$\bullet$}
\put(1.9,-2.1){$\bullet$}
\put(1.4,-1.6){$\bullet$}
\put(0.9,-1.1){$\bullet$}
\put(0.4,-0.6){$\bullet$}
\end{pspicture*}
\begin{pspicture*}(-0.5,-4.3)(4.5,3.1)
\put(0.5,1){Case 3}
\psline(0,0)(4,-4)
\psline(2,-1)(2,-2)
\put(-0.1,-0.1){$\bullet$}
\put(0.4,-0.6){$\bullet$}
\put(0.9,-0.6){$\bullet$}
\put(1.9,-2.1){$\bullet$}
\put(1.4,-1.6){$\bullet$}
\put(0.9,-1.1){$\bullet$}
\put(0.4,-0.6){$\bullet$}
\end{pspicture*}
\begin{pspicture*}(-0.5,-4.3)(4.1,3.1)
\put(0.5,1){Case 4}
\psline(0,0)(4,-4)
\psline(2,-1)(2,-2)
\put(-0.1,-0.1){$\bullet$}
\put(0.4,-0.6){$\bullet$}
\put(0.9,-0.6){$\bullet$}
\put(1.9,-2.1){$\bullet$}
\put(1.4,-1.6){$\bullet$}
\put(0.9,-1.1){$\bullet$}
\put(0.4,-0.6){$\bullet$}
\end{pspicture*}
\begin{pspicture*}(-0.1,-4.3)(4.1,3.1)
\put(0.5,1){Case 5}
\psline(1,-1)(4,-4)
\psline(2,-1)(2,-2)
\put(0.4,-0.6){$\bullet$}
\put(0.9,-0.6){$\bullet$}
\put(1.9,-2.1){$\bullet$}
\put(1.4,-1.6){$\bullet$}
\put(0.9,-1.1){$\bullet$}
\put(0.4,-0.6){$\bullet$}
\end{pspicture*}}
These quivers $Q$ are obtained from the quiver $Q_w$ by removing all the
vertices smaller 
than a hole $i$ of $Q_w$. It is
now easy to see that any subvariety $Z_K$ of the Bott-Samelson
resolution $\pit:\Xt(\wt)\to X(w)$ such that $\pit(Z_K)=X(w')$ is
contained in the divisor $Z_{i}$. We thus have $\pit^{-1}(X(w'))=Z_i$
and $\pih^{-1}(X(w'))$ is contained in the image of $Z_i$ in
$\Xh(\wh)$. Seeing $\Xh(\wh)$ as a configuration variety, the image of
$Z_i$ in $\Xh(\wh)$ is the configuration variety
$\overline{PuQ}\times^QX(v')$ where $Q_{v'}=Q\cap Q_v$. In particular,
the dimension of the fiber of $\pih$ above $X(w')$ is 1, 1, 2, 1, 1, 1
in the different cases and the morphism $\pih$ is always $IH$-small.

\subsubsection{The $E_7$ case}

We proceed in the same way in this case and list the quivers having at
least two pics and for which the fiber is at least two. 
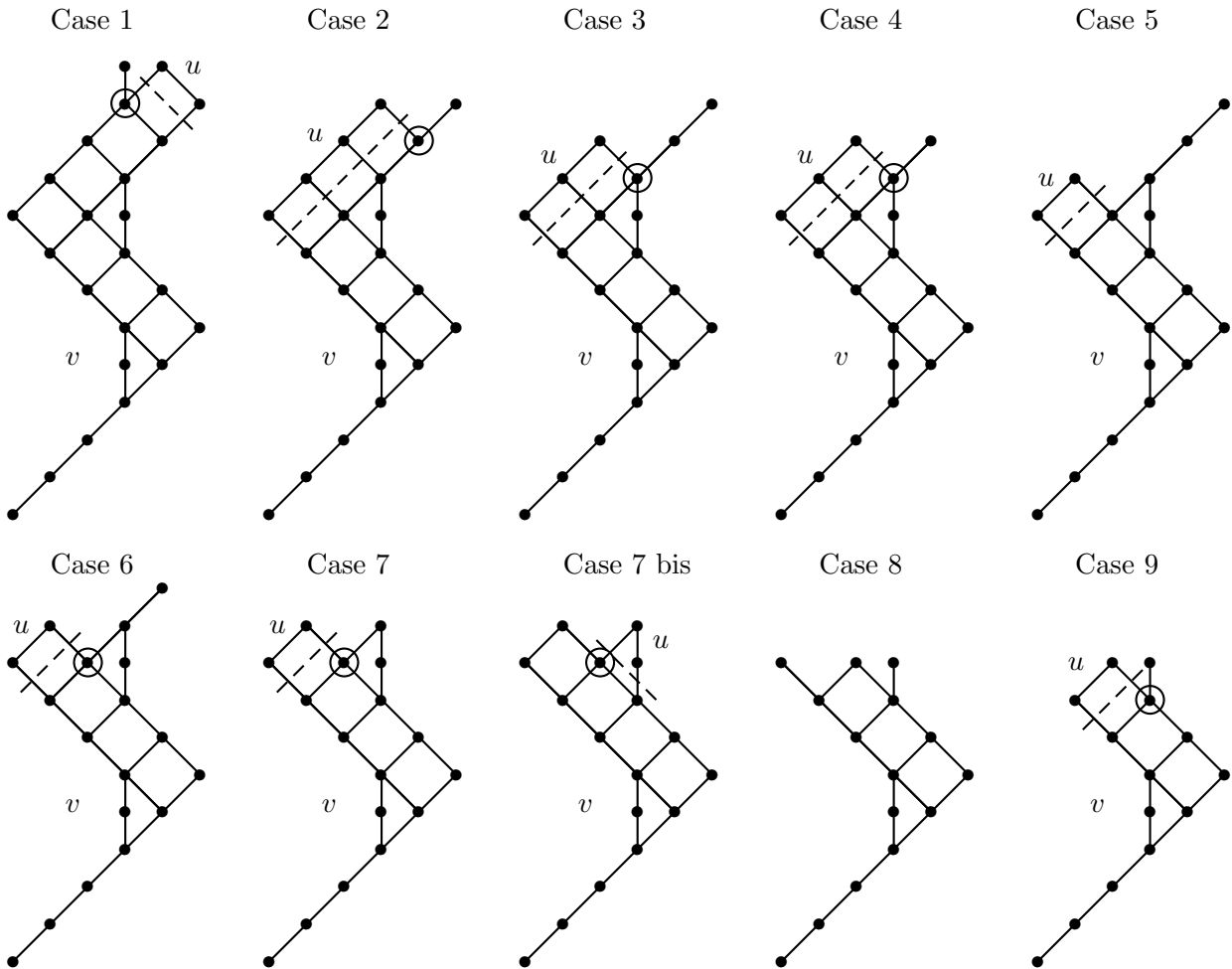
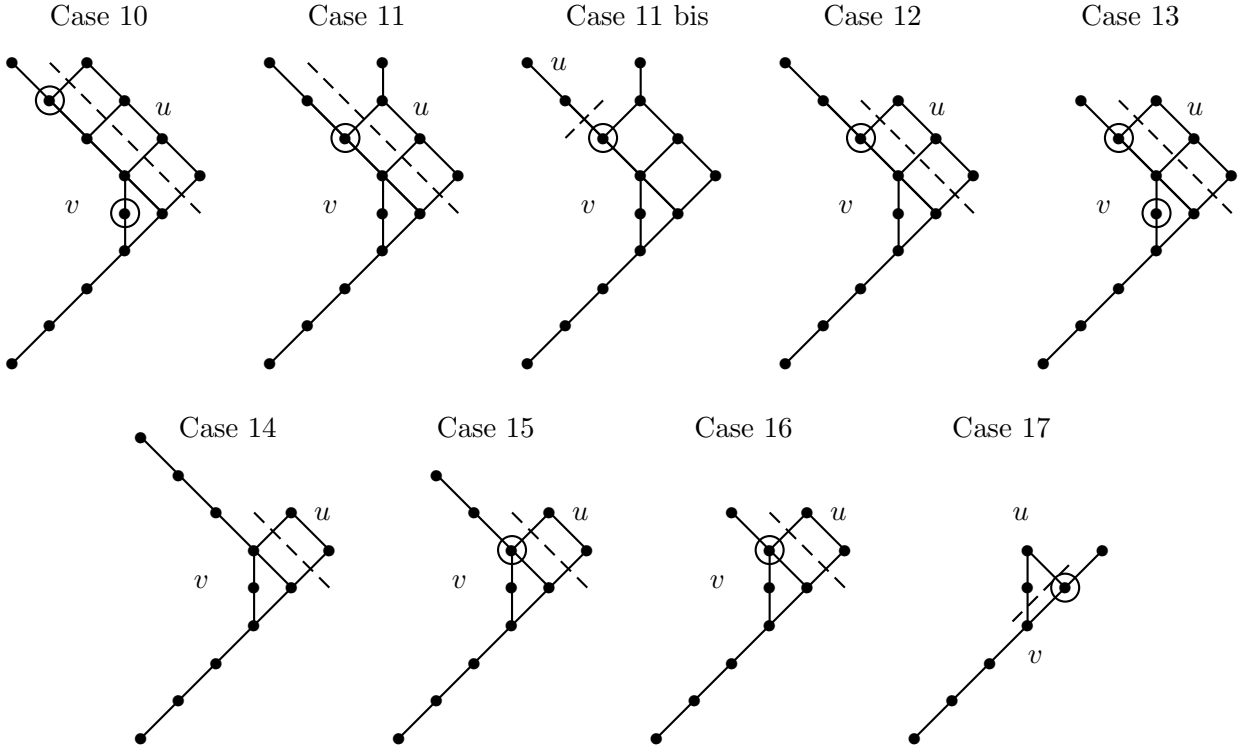
\vs 0.5 cm

\psset{xunit=0.5cm}
\psset{yunit=0.5cm}
\centerline{
}

All the
resolutions are of type $\overline{PuQ}\times^QX(v')$ except case
8. The maximal dimension $f$ of the fiber in all these cases is given
by
$$
\begin{array}{|c|c|c|c|c|c|c|c|c|c|c|c|c|c|c|c|c|c|c|}
\hline
{\rm case}&1&2&3&4&5&6&7&7\ {\rm bis}&9&10&11&11\ {\rm
  bis}&12&13&14&15&16&17\\
\hline
f&2&4&3&3&2&2&2&2&2&4&4&2&3&3&2&2&2&2\\
\hline
\end{array}$$

 We have
cercled the vertices $i$ such that the quivers $Q_{w'}$ obtained from
the quiver $Q_w$ by removing all the vertices smaller (for
$\preccurlyeq$) than the hole $i$ of $Q_w$ are the quivers of the
Schubert subvarieties $X(w')$ stable under $\stab(X(w))$ of
codimension not superior to $2f$. The codimension of $X(w')$ in $X(w)$
is given by:
$$
\begin{array}{|c|c|c|c|c|c|c|c|c|c|c|c|c|c|c|c|c|}
\hline
{\rm case}&1&2&3&4&6&7&7\ {\rm bis}&9&10&11&11\ {\rm
  bis}&12&13&15&16&17\\
\hline
{\rm codim}&3&3&4&3&4&3&3&3&3\ {\rm or}\ 8&5&5&4&3\ {\rm or}\ 6&4&3&3\\
\hline
\end{array}$$

Remark that in case 5 and 14 there is no such Schubert subvariety so
that the morphism is already small. In all the other cases and as for
the $E_6$ case, it is easy to see that any subvariety $Z_K$ of
the Bott-Samelson resolution $\pit:\Xt(\wt)\to X(w)$ such that
$\pit(Z_K)=X(w')$ is contained in the divisor $Z_{i}$. We thus have
$\pit^{-1}(X(w'))=Z_i$ and $\pih^{-1}(X(w'))$ is contained in the
image of $Z_i$ in $\Xh(\wh)$. Seeing $\Xh(\wh)$ as a configuration
variety, the image of $Z_i$ in $\Xh(\wh)$ is the configuration variety
$\overline{PuQ}\times^QX(v')$ where $Q_{v'}=Q_{w'}\cap Q_v$. In
particular, the dimension of the fiber of $\pih$ above $X(w')$ is
given by
$$
\begin{array}{|c|c|c|c|c|c|c|c|c|c|c|c|c|c|c|c|c|}
\hline
{\rm case}&1&2&3&4&6&7&7\ {\rm bis}&9&10&11&11\ {\rm
  bis}&12&13&15&16&17\\
\hline
{\rm dim}&1&1&1&1&1&1&1&1&1\ {\rm or}\ 3&2&2&1&1\ {\rm or}\ 2&1&1&1\\
\hline
\end{array}$$
and the morphism $\pih$ is always $IH$-small in these cases. 
We are left with case 8 for which the resolution is of the form
$\overline{PtQ}\times^Q\overline{RuS}\times^SX(v)$ and the partitions
of the quivers are given by:

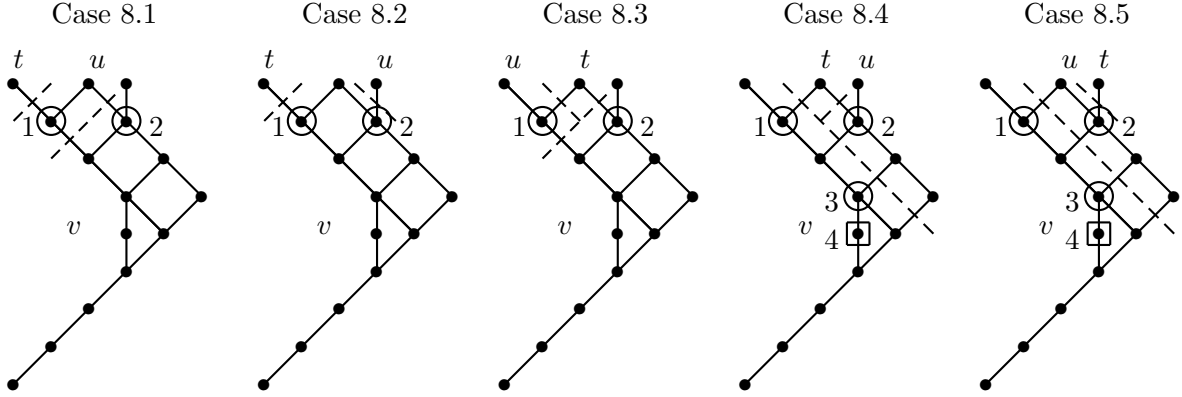
\centerline{\begin{pspicture*}(-0.2,-8.3)(6.2,3)
\put(0.5,0.8){Case 8.1}
\pscircle(1,-1){0.2}
\pscircle(3,-1){0.2}
\put(0.1,-0.7){1}
\put(1.8,-0.7){2}
\psline(1,-1)(2,0)
\psline(2,0)(5,-3)
\psline(4,-4)(0,0)
\psline(3,-3)(4,-2)
\psline(5,-3)(0,-8)
\psline(2,-2)(3,-1)
\psline(0,0)(4,-4)
\psline(3,-1)(3,0)
\psline(3,-3)(3,-5)
\put(-0.1,-0.1){$\bullet$}
\put(0.9,-0.1){$\bullet$}
\put(1.4,-0.1){$\bullet$}
\put(-0.1,-4.1){$\bullet$}
\put(0.4,-3.6){$\bullet$}
\put(0.9,-3.1){$\bullet$}
\put(1.4,-2.6){$\bullet$}
\put(1.4,-2.1){$\bullet$}
\put(1.4,-0.6){$\bullet$}
\put(1.9,-2.1){$\bullet$}
\put(1.4,-1.6){$\bullet$}
\put(0.9,-1.1){$\bullet$}
\put(2.4,-1.6){$\bullet$}
\put(1.9,-1.1){$\bullet$}
\put(0.4,-0.6){$\bullet$}
\psline[linestyle=dashed](1,0)(0,-1)
\psline[linestyle=dashed](2.7,-0.3)(1,-2)
\put(0.7,-2){$v$}
\put(0,0.2){$t$}
\put(1,0.2){$u$}
\end{pspicture*}
\begin{pspicture*}(-0.2,-8.3)(6.2,3)
\put(0.5,0.8){Case 8.2}
\pscircle(1,-1){0.2}
\pscircle(3,-1){0.2}
\put(0.1,-0.7){1}
\put(1.8,-0.7){2}
\psline(1,-1)(2,0)
\psline(2,0)(5,-3)
\psline(4,-4)(0,0)
\psline(3,-3)(4,-2)
\psline(5,-3)(0,-8)
\psline(2,-2)(3,-1)
\psline(0,0)(4,-4)
\psline(3,-1)(3,0)
\psline(3,-3)(3,-5)
\put(-0.1,-0.1){$\bullet$}
\put(0.9,-0.1){$\bullet$}
\put(1.4,-0.1){$\bullet$}
\put(-0.1,-4.1){$\bullet$}
\put(0.4,-3.6){$\bullet$}
\put(0.9,-3.1){$\bullet$}
\put(1.4,-2.6){$\bullet$}
\put(1.4,-2.1){$\bullet$}
\put(1.4,-0.6){$\bullet$}
\put(1.9,-2.1){$\bullet$}
\put(1.4,-1.6){$\bullet$}
\put(0.9,-1.1){$\bullet$}
\put(2.4,-1.6){$\bullet$}
\put(1.9,-1.1){$\bullet$}
\put(0.4,-0.6){$\bullet$}
\psline[linestyle=dashed](1,0)(0,-1)
\psline[linestyle=dashed](2.4,-0)(3.5,-1)
\put(0.7,-2){$v$}
\put(0,0.2){$t$}
\put(1.5,0.2){$u$}
\end{pspicture*}\begin{pspicture*}(-0.2,-8.3)(6.2,3)
\put(0.5,0.8){Case 8.3}
\pscircle(1,-1){0.2}
\pscircle(3,-1){0.2}
\put(0.1,-0.7){1}
\put(1.8,-0.7){2}
\psline(1,-1)(2,0)
\psline(2,0)(5,-3)
\psline(4,-4)(0,0)
\psline(3,-3)(4,-2)
\psline(5,-3)(0,-8)
\psline(2,-2)(3,-1)
\psline(0,0)(4,-4)
\psline(3,-1)(3,0)
\psline(3,-3)(3,-5)
\put(-0.1,-0.1){$\bullet$}
\put(0.9,-0.1){$\bullet$}
\put(1.4,-0.1){$\bullet$}
\put(-0.1,-4.1){$\bullet$}
\put(0.4,-3.6){$\bullet$}
\put(0.9,-3.1){$\bullet$}
\put(1.4,-2.6){$\bullet$}
\put(1.4,-2.1){$\bullet$}
\put(1.4,-0.6){$\bullet$}
\put(1.9,-2.1){$\bullet$}
\put(1.4,-1.6){$\bullet$}
\put(0.9,-1.1){$\bullet$}
\put(2.4,-1.6){$\bullet$}
\put(1.9,-1.1){$\bullet$}
\put(0.4,-0.6){$\bullet$}
\psline[linestyle=dashed](1,0)(2,-1)
\psline[linestyle=dashed](2.7,-0.3)(1,-2)
\put(0.7,-2){$v$}
\put(0,0.2){$u$}
\put(1,0.2){$t$}
\end{pspicture*}\begin{pspicture*}(-0.2,-8.3)(6.2,3)
\put(0.5,0.8){Case 8.4}
\psline(2.7,-4.3)(3.3,-4.3)
\psline(2.7,-3.7)(2.7,-4.3)
\psline(2.7,-3.7)(3.3,-3.7)
\psline(3.3,-4.3)(3.3,-3.7)
\put(1.05,-2.2){4}
\pscircle(1,-1){0.2}
\pscircle(3,-1){0.2}
\pscircle(3,-3){0.2}
\put(0.1,-0.7){1}
\put(1.8,-0.7){2}
\put(1.05,-1.7){3}
\psline(1,-1)(2,0)
\psline(2,0)(5,-3)
\psline(4,-4)(0,0)
\psline(3,-3)(4,-2)
\psline(5,-3)(0,-8)
\psline(2,-2)(3,-1)
\psline(0,0)(4,-4)
\psline(3,-1)(3,0)
\psline(3,-3)(3,-5)
\put(-0.1,-0.1){$\bullet$}
\put(0.9,-0.1){$\bullet$}
\put(1.4,-0.1){$\bullet$}
\put(-0.1,-4.1){$\bullet$}
\put(0.4,-3.6){$\bullet$}
\put(0.9,-3.1){$\bullet$}
\put(1.4,-2.6){$\bullet$}
\put(1.4,-2.1){$\bullet$}
\put(1.4,-0.6){$\bullet$}
\put(1.9,-2.1){$\bullet$}
\put(1.4,-1.6){$\bullet$}
\put(0.9,-1.1){$\bullet$}
\put(2.4,-1.6){$\bullet$}
\put(1.9,-1.1){$\bullet$}
\put(0.4,-0.6){$\bullet$}
\psline[linestyle=dashed](1,0)(5,-4)
\psline[linestyle=dashed](2.7,-0.3)(2,-1)
\put(0.7,-2){$v$}
\put(1.5,0.2){$u$}
\put(1,0.2){$t$}
\end{pspicture*}\begin{pspicture*}(-0.2,-8.3)(6.2,3)
\put(0.5,0.8){Case 8.5}
\psline(2.7,-4.3)(3.3,-4.3)
\psline(2.7,-3.7)(2.7,-4.3)
\psline(2.7,-3.7)(3.3,-3.7)
\psline(3.3,-4.3)(3.3,-3.7)
\put(1.05,-2.2){4}
\pscircle(1,-1){0.2}
\pscircle(3,-1){0.2}
\pscircle(3,-3){0.2}
\put(0.1,-0.7){1}
\put(1.8,-0.7){2}
\put(1.05,-1.7){3}
\psline(1,-1)(2,0)
\psline(2,0)(5,-3)
\psline(4,-4)(0,0)
\psline(3,-3)(4,-2)
\psline(5,-3)(0,-8)
\psline(2,-2)(3,-1)
\psline(0,0)(4,-4)
\psline(3,-1)(3,0)
\psline(3,-3)(3,-5)
\put(-0.1,-0.1){$\bullet$}
\put(0.9,-0.1){$\bullet$}
\put(1.4,-0.1){$\bullet$}
\put(-0.1,-4.1){$\bullet$}
\put(0.4,-3.6){$\bullet$}
\put(0.9,-3.1){$\bullet$}
\put(1.4,-2.6){$\bullet$}
\put(1.4,-2.1){$\bullet$}
\put(1.4,-0.6){$\bullet$}
\put(1.9,-2.1){$\bullet$}
\put(1.4,-1.6){$\bullet$}
\put(0.9,-1.1){$\bullet$}
\put(2.4,-1.6){$\bullet$}
\put(1.9,-1.1){$\bullet$}
\put(0.4,-0.6){$\bullet$}
\psline[linestyle=dashed](1,0)(5,-4)
\psline[linestyle=dashed](2.4,-0)(3.5,-1)
\put(0.7,-2){$v$}
\put(1.5,0.2){$t$}
\put(1,0.2){$u$}
\end{pspicture*}}

The maximal dimension $f$ of the fiber in all these cases is given
by
$$
\begin{array}{|c|c|c|c|c|c|}
\hline
{\rm case}&8.1&8.2&8.3&8.4&8.5\\
\hline
f&3&2&3&5&5\\
\hline
\end{array}$$

We have cercled and numeroted the vertices $i$ such that the quivers
$Q_{w'}$ obtained from the quiver $Q_w$ by removing all the vertices
smaller (for $\preccurlyeq$) than a fixed subset of the holes of $Q_w$
are the quivers of the Schubert subvarieties $X(w')$ stable under
$\stab(X(w))$ of codimension not superior to $2f$. Let $A$ be a non
empty subset of $\{1,2, 3\}$ and let $Q_{w'}$ be the quiver obtained
by removing the vertices smaller than the vertices in $A$. The
codimension of $X(w')$ in $X(w)$ is given by (here $A$ is $\{1\}$,
$\{2\}$, $\{1,2\}$ or in the last two cases $\{3\}$):
$$
\begin{array}{|c|c|c|c|c|c|}
\hline
{\rm case}&8.1&8.2&8.3&8.4&8.5\\
\hline
{\rm codim}&3,\ 3\ {\rm or}\ 5&3,\ 3\ {\rm or}\ 5&3,\ 3\ {\rm or}\
5&3,\ 3,\  5\ {\rm or}\ 8&3,\ 3,\  5\ {\rm or}\ 8\\
\hline
\end{array}$$

 Suppose that $A$ is a subset of $\{1,2\}$, it is easy to see that any
 subvariety $Z_K$ of the Bott-Samelson
resolution $\pit:\Xt(\wt)\to X(w)$ such that $\pit(Z_K)=X(w')$ is
contained in the variety $Z_A$. The fiber $\pih^{-1}(X(w'))$ is thus
contained in the image of $Z_A$ in $\Xh(\wh)$. Seeing $\Xh(\wh)$ as a
configuration variety, this image in $\Xh(\wh)$ is the configuration
variety $\overline{PtQ}\times^Q\overline{Ru'S}\times^SX(v')$ where $Q_{u'}$
and $Q_{v'}$ are obtained respectively from $Q_u$ and $Q_v$ by
removing the vertices smaller than one vertex in $A\cap Q_u$
respectively $A\cap Q_v$. In particular, the dimension of the fiber of
$\pih$ above $X(w')$ is given by
$$\begin{array}{|c|c|c|c|c|c|}
\hline
{\rm case}&8.1&8.2&8.3&8.4&8.5\\
\hline
{\rm dim}&1,\ 1\ {\rm or}\ 1&1,\ 1\ {\rm or}\ 2&1,\ 1\ {\rm
  or}\ 1&1,\ 1\ {\rm or}\ 1&1,\ 1\ {\rm or}\ 1\\
\hline
\end{array}$$

and the morphism $\pih$ is always $IH$-small in these cases. We are
left with the case where $A=\{3\}$. In this case it is not hard to see
that any subvariety $Z_K$ of the Bott-Samelson
resolution $\pit:\Xt(\wt)\to X(w)$ such that $\pit(Z_K)=X(w')$ is
contained in the variety $Z_{\{2,3\}}$ or in $Z_4$. The fiber
$\pih^{-1}(X(w'))$ is thus contained in the image of $Z_{\{2,3\}}$ or
of $Z_4$ in $\Xh(\wh)$. Seeing $\Xh(\wh)$ as a configuration variety,
the image of $Z_{\{2,3\}}$ in $\Xh(\wh)$ is the configuration variety
$\overline{PtQ}\times^Q\overline{Ru'S}\times^SX(v')$ where $Q_{u'}$
and $Q_{v'}$ are obtained respectively from $Q_u$ and $Q_v$ by
removing the vertices smaller than one vertex in $\{2,3\}\cap Q_u$
respectively $\{2,3\}\cap Q_v$. The image of $Z_4$ is the
configuration variety
$\overline{PtQ}\times^Q\overline{RuS}\times^SX(v')$ where $Q_{v'}$ is
obtained from $Q_v$ by removing the vertices smaller than the vertex 4. In
particular, the fiber of $\pih$ above $X(w')$ has two components whose
dimensions are given by
$$\begin{array}{|c|c|c|c|c|}
\hline
{\rm case}&8.4\ ; Z_{\{2,3\}}&8.4\ ; Z_4 &8.5\ ; Z_{\{2,3\}}&8.5\ ; Z_4\\
\hline
{\rm dim}&2&3&2&3\\
\hline
\end{array}$$
and the morphism $\pih$ is always $IH$-small.

\subsection{Small resolution}

Let us now describe all $IH$-small resolutions of minuscule Schubert
varieties whenever they exist. Having adopted (and described) the
relative minimal models point of view, we use the following result of
B. Totaro \cite{Totaro} using a key result of J. Wisniewski
\cite{wisniewski}:

\begin{theo}
\label{totaro}
  Any $IH$-small resolution of $X$ is a small relative
  minimal model for $X$.
\end{theo}

Looking for $IH$-small resolution we only have to check in our list of
minimal models. Furthermore, because of theorem \ref{petit}, the
morphism $\pih:\Xh(\wh)\to X(w)$ from any minimal model to $X(w)$ is
$IH$-small so that we get the following

\begin{coro}
The $IH$-small resolutions of $X(w)$ are given by the morphisms
$\pih:\Xh(\wh)\to X(w)$ obtained from construction \ref{lespetites}
with $\Xh(\wh)$ smooth.  
\end{coro}

We now give a combinatorial description of these varieties. Let $Q_v$
be a quiver associated to a minuscule Schubert variety $X(v)$ and $i$
a vertex of $Q_v$. 

\begin{defi}
\label{pic-min}
The vertex $i$ of $Q_v$ is called minuscule if $\b(i)$ is a minuscule
simple root of the sub-Dynkin diagram of $G$ defined by $\supp(v)$.
\end{defi}

Construction \ref{lespetites} gives a partition of the quiver $Q_w$
into subquivers $Q_{w_i}$ which are quivers of minuscule Schubert
varieties having only one pic. We have

\begin{theo}
\label{lisse}
  The variety $\Xh(\wh)$ obtained from construction \ref{tous}
  is smooth if and only if for all $i$, the unique pic $p_i$ of $Q_{w_i}$
  is minuscule in $Q_{w_i}$.
\end{theo}

\begin{preu}
  We have seen that the variety $\Xh(\wh)$ is a sequence of locally
  trivial fibrations with fiber Schubert varieties $X(w_i)$. The
  theorem will follow from:
\end{preu}

  \begin{prop}
    A minuscule Schubert variety $X(w)$ is smooth if and only if $Q_w$
    has a unique pic $p$ and $p$ is minuscule in $Q_w$.
  \end{prop}

  \begin{preu}
    We know from \cite{BrionPolo} that a minuscule Schubert variety
    $X(w)$ is smooth if and only if it is homogeneous under its
    stabiliser. It is easy to verify that the quiver of any minuscule
    homogeneous variety has a unique pic which is minuscule. 

Conversely, according to proposition \ref{stab}, the variety is
homogeneous under its stabiliser if and only if the quiver $Q_w$ has
no non virtual hole. Now we have seen that for $A_n$ the quiver of any
Schubert variety is of the form (we have cercled the non virtual holes
of the quiver):

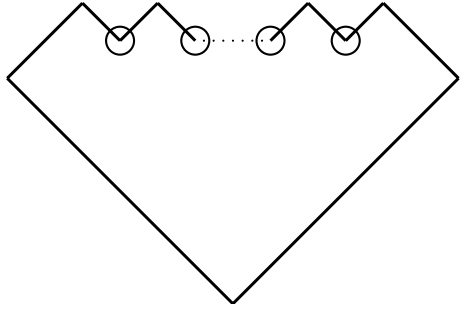
\psset{xunit=0.5cm}
\psset{yunit=0.5cm}
\centerline{\begin{pspicture*}(-9,-7)(9,2)
\psline[linewidth=0.04](-6,-1)(-4,1)
\psline[linewidth=0.04](-3,0)(-4,1)
\psline[linewidth=0.04](-3,0)(-2,1)
\psline[linewidth=0.04](-1,0)(-2,1)
\psline[linewidth=0.03,linestyle=dotted](-0.8,0)(0.8,0)
\psline[linewidth=0.04](6,-1)(4,1)
\psline[linewidth=0.04](3,0)(4,1)
\psline[linewidth=0.04](3,0)(2,1)
\psline[linewidth=0.04](1,0)(2,1)
\psline[linewidth=0.04](-6,-1)(0,-7)
\psline[linewidth=0.04](0,-7)(6,-1)
\pscircle(-3,0){0.2}
\pscircle(3,0){0.2}
\pscircle(-1,0){0.2}
\pscircle(1,0){0.2}
\end{pspicture*}}

and the only case where there is a unique pic is when there is no
hole. In this case we have the quiver of a grassmannian and it is
smooth see appendix). For the case of maximal isotropic subspaces (say
associated to the simple root $\a_n$ with the notations of
\cite{bourb}), the quiver is of the form (we have cercled the non
virtual holes of the quiver):

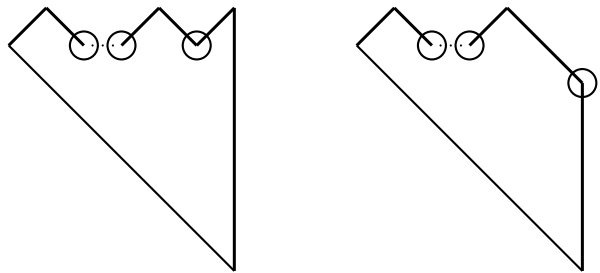
\psset{xunit=0.5cm}
\psset{yunit=0.5cm}
\centerline{\begin{pspicture*}(-5,-7)(4,1.5)
\pscircle(-1,0){0.2}
\pscircle(2,0){0.2}
\pscircle(0,0){0.2}
\psline(-3,0)(3,-6)
\psline[linewidth=0.04](-3,0)(-2,1)
\psline[linewidth=0.04](-1,0)(-2,1)
\psline[linewidth=0.03,linestyle=dotted](-0.8,0)(-0.2,0)
\psline[linewidth=0.04](3,-6)(3,1)
\psline[linewidth=0.04](2,0)(3,1)
\psline[linewidth=0.04](2,0)(1,1)
\psline[linewidth=0.04](0,0)(1,1)
\end{pspicture*}
\begin{pspicture*}(-5,-7)(4,1.5)
\pscircle(-1,0){0.2}
\pscircle(3,-1){0.2}
\pscircle(0,0){0.2}
\psline(-3,0)(3,-6)
\psline[linewidth=0.04](-3,0)(-2,1)
\psline[linewidth=0.04](-1,0)(-2,1)
\psline[linewidth=0.03,linestyle=dotted](-0.8,0)(-0.2,0)
\psline[linewidth=0.04](3,-6)(3,-1)
\psline[linewidth=0.04](3,-1)(1,1)
\psline[linewidth=0.04](0,0)(1,1)
\end{pspicture*}}

and there are three cases when there is only one pic namely 

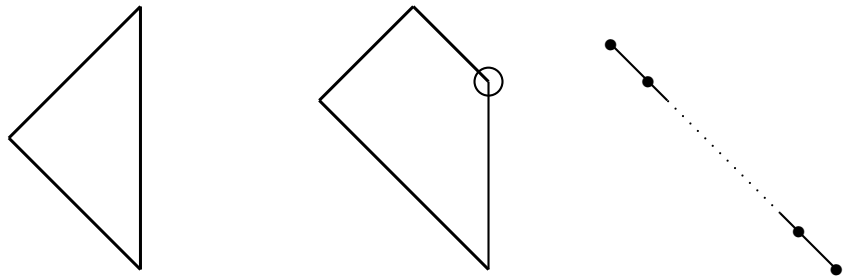
\psset{xunit=0.5cm}
\psset{yunit=0.5cm}
\centerline{\begin{pspicture*}(-5,-7)(4,1.5)
\psline[linewidth=0.04](3,1)(-0.5,-2.5)
\psline[linewidth=0.04](3,-6)(-0.5,-2.5)
\psline[linewidth=0.04](3,-6)(3,1)
\end{pspicture*}
\begin{pspicture*}(-5,-7)(4,1.5)
\pscircle(3,-1){0.2}
\psline(3,-1)(3,-6)
\psline[linewidth=0.04](3,-6)(-1.5,-1.5)
\psline[linewidth=0.04](-1.5,-1.5)(1,1)
\psline[linewidth=0.04](3,-1)(1,1)
\end{pspicture*}
\begin{pspicture*}(-5,-7)(4,1.5)
\psline(-3,0)(-1.5,-1.5)
\psline(1.5,-4.5)(3,-6)
\psline[linewidth=0.03,linestyle=dotted](-1.5,-1.5)(1.5,-4.5)
\put(-1.6,-0.1){$\bullet$}
\put(-1.1,-0.6){$\bullet$}
\put(0.9,-2.6){$\bullet$}
\put(1.4,-3.1){$\bullet$}
\end{pspicture*}}

In the second case, one of the two vertices $i_{n-1}$ and $i_n$ such
that $i_k$ is the smallest element (for $\preccurlyeq$) with
$\b(i_k)=\a_k$ with the notations of \cite{bourb} is a hole of the
quiver. In the first case the quiver is the quiver of the isotropic
grassmannian and in the third one it is the quiver of a projective
space. For the quadric case, the quiver has one of the four following
forms:

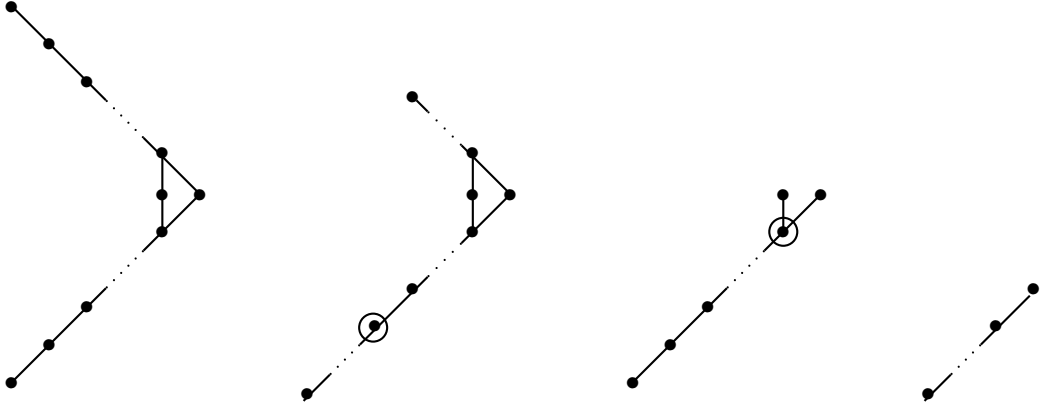
\psset{xunit=0.5cm}
\psset{yunit=0.5cm}
\centerline{\begin{pspicture*}(-1,-12)(7,0)
\psline(2,-2)(3.5,-3.5)
\psline[linewidth=0.03,linestyle=dotted](3.5,-3.5)(4.5,-4.5)
\psline(6,-6)(4.5,-4.5)
\psline(5,-5)(5,-6)
\psline(5,-7)(5,-6)
\psline(2,-10)(3.5,-8.5)
\psline[linewidth=0.03,linestyle=dotted](3.5,-8.5)(4.5,-7.5)
\psline(6,-6)(4.5,-7.5)
\put(0.9,-1.1){$\bullet$}
\put(1.4,-1.6){$\bullet$}
\put(2.9,-3.1){$\bullet$}
\put(2.4,-3.1){$\bullet$}
\put(2.4,-3.6){$\bullet$}
\put(2.4,-2.55){$\bullet$}
\put(0.9,-5.1){$\bullet$}
\put(1.4,-4.6){$\bullet$}
\psline(2,-10)(1,-11)
\put(0.4,-5.6){$\bullet$}
\psline(2,-2)(1,-1)
\put(0.4,-.6){$\bullet$}
\end{pspicture*}
\begin{pspicture*}(-1,-12)(7,0)
\psline(3.8,-3.8)(3.5,-3.5)
\psline[linewidth=0.03,linestyle=dotted](3.8,-3.8)(4.7,-4.7)
\psline(6,-6)(4.7,-4.7)
\psline(5,-5)(5,-6)
\psline(5,-7)(5,-6)
\psline(2,-10)(3.8,-8.2)
\psline[linewidth=0.03,linestyle=dotted](2,-10)(1.2,-10.8)
\psline[linewidth=0.03,linestyle=dotted](3.8,-8.2)(4.7,-7.3)
\psline(6,-6)(4.7,-7.3)
\psline(1.2,-10.8)(0.5,-11.5)
\put(1.6,-1.8){$\bullet$}
\put(1.6,-4.35){$\bullet$}
\put(2.9,-3.1){$\bullet$}
\put(2.4,-3.1){$\bullet$}
\put(2.4,-3.6){$\bullet$}
\put(2.4,-2.55){$\bullet$}
\put(0.2,-5.75){$\bullet$}
\pscircle(2.35,-9.55){0.2}
\put(1.1,-4.85){$\bullet$}
\end{pspicture*}
\begin{pspicture*}(-1,-12)(7,0)
\psline(5,-7)(5,-6)
\psline(2,-10)(3.5,-8.5)
\psline[linewidth=0.03,linestyle=dotted](3.5,-8.5)(4.5,-7.5)
\psline(6,-6)(4.5,-7.5)
\put(2.9,-3.1){$\bullet$}
\put(2.4,-3.1){$\bullet$}
\put(2.4,-3.6){$\bullet$}
\pscircle(5,-7){0.2}
\put(0.9,-5.1){$\bullet$}
\put(1.4,-4.6){$\bullet$}
\psline(2,-10)(1,-11)
\put(0.4,-5.6){$\bullet$}
\end{pspicture*}
\begin{pspicture*}(-1,-12)(7,0)
\psline(2,-10)(3.3,-8.7)
\psline[linewidth=0.03,linestyle=dotted](2,-10)(1.2,-10.8)
\psline(1.2,-10.8)(0.5,-11.5)
\put(1.6,-4.35){$\bullet$}
\put(0.2,-5.75){$\bullet$}
\put(1.1,-4.85){$\bullet$}
\end{pspicture*}}
and in the first and last cases we get respectively the quiver of a
quadric or the quiver of a projective space. In the two intermediate
cases, there is one hole in the quiver.

Finaly, it is an easy verification on the quivers of $E_6$ and $E_7$
to check that the proposition is true (cf. appendix).
\end{preu}

\subsection{Stringy polynomials}

Another way of proving the non existence of $IH$-small resolutions is
the following: because of theorem \ref{totaro} of 
any $IH$-small resolution $\Xt$ of a variety $X$ if it exists will
factor through the relative canonical model $\Xh$ (which will always
exists if $\Xt$ does). Furthermore, the resolution
$\Xt\to\Xh$ will be $IH$-small and in particular crepant. We can thus
use the stringy polynomial $E(\Xh,u,v)$ defined by V. Batyrev in
\cite{batyrev}. If $\Xh$ admits a crepant resolution then this
polynomial (which in general is a formal power serie) is a true
polynomial. To prove the non existence of $IH$-resolution, it would be
enough to prove that $E(\Xh,u,v)$ is not a polynomial.

Let us give an example where we make the full calculation. Let us
first recall the following definitions (for more details and more
general definitions, see \cite{batyrev}).

Let $X$ be a normal irreducible variety, we define the following
notations:
$$E(X,u,v)=\sum_{u,v}e^{p,q}(X)u^pv^q\ \ \ {\rm with}\ \ \
e^{p,q}(X)=\sum_{i}(-1)^ih^{p,q}(H^i_c(X,\mathbb{C}))$$
where $H^i_c(X,\mathbb{C})$ is the $i^{\rm th}$ cohomology group with
compact support and $h^{p,q}((H^i_c(X,\mathbb{C}))$ is the dimension
of its $(p,q)$-type component. The polynomial $E(X,u,v)$ is what
V. Batyrev call the Euler polynomial (or $E$-polynomial).

Assume now that $X$ is a gorenstein normal irreducible variety with at
  worst terminal singularities. Let $\pi:Y\to X$ a resolution of
  singularities such that the exceptional locus is a divisor $D$ whose
  irreducible components $(D_i)_{i\in I}$ are smooth divisors with
  only normal crossing. We then have
$$K_Y=\pi^* K_X+\sum_{i\in I}a_iD_i\ \ \ {\rm with}\ \ \ a_i>0.$$ 

For any subset $J\subset I$ we define
$$D_J=\left\{
  \begin{array}{cl}
\displaystyle{\bigcap_{j\in J}D_j}& \textrm{ if } J\neq\emptyset\\
Y& \textrm{ if } J=\emptyset\\
  \end{array}\right.\ \ \ \ \ \textrm{and}\ \ \ \ \
D_J^o=D_J\setminus\bigcup_{i\in I\setminus J}(D_J\cap D_i).$$

\begin{defi}
  The stringy function associated with the resolution $\pi:Y\to X$ is
  the following:
$$E_{\rm st}(X,u,v)=\sum_{J\subset I}E(D^o_J,u,v)\prod_{j\in
  J}\frac{uv-1}{(uv)^{a_j+1}-1}.$$ 
\end{defi}

Then V. Batyrev proves the following

\begin{theo}
  (\i) The function $E_{\rm st}(X,u,v)$ is independent of the
  resolution $\pi:Y\to X$ with exceptional locus of pure codimension 1
  given by smooth irreducible divisors with normal crossing.

(\i\i) If $X$ admits a crepant resolution $\pi:Y\to X$ (that is to say
$\pi^*K_X=K_Y$) then $E_{\rm st}(X,u,v)=E(Y,u,v)$ and it is a
polynomial.

(\i\i\i) In particular, if $X$ admits a crepant resolution, the stingy
Euler number 
$$e_{\rm st}(X)=\lim_{u,v\to0}E_{\rm st}(X,u,v)=\sum_{J\subset
  I}e(D^o_J)\prod_{j\in J}\frac{1}{1+a_j}$$
is an integer.
\end{theo}

We now give an example of a minuscule Schubert variety which is
singular non locally factorial and does not admit an $IH$-small
resolution.

\begin{exem}
  \label{example}
Let $G$ be $SO(12)$ and $w$ given by the following reduced writing
$\wt$ (the symetry $s_i$ is the simple reflection associated to the
$i^{\rm th}$ simple root with the notation of \cite{bourb}):
$$w=s_2s_4s_1s_3s_6s_2s_4s_3s_4s_4s_6.$$
The associated Schubert variety is the following ($\G_{iso}(k,12)$ is
the isotropic grassmannian, we denote by $\G_{iso}^1(6,12)$ and
$\G_{iso}^2(6,12)$ the homogenous varieties associated to the simple
roots $\a_5$ and $\a_6$):
$$X(w)=\left\{V\in \G_{iso}^2(6,12)\ /\ \dim(V\cap W_3)\geq1\ {\rm
    and}\ \dim(V\cap W_6)\geq3 \right\}$$
where $W_3\in\G_{iso}(3,12)$ and $W_6\in \G_{iso}^1(6,12)$. The quiver
    $Q_w$ is

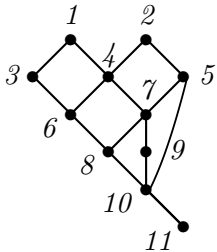
\psset{xunit=1cm}
\psset{yunit=1cm}
\centerline{\begin{pspicture*}(-1,-2.5)(5,1.2)
\psline(0,0)(2,-2)
\psline(0.5,0.5)(0,0)
\psline(0.5,0.5)(1.5,-0.5)
\psline(0.5,-0.5)(1.5,0.5)
\psline(1.5,-0.5)(1.5,-1.5)
\psline(1.5,-1.5)(2,-2)
\psline(1.5,0.5)(2,0)
\psline(2,0)(1,-1)
\parametricplot[linewidth=0.7pt,plotstyle=curve]{1.5}{2}{t 0.05 add t
  2 mul 2 sub 2 exp 2 div 2 sub}
\put(-0.1,-0.1){$\bullet$}
\put(0.4,0.4){$\bullet$}
\put(0.4,-0.6){$\bullet$}
\put(0.9,-0.1){$\bullet$}
\put(1.4,0.4){$\bullet$}
\put(1.4,-0.6){$\bullet$}
\put(0.9,-1.1){$\bullet$}
\put(1.4,-1.6){$\bullet$}
\put(1.9,-2.1){$\bullet$}
\put(1.9,-0.1){$\bullet$}
\put(1.4,-1.1){$\bullet$}
\put(-0.4,-0.1){3}
\put(0.4,0.7){1}
\put(0.1,-0.8){6}
\put(0.9,0.2){4}
\put(1.4,0.7){2}
\put(1.4,-0.3){7}
\put(0.6,-1.3){8}
\put(0.9,-1.8){10}
\put(1.45,-2.3){11}
\put(2.2,-0.1){5}
\put(1.8,-1.1){9}
\end{pspicture*}}

The variety $X(w)$ has for resolution the Bott-Samelson
resolution $\Xt(\wt)$. Moreover, because the morphism
$\pi:\Xt(\wt)\to X(w)$ is $B$-equivariant, the exceptional locus has
to be $B$-invariant and thus an union of $Z_K$. The only non
contracted divisors $Z_i$ of $\Xt(\wt)$ in $X(w)$ are $Z_1$ and
$Z_2$. Furthermore, the variety $Z_{\{1,2\}}$ is not contracted so
that the exceptional locus $D$ is the union
$$D=\cup_{i=3}^{11}Z_i.$$
All $Z_i$ are smooth and intersect transversally. Denote by $D_1$ and
$D_2$ the images of $Z_1$ and $Z_2$ in $X(w)$. The ample generator of
the Picard group of $X(w)$ is given by $\L=D_1+D_2$. We have
$$\pi^*\L=\sum_{i=1}^{11}Z_i.$$
Formulae of paragraph \ref{canonique} and lemma
\ref{som-haut} give us:
$$-K_{\Xt(\wt)}=\sum_{i=1}^{11}(h(i)+1)Z_i$$
and proposition \ref{cano-ample} gives us:
$$-K_{X(w)}=7D_1+7D_2=7\L.$$
In particular, we have
$$K_{\Xt(w)}-\pi^*K_{X(w)}=\sum_{i=1}^{11}(6-h(i))Z_i=
(Z_3+Z_4+Z_5)+2(Z_6+Z_7)+3(Z_8+Z_9)+4Z_{10}+5Z_{11}.$$ 
Remark that for $J\subset[3,9]$, the variety $Z_J^o$ is a sequence of
9 locally trivial fibrations in $\mathbb{A}^1$ of in points (there are
exactely $\vert J\vert$ points) over $\pu\times\pu$. In particular, we
have $e(Z^o_J)=4$ for all $J\subset[3,9]$.

Now we have the easy formula
$$\sum_{J\subset I}\prod_{j\in J}x_j=\prod_{i\in I}(1+x_i).$$
We can thus calculate in our situation:
$$e_{\rm st}(X(w))=4(1+\frac{1}{2})^3(1+\frac{1}{3})^2
(1+\frac{1}{4})^2(1+\frac{1}{5})(1+\frac{1}{6})=\frac{105}{2}.$$ 
We conclude that $X(w)$ has no $IH$-small resolution as given by
theorem \ref{lisse}.

This kind of calculation can be generalised, this with be done in a
subsequent paper. For example, the same calculation in the general
case where $X(w)$ is gorenstein (or equivalentely all the pics $p\in
p(Q_w)$ have the same height $h(w)$) gives the
following result: let us define for $i\in Q_w$ its coheight
$coh(i)=h(w)-h(i)$. Then we have
$$e_{\rm st}(X(w))=
\prod_{i\in
  Q_w}\left(1+\frac{1}{1+coh(i)}\right).$$
\end{exem}

\section{Appendix}

In this appendix we give the quivers of minuscule homogeneous
varieties and describe the quivers of minuscule Schubert varieties.

\subsection{Quivers of minuscule homogeneous varieties}

The following quiver is the quiver of the grassmannian of
$p$-dimensional subvector spaces of an $n$-dimensional vector
space. The morphism $\b$ associating to any vertex a simple root is
simply the vertical projection on the Dynkin diagram.

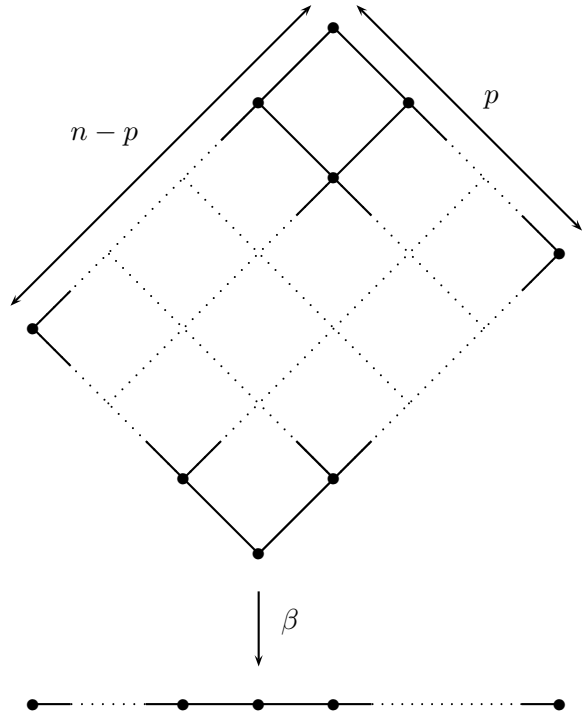
\psset{xunit=1cm}
\psset{yunit=1cm}
\centerline{\begin{pspicture*}(-5,-6)(5,5)
\psline[linewidth=0.03](0,4)(-1.5,2.5)
\psline(0,4)(1.5,2.5)
\psline(0.5,1.5)(-1,3)
\psline(-0.5,1.5)(1,3)
\psline{<->}(-0.3,4.3)(-4.3,0.3)
\psline{<->}(0.3,4.3)(3.3,1.3)
\put(2,3){$p$}
\put(-3.5,2.5){$n-p$}
\psline[linestyle=dotted](-1,3)(-4,0)
\psline[linestyle=dotted](0,2)(-3,-1)
\psline[linestyle=dotted](2,2)(-2,-2)
\psline[linestyle=dotted](3,1)(0,-2)
\psline[linestyle=dotted](-2,-2)(-4,0)
\psline[linestyle=dotted](-3,1)(0,-2)
\psline[linestyle=dotted](-2,2)(1,-1)
\psline[linestyle=dotted](0,2)(2,0)
\psline[linestyle=dotted](1,3)(3,1)
\psline(-2.5,-1.5)(-1,-3)
\psline(-2,-2)(-1.5,-1.5)
\psline(3,1)(2.5,0.5)
\psline(3,1)(2.5,1.5)
\psline(-4,0)(-3.5,0.5)
\psline(-4,0)(-3.5,-0.5)
\psline[linewidth=0.03](-1,-3)(0.5,-1.5)
\psline[linewidth=0.03](-0.5,-1.5)(0,-2)
\put(-0.1,3.9){$\bullet$}
\put(-0.1,1.9){$\bullet$}
\put(-1.1,2.9){$\bullet$}
\put(0.9,2.9){$\bullet$}
\put(-4.1,-0.1){$\bullet$}
\put(-0.1,-2.1){$\bullet$}
\put(-1.1,-3.1){$\bullet$}
\put(2.9,0.9){$\bullet$}
\put(-2.1,-2.1){$\bullet$}
\psline{->}(-1,-3.5)(-1,-4.5)
\psline(-4,-5)(-3.5,-5)
\psline[linestyle=dotted](-2.5,-5)(-3.5,-5)
\psline(-2.5,-5)(0.5,-5)
\psline[linestyle=dotted](0.5,-5)(2.5,-5)
\psline(2.5,-5)(3,-5)
\put(-4.1,-5.1){$\bullet$}
\put(-2.1,-5.1){$\bullet$}
\put(-1.1,-5.1){$\bullet$}
\put(-0.1,-5.1){$\bullet$}
\put(2.9,-5.1){$\bullet$}
\put(-0.7,-4){$\b$}
\end{pspicture*}}

It is easy to verify that this diagramm satifies the geometric
conditions of proposition \ref{geom-carq} so that it correponds to a
Schubert variety of dimension $p(n-p)$ of the grassmannian. It must be
the quiver of the grassmannian. 

In the same way, the quiver of the grassmannian of maximal
isotropic subspaces in a $2n$-dimensional vector space endowed with a
non degenerate quadratic form is given by (one more time, the morphism
$\b$ is given by the vertical projection on the Dynkin diangram) one
of the followinf form depending on the parity of $n$:

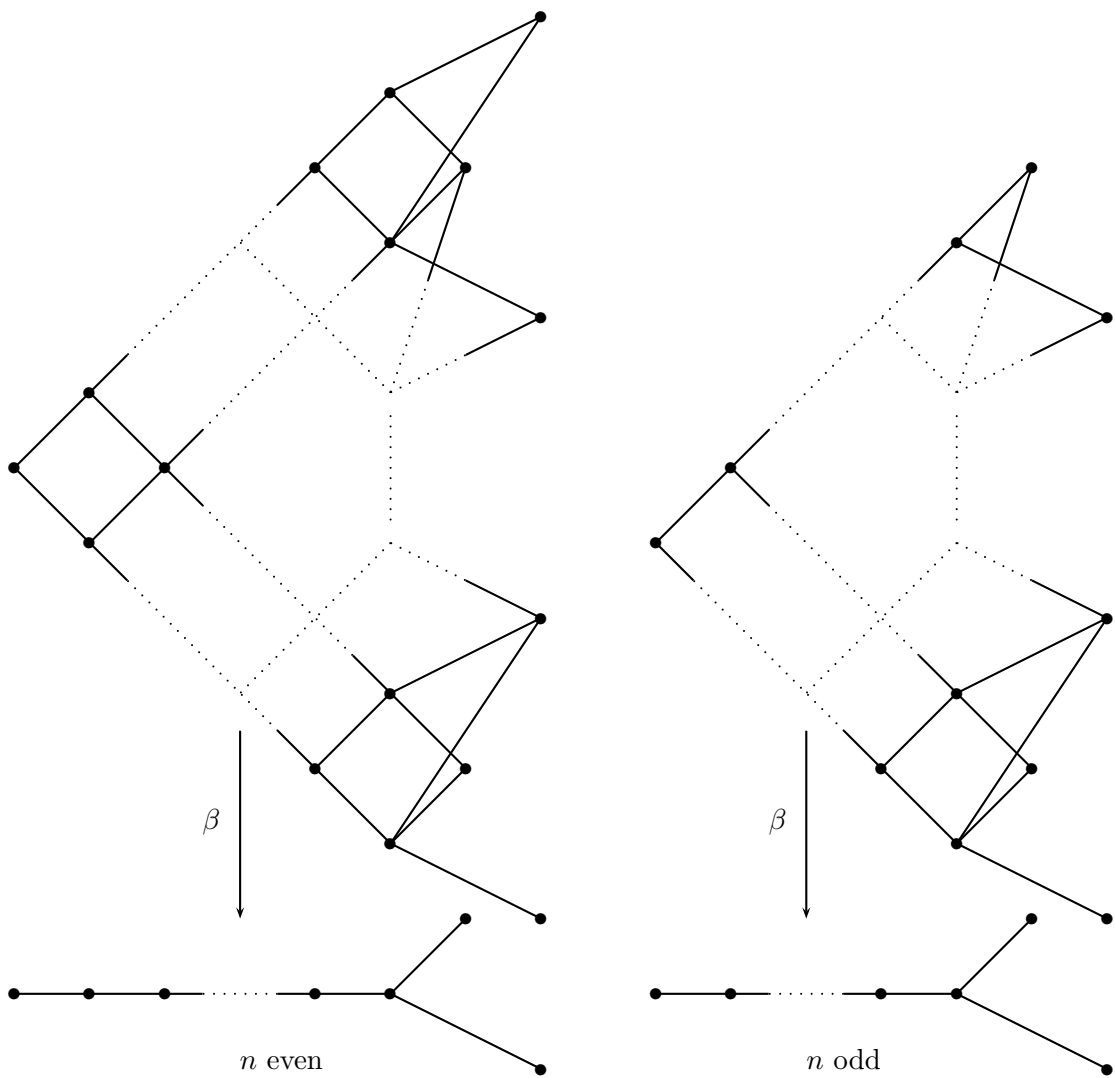
\psset{xunit=1cm}
\psset{yunit=1cm}
\centerline{
}

In the text we use for the quiver of the grassmannian and for the
quiver of the grassmannian of maximal isotropic subspaces the
following schematic version of the quivers:

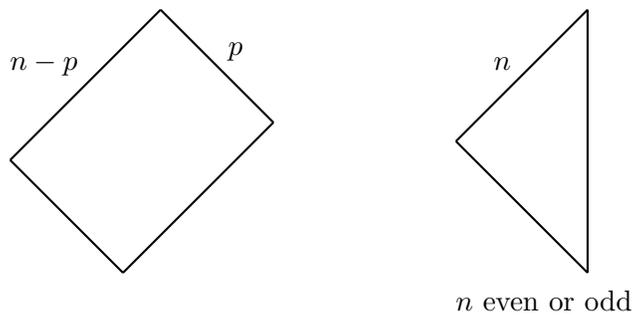
\psset{xunit=1cm}
\psset{yunit=1cm}
\centerline{%
}

For the even dimensional quadrics, we get the following quiver (we
often draw it like the diagram on the right even if $\b$ is not
exactely given by the projection and if we forget some arrows):

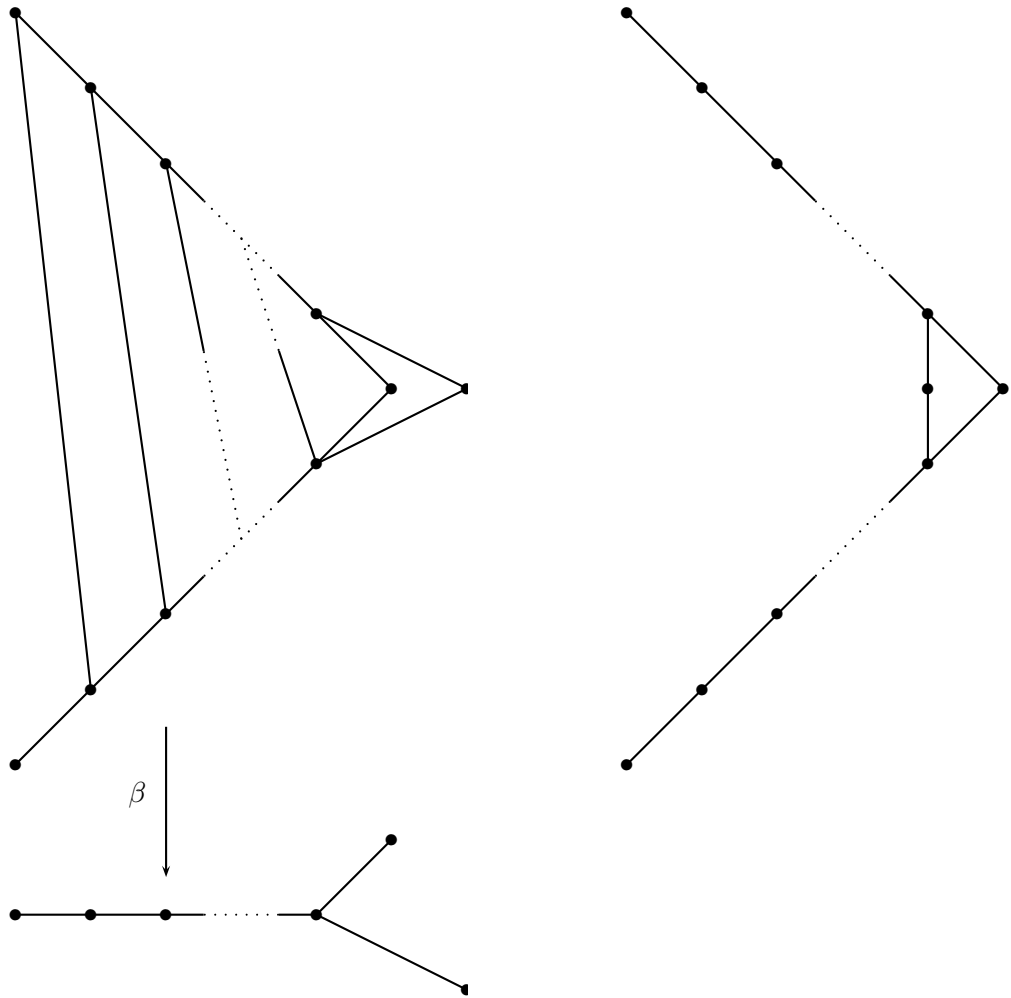
\psset{xunit=1cm}
\psset{yunit=1cm}
\centerline{%
}

Finally for $E_6$ and $E_7$ we only draw the simplified versions where
some arrows (easily detected) have been omitted and where the map $\b$
is not exactely the projection:

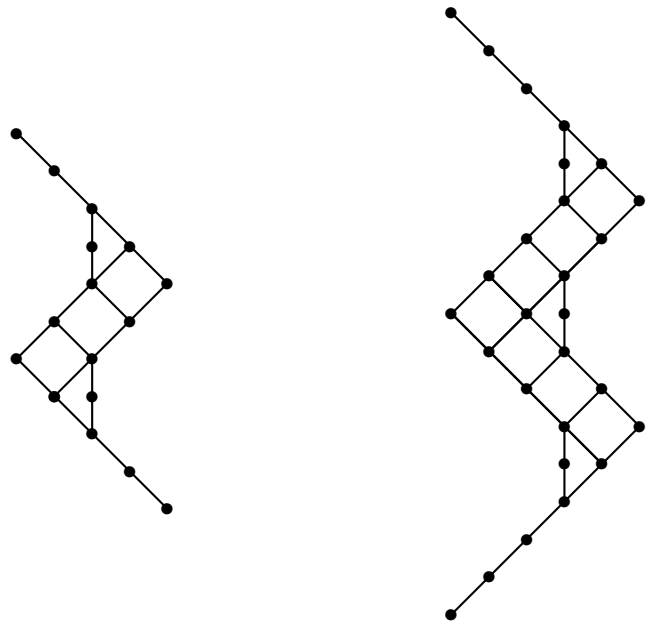
\psset{xunit=0.5cm}
\psset{yunit=0.5cm}
\centerline{%
}

\subsection{Quivers of minuscule Schubert varieties}

Thanks to the description of the quivers of minuscule homogeneous
varieties and the proposition \ref{description}, we know that the
quivers of a minuscule Schubert variety is of the following form (we
have cercled the successors of elements in the set $A$ described in
proposition \ref{description}):

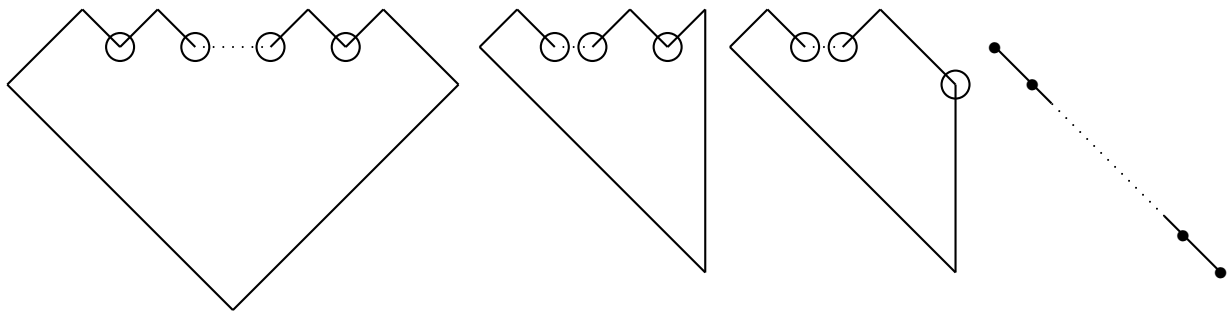
\psset{xunit=0.5cm}
\psset{yunit=0.5cm}
\centerline{%
}

or for the quadrics:

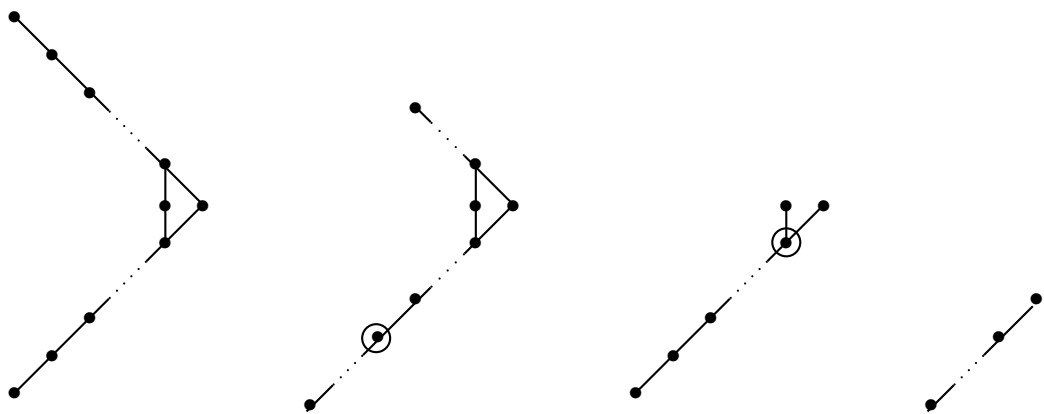
\psset{xunit=0.5cm}
\psset{yunit=0.5cm}
\centerline{%
}

and finally for $E_6$ and $E_7$ let us give a complete list of the
quivers (except the empty quiver):

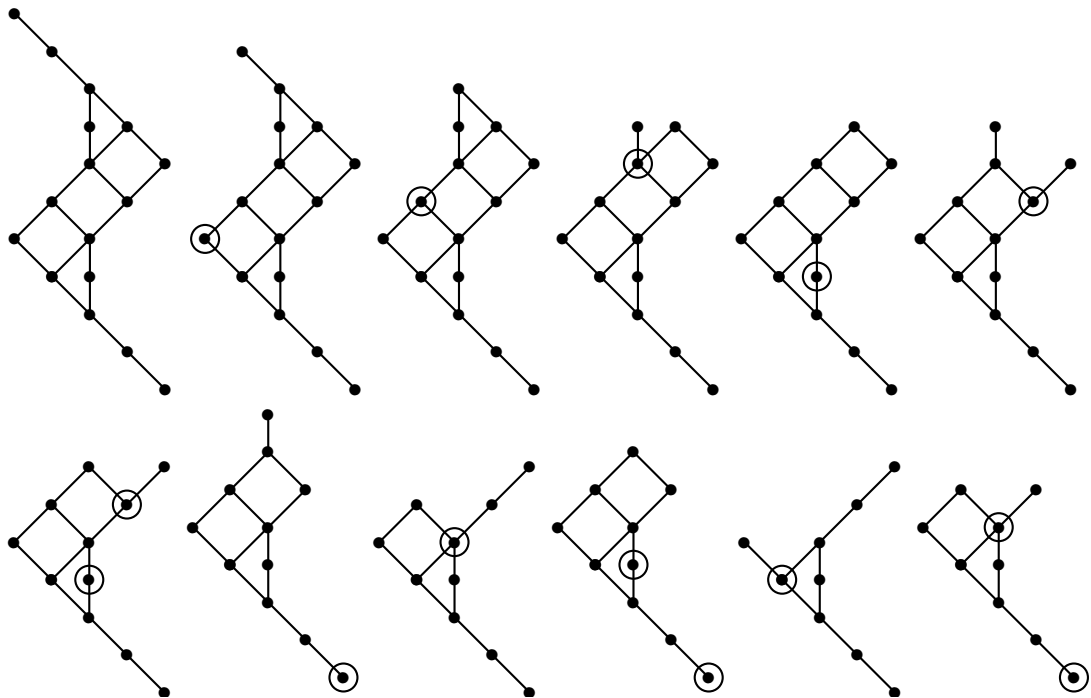
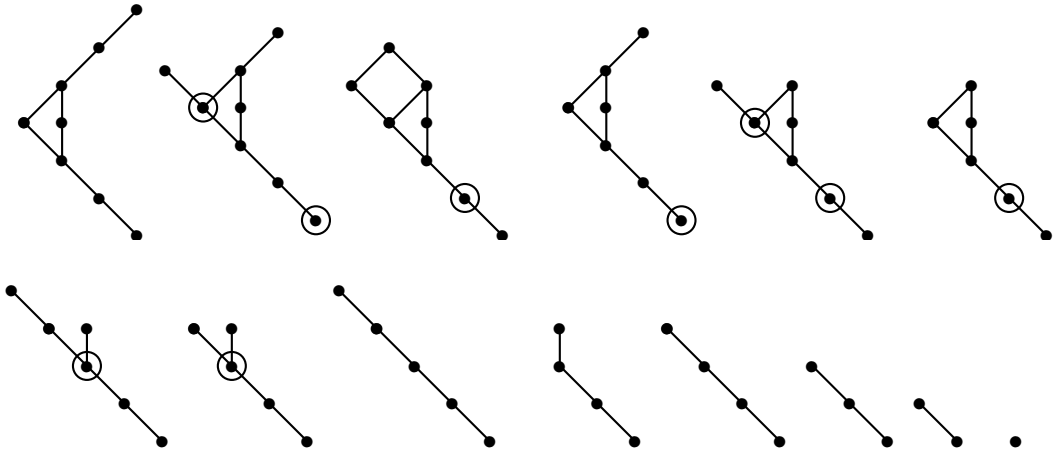
\psset{xunit=0.5cm}
\psset{yunit=0.5cm}
\centerline{%
}

It is easy to verify (thanks to our results) that the only Schubert
varieties admetting a $IH$-small resolution are the following (we only
list there number in the previous list):
1, 6, 7, 9, 11, 13, 17, 19, 20, 21, 22, 23, 24, 25, 26 and the
0-dimensional one.

Let us now list the Schubert varieties for the $E_7$ case:

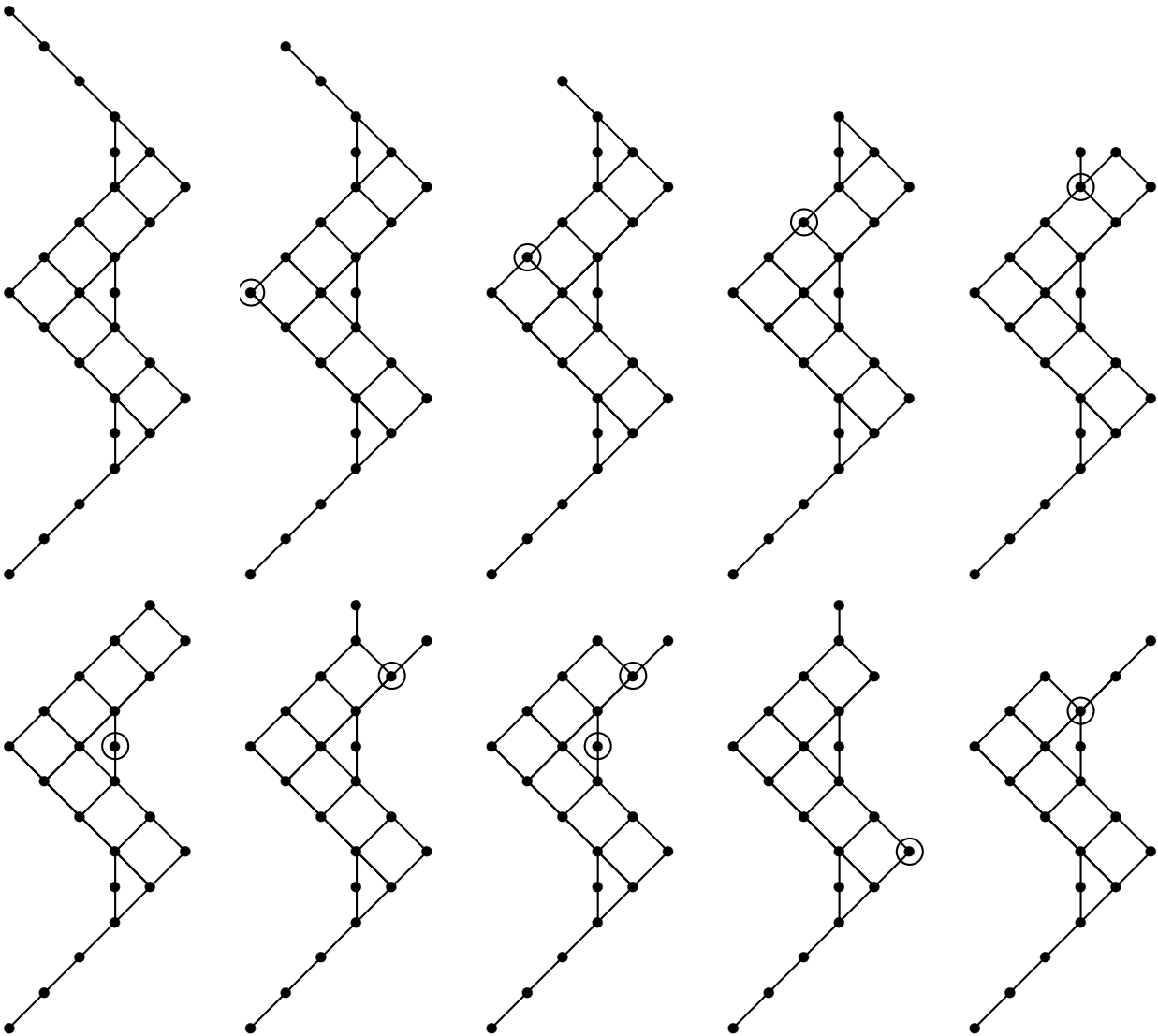
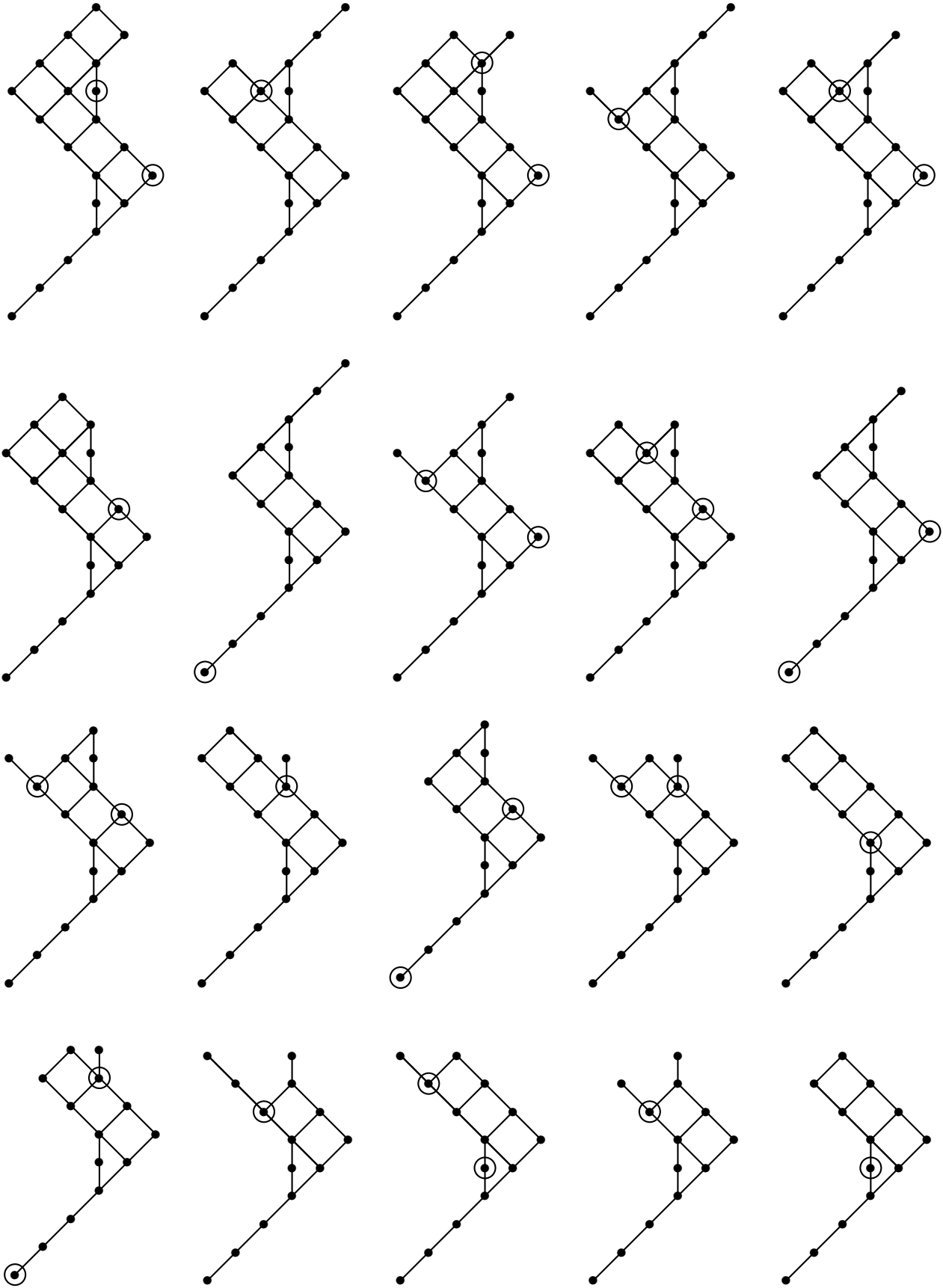
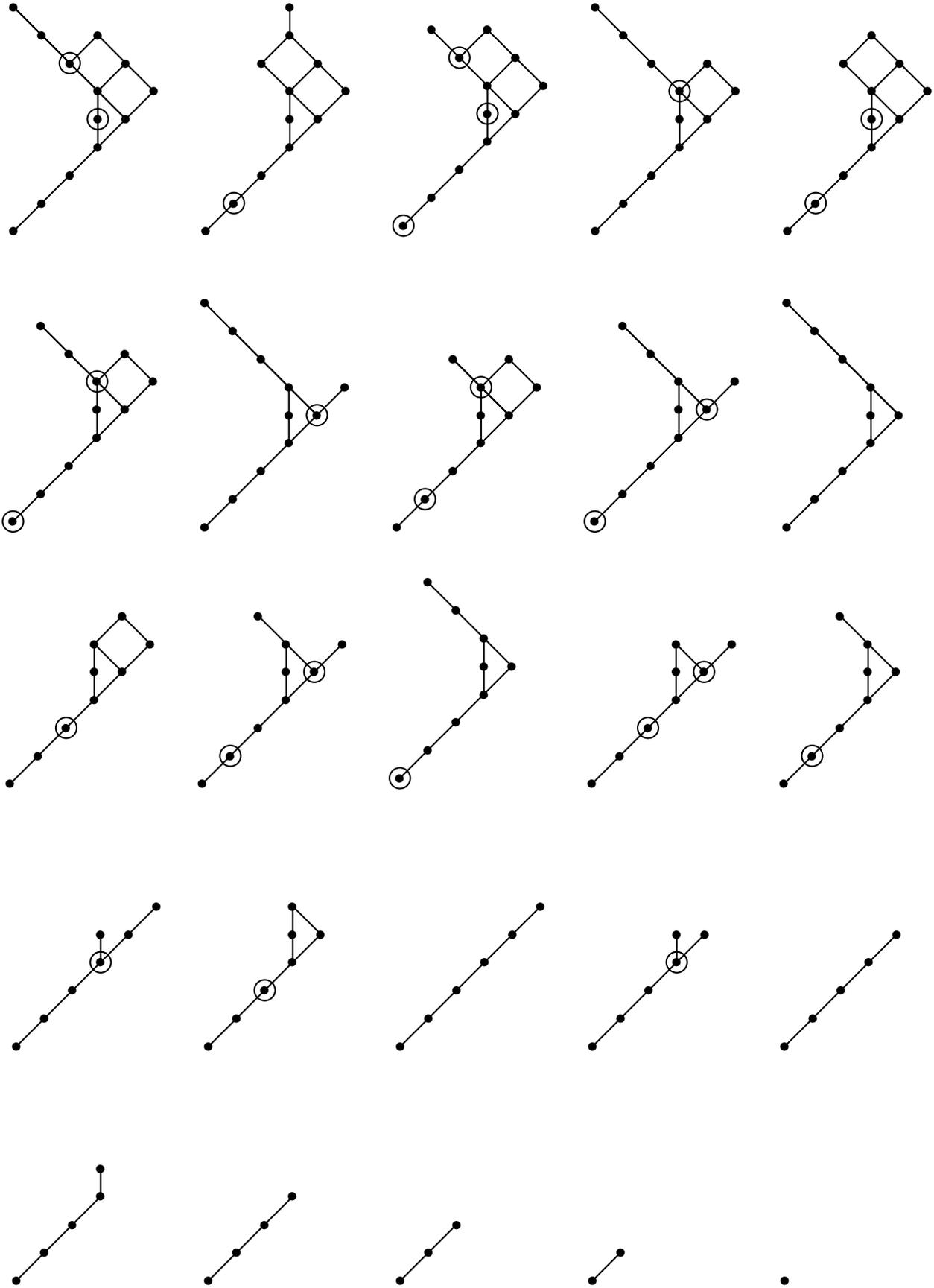
\psset{xunit=0.5cm}
\psset{yunit=0.5cm}
\centerline{%
}

It is easy to verify (thanks to our results) that the only Schubert
varieties admetting a $IH$-small resolution are the following (we only
list there number in the previous list):
1, 24, 27, 28, 31, 34, 37, 40, 44, 46, 48, 49, 50, 51, 52, 53, 54, 55 and the
0-dimensional one.

\begin{small}

\vs 0.2 cm

\noi
{\textsc{Institut de Math{\'e}matiques de Jussieu}}

\vs -0.1 cm

\noi
{\textsc{175 rue du Chevaleret}}

\vs -0.1 cm

\noi
{\textsc{75013 Paris,}} \hs 0.2 cm{\textsc{France.}}

\vs -0.1 cm

\noi
{email : \texttt{nperrin@math.jussieu.fr}}

\end{small}

\end{document}